\documentclass[12pt,leqno]{amsart}
\usepackage{amsmath,amsthm,amssymb,amscd,verbatim}
\newtheoremstyle{fancy}{}{}{\itshape}{}{\textsc\bgroup}{.\egroup}{ }{}
\newtheoremstyle{fancy2}{}{}{\rm}{}{\textsc\bgroup}{.\egroup}{ }{}
\theoremstyle{fancy}
\newcounter{intro}

\newtheorem{theo}[intro]{Theorem}

\numberwithin{equation}{section}    \swapnumbers
\newtheorem{cor}[equation]{Corollary}
\newtheorem{lem}[equation]{Lemma}
\newtheorem{prop}[equation]{Proposition}
\newtheorem{thm}[equation]{Theorem}
\newtheorem{named}[equation]{\name}
\newcommand{\name}{Proof of}
\newenvironment{nameit}[1]{\renewcommand{\name}{#1}
     \begin{named}}{\end{named}}
\theoremstyle{fancy2}
\newtheorem{dfn}[equation]{Definition}
\newcounter{axa}

\newtheorem{axia}[axa]{Axiom}
\newtheorem{exa}[equation]{Example}

\newtheorem{rem}[equation]{Remark}
\newtheorem{rems}[equation]{Remarks}
\newcommand{\cref}[1]{Corollary~\ref{#1}}
\newcommand{\lref}[1]{Lemma~\ref{#1}}
\newcommand{\pref}[1]{Proposition~\ref{#1}}
\newcommand{\tref}[1]{Theorem~\ref{#1}}
\def\C{\mathbb C}       \def\N{\mathbb{N}}
\def\R{\mathbb R}       
\def\Z{\mathbb{Z}}

 \def\la{\langle} \def\ra{\rangle}

\newcommand{\coker}{\operatorname{coker}}
\newcommand{\codim}{\operatorname{codim}}
\newcommand{\defi}{\operatorname{def}}

\newcommand{\ext}{\operatorname{ext}}

\newcommand{\ind}{\operatorname{ind}}

\newcommand{\graph}{\operatorname{graph}}
\newcommand{\im}{\operatorname{im}}
\newcommand{\Lip}{\operatorname{Lip}}
\newcommand{\loc}{\operatorname{loc}}

\newcommand{\nuli}{\operatorname{null}}
\newcommand{\rk}{\operatorname{rk}}
\renewcommand{\Re}{\operatorname{Re}}

\newcommand{\spec}{\operatorname{spec}}

\newcommand{\supp}{\operatorname{supp}}

\begin{document}


\title[Dirac-Schr\"odinger systems]
{Regularity and index theory \\ for Dirac-Schr\"odinger systems \\
with Lipschitz coefficients}
\author{Werner Ballmann, Jochen Br\"uning, \and Gilles Carron}
\address{
Mathematisches Institut\\
Rheinische Friedrich-Wilhelms-Universi\-t\"at Bonn\\
Beringstra\ss e 1, 53115 Bonn, Deutschland\\}
\email{hwbllmnn\@@math.uni-bonn.de}
\address{
Institut f\"ur Mathematik\\
Humboldt--Universit\"at\\
Rudower Chaussee 5, 12489 Berlin, Germany\\}
\email{bruening\@@mathematik.hu-berlin.de}
\address{
D\'epartement de Math\'ematiques\\
Universit\'e de Nantes\\
2 rue de la Houssini\'ere, BP 92208\\
44322 Nantes Cedex 03, France\\}
\email{Gilles.Carron\@@math.univ-nantes.fr}


\dedicatory{Dedicated to Robert Seeley on the occasion
of his 75. birthday.}

\date{\today}


\begin{abstract}
Dirac-Schr\"odinger systems play a central role 
when modeling Dirac bundles and 
Dirac-Schr\"odinger operators near the boundary, 
along ends or near other
singularities of Riemannian manifolds.
In this article we develop the Fredholm theory 
of Dirac-Schr\"odinger systems
with Lipschitz coefficients.
\end{abstract}

\maketitle


\section*{Introduction}

A {\em Dirac system} $d$ consists of a bundle
$\mathcal H\to\R_+$ of separable complex Hilbert spaces
together with a differential operator,
its {\em Dirac operator}
\begin{equation}\label{intdo}
  D = \gamma (\partial + A) , 
\end{equation}
where $\gamma=(\gamma_t)_{t\in\R_+}$ is a family of
unitary  operators on the fibers $H_t$
of $\mathcal H$ with $\gamma_t^{-1}=-\gamma_t$,
$(A_t)_{t\in\R_+}$ is a family of
self-adjoint operators on the fibers $H_t$
with discrete spectrum and anti-commuting with $\gamma$,
and $\partial$ is a metric connection on $\mathcal H$
derived from these data.
The Dirac operator is symmetric on sections 
with compact support in $(0,\infty)$.

The notion of Dirac system is traditionally connected 
with the finite dimensional versions of \eqref{intdo}
which derived from separating variables in Dirac's
original equation describing the relativistic electron.
A very influential discussion of an infinite
dimensional case was carried out in the celebrated work 
of Atiyah, Patodi, and Singer \cite{APS},
where manifolds with cylindrical ends are considered.
More generally,
Dirac systems arise in the study of Dirac operators
on Dirac bundles in the sense of Gromov-Lawson
when studying boundary value problems
or ends with special geometric features.
This is the motivation underlying the investigation 
of Dirac systems we present here.

In many situations encountered in geometry,
the data of the relevant Dirac system do not depend
smoothly on the parameter $t\in\R_+$.
For example, if $M$ is a complete, non-compact
Riemannian manifold with finite volume and
pinched negative sectional curvature,
then the Busemann functions associated to the
ends of the manifold are only $C^2$,
so that the tangent and normal bundles
of their level surfaces are only $C^1$.
This is the situation studied in \cite{BB1}
and \cite{BB2}.
The natural setup seems to be Dirac systems
with (locally) Lipschitz coefficients
as we consider them here.
The present work leads to generalizations of the
results in \cite{BB1} and \cite{BB2}.
We will discuss this in a future publication.

After \cite{APS},
where the so-called APS-projection is introduced,
it became customary to state boundary conditions
for Dirac systems in terms of orthogonal
projections in $H=H_0$.
The regularity theory of boundary conditions
defined by orthogonal projections in $H$
plays a central role in \cite{BL2},
see for example Theorem 4.3 in \cite{BL2},
an important predecessor of this article regarding 
the regularity theory of boundary conditions.

The first main contribution of the present work
consists in a new way of looking at 
boundary value problems for Dirac systems.
Let $D_0$ be the restriction of $D$ to Lipschitz 
sections of $\mathcal H$ which vanish at $t=0$.
Then $D_0$ is symmetric and contained in $D_{\max}:=D_0^*$, 
the {\em maximal extension} of $D_0$,
with domain $\mathcal D_{\max}$.
Denote by $H^s$, $s\in\R$,
the domain of the operator $(I+|A_0|^2)^{s/2}$.
For $I\subset\R$, 
denote by $Q_I$ the spectral projection of $A_0$
associated to $I\cap\spec A_0$ 
and set $H^s_I:=Q_I(H^s)$.
We show that the space
$\check H:=\{\sigma(0)\,:\,\sigma\in\mathcal D_{\max}\}$
of admissible initial values is the hybrid Sobolev space
\begin{equation}\label{intdm}
  \check H
  = H^{1/2}_{(-\infty,0]}\oplus H^{-1/2}_{(0,\infty)} . 
\end{equation}
This leads us to say that a {\em boundary value problems}
or a {\em boundary condition} for $D$ is a closed subspace
of $\check H$.
By \eqref{intdm}, the topology of the space $\check H$
is a mixture of the topologies of the spaces $H^{1/2}$
and $H^{-1/2}$ and is therefore not compatible with 
the topology of $H$ or the Sobolev spaces $H^s$, 
which causes considerable technical problems when
discussing boundary value problems given by projections.

Our first observation is that the closed 
extensions of $D_{0}$ are precisely the operators 
$D_{B,\max}$ with domain
\begin{equation}\label{intdb}
  \mathcal D_{B,\max}
  := \{ \sigma\in\mathcal D_{\max} \,:\, \sigma(0) \in B\} ,
\end{equation}
given by boundary conditions $B\subset\check H$ 
as defined above.
We show this in our discussion of constant coefficient
Dirac systems (\pref{bcbclex}), 
but the same arguments also apply in the case of Dirac
systems with Lipschitz coefficients,
cf. \tref{intreg} below.
This characterization of closed extensions of $D_0$ is 
a first confirmation that our way of defining boundary 
value problems is the appropriate one.

The adjoint operator, $D_{B,\max}^*$, arises from the
boundary form
\begin{equation}\label{intbf}
\begin{split}
  (D_{\max}\sigma_1,\sigma_2)_{L^2(\mathcal H)}
  - (\sigma_1,D_{\max}\sigma_2)_{L^2(\mathcal H)}
  &= \la\sigma_1(0),\gamma_0\sigma_2(0)\ra_H \\
  &=: \omega(\sigma_1(0),\sigma_2(0)) ,
\end{split}
\end{equation}
a non-degenerate skew-Hermitian form on $\check H$.
We show that
\begin{equation}\label{intad}
  D_{B,\max}^* = D_{B^a,\max} ,
\end{equation}
where $B^a$ denotes the annihilator of $B$ 
with respect to $\omega$.

With $H_{\loc}^1(d)$ the natural Sobolev space 
associated to $d$, we show an important regularity 
property of $\mathcal D_{\max}$,
\begin{equation}\label{intdh}
  \mathcal D_{\max} \cap H^1_{\loc}(d)
  = \{ \sigma\in\mathcal D_{\max} 
    \,:\, \sigma(0)\in H^{1/2}\} .
\end{equation}
Consequently we say that a boundary value problem $B$
for $D$ is {\em regular} if $B\subset H^{1/2}$.
We say that a boundary value problem $B$
is {\em elliptic} if $B$ and its adjoint
boundary value problem $B^a$ are both regular.
We prove next that elliptic boundary conditions
coincide with the boundary conditions introduced
in \cite{BB} (\pref{bcbsrdata}).

We say that a boundary condition $B$ is {\em self-adjoint}
if $B=B^a$. 
By definition, a self-adjoint boundary condition is
elliptic if it is regular.
In one of our main results on boundary value problems
we characterize elliptic self-adjoint boundary conditions
(\tref{bcsrsa} and \cref{bcsrsacor}).
Part of this characterization is the following result.

\begin{theo}\label{intbcsa}
Let $H^\pm := \{x\in H \,:\, i\gamma x = \pm x\}$.
Then $\check H$ contains an elliptic self-adjoint boundary
condition if and only if the restriction $A_0^+$ of $A_0$ 
to $H^+$ is a Fredholm operator to $H^-$ 
(in general unbounded) with index $\ind A_0^+ = 0$.
\end{theo}

Let $d=((H_t),(A_t),(\gamma_t))$ be a Dirac system 
with Lipschitz coefficients, and denote
by $d^0$ the Dirac system with constant coefficients
$(H_0,A_0,\gamma_0)$ and associated Dirac operator $D^0$. 
Our second main contribution to Dirac systems is the 
regularity theory for Dirac systems with Lipschitz
coefficients.
The first part of our work in this direction is concerned
with the regularity theory of the maximal domain 
(\tref{dcomp}):

\begin{theo}\label{intreg}
Let $\mathcal D_{\max}$ and $\mathcal D^0_{\max}$ be the 
domains of the maximal extension of $D$ and $D^0$,
respectively.
If $\sigma\in L^2(\mathcal H)$ has compact support, 
then $\sigma\in\mathcal D_{\max}$ if and only if
$\sigma\in\mathcal D^0_{\max}$.
\end{theo}

This result underlies the asserted equalities
in \eqref{intdm} and \eqref{intdh} above
which we show for constant coefficients first
and then extend to Lip\-schitz coefficients, 
by \tref{intreg}.

For a satisfactory analysis of the index theory of Dirac
systems it is necessary to consider extended solutions.
This goes back to the work of Atiyah, Patodi, and Singer
in \cite{APS}.
Here we rely on the approach of the third author
and his related notion of non-parabolicity,
compare \cite{Ca1} and \cite{Ca2}.
The domain of the corresponding extended Dirac
operator $D_{\ext}$ is denoted $W$,
the operator and subdomain associated to a boundary 
condition $B$ by $D_{B,\ext}$ and $W_B$, respectively.

In the second part of our work on the regularity theory
of Dirac systems we study the space of Cauchy data
of the spaces $\ker D_{\max}$ and $\ker D_{\ext}$.
Before we formulate our results in this direction,
some comments seem in order.
Let $M$ be a smooth compact manifold with boundary $N$
and $E\to M$ be a smooth Hermitian vector bundle.
Let $D:C^\infty(M,E)\to C^\infty(M,E)$
be an elliptic pseudo-differential operator of order one. 
In \cite{Calderon,Seeley}, A. Calder\'on and R. Seeley 
studied the space of Cauchy data of $\ker D$. 
Let $\mathcal{C}^s$ be the space of such data 
which belong to the Sobolev space $H^{s+1/2}(M,E)$. 
By the Trace Theorem, 
$\mathcal{C}^s$ is a subspace of $H^{s}(N,E)$.
Calder\'on and Seeley showed that there is a 
pseudo-differential projector $P$ in $H^{s}(N,E)$  
(of order $0$) onto $\mathcal{C}^s$
and that the principal symbol of $P$ is the projection 
onto the positive eigenspace of a certain operator
derived from the symbol of $D$ \footnote{
Actually, Calder\'on and Seeley considered also elliptic 
operators of higher order and treat the $L^p$ theory as well,
see Theorem 2 in \cite[p. 287]{Palais} 
or Theorem 12.4 in \cite{BW}.}. 
The projection $P$ is obtained with a single layer potential
and is not the orthogonal projection onto the $L^2$-closure 
of $\mathcal{C}^s$. 
However, B.~Boo\ss{}-Bavnbek and K.~Wojciechowski remarked
that the $L^2$-orthogonal projection has the same properties,
see Lemma 12.8 in \cite{BW}.
Our result for Dirac systems with Lipschitz coefficients 
(and its adaptation to manifolds in Chapter \ref{nsb})
is a generalization of this result to a non-smooth setting
(Theorems \ref{cald2}, \ref{cald3}, and \ref{cmaxthm});
we emphasize that this generalization is achieved without
any recourse to pseudo-differential techniques.

\begin{theo}\label{theocald}
Let $d$ be a non-parabolic Dirac system
with Lip\-schitz coefficients.
Let $\check{\mathcal C}_{\max}$ 
and $\check{\mathcal C}_{\ext}$
be the Calder\'on spaces of Cauchy data
$\sigma(0)\in\check H$ with $\sigma\in\ker D_{\max}$
and $\sigma\in\ker D_{\ext}$, respectively.
Then 
\[
  \mathcal C_{\max}^{1/2}
  := \check{\mathcal C}_{\max}\cap H^{1/2}
  \quad\text{and}\quad
  \mathcal C_{\ext}^{1/2}
  := \check{\mathcal C}_{\ext}\cap H^{1/2}
\]
are mutually adjoint elliptic boundary conditions.

Let $C_{\max}$ and $C_{\ext}$ be the orthogonal projections
in $H$ onto the closure of 
$\mathcal C_{\max}:=\check{\mathcal C}_{\max}\cap H$
and onto
$\mathcal C_{\ext}:=\check{\mathcal C}_{\ext}\cap H$, 
respectively.
Then $C_{\max}$ and $C_{\ext}$ restrict to $H^s$
and extend to $H^{-s}$, $0\le s\le1/2$, 
and 
\[
  C_{\max} - Q_{(0,\infty)}
  \quad\text{and}\quad  
  C_{\ext} - Q_{(0,\infty)}
\]
are compact in $H^s$ for all $|s|\le1/2$.
\end{theo}

Recall Kato's notion of a Fredholm pair of closed 
subspaces in a Banach space \cite[Section IV.4]{Ka}, 
compare Appendix \ref{fredpair}.
Our main index formula is formulated in terms
of such pairs (\tref{windgen}).

\begin{theo}\label{theoind}
Let $d$ be a non-parabolic Dirac system
with Lip\-schitz coefficients
and $B$ be an elliptic boundary condition.
Then $(\bar B,\mathcal C_{\ext})$ is a Fredholm pair in $H$
and
\begin{equation*}
  \ind D_{B,\ext} = \ind(\bar B,\mathcal C_{\ext}) ,
\end{equation*}
where $\bar B$ denotes the closure of $B$ in $H$.
\end{theo}

The boundary value problem considered by Atiyah, Patodi,
and Singer corresponds to $B_{APS}:=H^{1/2}_{(-\infty,0]}$.
Another main index formula is of Agranovi\v{c}-Dynin type
and shows the fundamental character of the
Atiyah-Patodi-Singer boundary condition
(\tref{relindadt}):

\begin{theo}\label{theoadt}
Let $d$ be a non-parabolic Dirac system
with Lip\-schitz coefficients
and $B$ be an elliptic boundary condition.
Then
\begin{equation*}
  \ind D_{B,\ext} 
  = \ind D_{B_{APS},\ext} + \ind(\bar B,H_{(0,\infty)}) .
\end{equation*}
\end{theo}

The Cobordism Theorem for the chiral Dirac operator $D^+$
on the space of spinor fields of a closed spin 
manifold $M$ of even dimension states that $\ind D^+=0$
if $M$ is cobordant to a compact spin manifold, 
compare \cite[Corollary 21.6]{BW}.
We prove a version of the Cobordism Theorem
for Dirac systems with Lipschitz coefficients (\tref{cothm}).
As in \tref{intbcsa} above, 
let $H^\pm := \{x\in H \,:\, i\gamma x = \pm x\}$
and $A_0^+$ be the restriction of $A_0$ to $H^+$,
a Fredholm operator to $H^-$.

\begin{theo}[Cobordism Theorem]\label{intcoth}
Let $d$ be a Dirac system with Lipschitz coefficients.
If the associated Dirac operator $D$ 
is of Fredholm type in the sense that $d$ is non-parabolic 
with $W=\mathcal D_{\max}$, then \[\ind A_0^+ = 0 . \]
\end{theo}

When cutting a manifold $M$ into pieces $M_1$ and $M_2$
along a compact hypersurface $N=M_1\cap M_2$,
we may ask for the index of a Dirac operator $D$
on sections of a Hermitian bundle $E$ over $M$ 
in terms of its restrictions to the pieces.
The corresponding boundary condition along $N$,
the so-called {\em transmission boundary condition}, 
requires that sections $\sigma_1$ and $\sigma_2$ of $E$
over $M_1$ and $M_2$, respectively, coincide along $N$.
In terms of Dirac systems, the decomposition
of $M$ and $D$ corresponds to the direct sum
of two Dirac systems which have compatible
initial conditions at $t=0$.
Our first result concerning this type of boundary value
problem is of Bojarski type (\tref{botythm}):

\begin{theo}\label{intbotythm}
Let $d_1$ and $d_2$ be non-parabolic Dirac systems
with Lipschitz coefficients and Calder\'on spaces
$\mathcal C_{1,\ext}$ and $\mathcal C_{2,\ext}$,
respectively.
Suppose that the initial conditions of $d_1$ and $d_2$
satisfy 
\[
  H:= H_{1,0} = H_{2,0} , \quad
  A:= A_{1,0} = -A_{2,0} , \quad\text{and}\quad 
  \gamma_{1,0} = -\gamma_{2,0} .
\]
Then $(\mathcal C_{1,\ext},\mathcal C_{2,\ext})$
is a Fredholm pair in $H$.

Consider the Dirac operator $D$ on $d=d_1\oplus d_2$
with transmission boundary condition
$B=\{(x,x)\,:\, x\in H^{1/2}\}$.
Then $B$ is an elliptic and self-adjoint boundary
condition and
\[
  \ind D_{B,\ext} 
  = \ind (\mathcal C_{1,\ext},\mathcal C_{2,\ext}) .
\]
\end{theo}

Another convenient way of determining the index
of a Dirac operator via decompositions is by
decoupling the boundary conditions on the pieces
$M_1$ and $M_2$.
Our relevant result in this direction (\tref{decothm})
generalizes Theorem 4.3 of \cite{BL1}.

\begin{theo}\label{intdecothm}
Let $d_1$ and $d_2$ be non-parabolic Dirac systems
with Lipschitz coefficients as in \tref{intbotythm}.
Then
\[
  \ind D_{B,\ext}
  = \ind D_{1,B_1,\ext} + \ind D_{2,B_2,\ext} ,
\]
where $B$ is the transmission boundary condition,
$B_1$ is any elliptic boundary condition for $d_1$, 
and $B_2=B_1^\perp\cap H^{1/2}$.
\end{theo}

The above results are discussed and proved
in Chapters 1--3 of the text.
Many of our arguments and results here
extend and simplify what is known from the literature.
In Chapter 4, 
we discuss supersymmetric Dirac systems
and derive the corresponding index formulas.
In Chapter 5, we describe a geometric setup
for non-smooth boundary value problems
for differential operators of Dirac type
and explain how our results extend
to this situation.
This will be important for our geometric applications
in a forthcoming article, 
in which we will extend the results from \cite{BB1,BB2}.
We believe that it will also be useful in further work
on boundary value problems and index theory
of Dirac type operators.
We derive our results not only for Dirac operators,
but for the more general class of Dirac-Schr\"odinger
operators, that is, operators of the form $D+V$,
where $D$ is a Dirac operator and $V$ is a symmetric
potential, see Definition \ref{dsdef}.

In two appendices, 
we derive some results which are used in the main text
and seem to be of independent interest,
but are not closely connected with the program
we are pursuing here.

In all our estimates,
generic constants may change from line to line.

WB and JB would like to use this
occasion to refer to the article \cite{Kas},
which already contains one of the main observations
underlying the proof of Theorem B of \cite{BB2}
and also similar applications.
We would like to thank Tobias Ebel
for pointing this out to us.

WB, JB, and GC would like to thank the MPI f\"ur
Mathematik in Bonn for its hospitality.
WB and JB enjoyed the hospitality of the MSRI in Berkeley
and the FIM at the ETH in Z\"urich.
WB appreciated helpful discussions with Charles Epstein
and is grateful to the ESI in Vienna for its hospitality.
JB wants to thank Bob Seeley and Jean-Michel Bismut
for helpful conversations,
and he is indebted to the Universit\'e Paris-Sud
and the University of Bergen for their hospitality;
he acknowledges the financial support of the SFB 647
gratefully.

\tableofcontents

\newpage
\section{Dirac systems with constant coefficients}\label{sec:cc}
\subsection{Generalities}\label{subsecgen}
Let $H$ be a separable complex Hilbert space
with Hermitian inner product
$\langle\cdot,\cdot\rangle_{H}=\langle\cdot,\cdot\rangle$
and norm $||\cdot||_{H}=||\cdot||$.
Let $A$ be a self-adjoint operator in $H$
with domain $H_{A}$ such that,
with respect to the graph norm $||\cdot||_{A}$,
the embedding $H_{A}\to H$ is compact;
equivalently, $A$ is discrete in the sense that $\spec A$
consists only of isolated eigenvalues with finite multiplicity.
The pair $e:=(H,A)$ will be referred to as
an {\em evolution system}
since we will associate an evolution operator to it.
To that end we note first that any local Lipschitz function
$\sigma: \R_{+} := [0,\infty)\to H$ is weakly differentiable
with locally uniformly bounded weak derivative $\sigma'$ a.e.;
this is a well known fact,
but for the sake of completeness we will give a proof below.
Then we can introduce the space
\begin{equation}\label{eq:evspace}
  \mathcal{L}_{\loc}(e)
  := \Lip_{\loc}(\R_{+},H)\cap L^\infty_{\loc}(\R_{+},H_{A})
\end{equation}
and the operator
\begin{equation}
     L = L(e) :=\partial_{t} + A:
     \mathcal{L}_{\rm loc}(e)\to
     L^{\infty}_{\rm loc}(\R_{+},H),
     \label{eq:evoperator}
\end{equation}
where $\partial_{t}\sigma=\sigma'$ denotes the derivative
of $\sigma$ with respect to $t$.
We call $L$ the {\em evolution operator} associated to
the evolution system $e = (H,A)$.

\begin{lem}\label{lem:weakdiff}
If $f: \R_{+} \to H$ is locally Lipschitz,
then $f$ is weakly differentiable almost everywhere
with locally uniformly bounded derivative.
More precisely,
if $L_{[a,b]}(f)$ denotes the Lipschitz constant of $f$ on $[a,b]$,
then
\[
     ||f'(t)||_{H} \le L_{[a,b]}(f) ,
\]
for almost all $t\in [a,b]$.
\end{lem}

\begin{proof}
Since $H$ is separable,
there is a countable orthonormal basis $(e_{n})_{n\in \N}$ of $H$.
By Lebesgue's Theorem, there exists a  measurable subset
$\mathcal{R}\subset\R_{+}$ of full measure such that the functions
$t\mapsto\langle f(t), e_{n}\rangle$ are differentiable
in all points of $\mathcal{R}$ for all $n\in \N$.
Hence the functions $t\mapsto\langle f(t), u\rangle$,
where $u$ is in the dense subspace of $H$ generated
by the chosen basis,
are also differentiable in all points of $\mathcal{R}$.
We have
\[
  |\langle h^{-1}(f(t+h)- f(t)), u\rangle|
  \le L_{[0,T]}(f)||u|| ,
\]
for all $u\in H$ and $t,h\in\R_{+}$ with $t,t+h\in[0,T]$.
It follows that $t\mapsto\langle f(t), u\rangle$ is
differentiable in $\mathcal{R}$ for {\em all} $u\in H$
and thus that the function $f$ has a weak derivative,
$f'(t)\in H$, in each $t\in\mathcal{R}$
and with the asserted norm estimate.
\end{proof}

We will also need the spaces
\begin{align}
     \mathcal{L}_{c}(e)
     &:= \{\sigma\in\mathcal{L}_{\rm loc}(e) \,:\,
     \text{$\supp\sigma$ compact}\},
     \label{eq:evspacec}  \\
     \mathcal{L}_{0,c}(e)
     &:= \{\sigma\in\mathcal{L}_{c}(e) \,:\, \sigma(0) = 0\} .
     \label{eq:evspaceco}
\end{align}
On $\mathcal{L}_{c}(e)$, we define the scalar product
\begin{equation}
  (\sigma_{1}, \sigma_{2})  := \int_{0}^{\infty}
     \langle \sigma_{1}(t), \sigma_{2}(t)\rangle\, dt ,
     \label{eq:scalprod}
\end{equation}
and we denote by $L^{2}(\R_{+}, H)$ the Hilbert space
arising by completion.

The formal adjoint of $L$ in $L^{2}(\R_{+}, H)$
is $-\partial_t+A$,
hence $L$ does not induce a symmetric operator
on $\mathcal{L}_{0,c}(e)$.
This defect can be cured if there is an operator
$\gamma\in\mathcal{L}(H)\cap\mathcal{L}(H_{A})$
which satisfies the following two relations:
\begin{alignat}{2}
    -\gamma = \gamma^* &= \gamma^{-1} \quad
    &&\mbox{on $H$} , \label{cc:gamma1} \\
    A\gamma + \gamma A &= 0
    &&\mbox{on $H_{A}$} . \label{cc:gamma2}
\end{alignat}
Note that \eqref{cc:gamma2} implies that $\spec A$
is symmetric with respect to $0$.
Then the triple $d:=(H,A,\gamma)$ is called a {\em Dirac system}.
The associated {\em Dirac operator} is defined as
\begin{equation}\label{cc:diracsys}
    D = D(d) :=\gamma (\partial_{t} + A):
       \mathcal{L}_{\rm loc}(e)
       \to L^{\infty}_{\rm loc}(\R_{+}, H) .
\end{equation}
We find, for $\sigma_{1},\sigma_{2}\in\mathcal{L}_{c}(e)$,
\begin{equation}
     \langle \gamma\sigma_{1},\sigma_{2}\rangle'
     =\langle D\sigma_{1}, \sigma_{2}\rangle
     - \langle \sigma_{1}, D\sigma_{2}\rangle,
     \label{eq:intpart0}
\end{equation}
hence
\begin{equation}
     (D\sigma_{1}, \sigma_{2}) - (\sigma_{1}, D\sigma_{2})
     = \langle\sigma_{1}(0), \gamma\sigma_{2}(0)\rangle
     =: \omega(\sigma_{1}(0), \sigma_{2}(0)) ,
     \label{eq:intpart1}
\end{equation}
and therefore the restriction $D_{0,c}$ of $D$
to $\mathcal{L}_{0,c}(e)$ is symmetric.
The adjoint operator $D_{\max}:=(D_{0,c})^*$ of $D_{0,c}$ is
called the {\em maximal extension} of $D_{0,c}$;
its domain is denoted by $\mathcal D_{\max}$.
The closure $D_{\min}$ of  $D_{0,c}$ is called
the {\em minimal extension} of  $D_{0,c}$,
the domain of  $D_{\min}$ is denoted by $\mathcal D_{\min}$.
By definition,
\begin{equation}\label{maxop0}
  D_{0,c} \subset D_{\min} = (D_{\max})^* \subset D_{\max} .
\end{equation}
For later purposes it is useful to note that norm estimates
for $L\sigma$ also hold for $D\sigma$,
\begin{equation}\label{cc:normest}
   ||D\sigma(t)||_{H} = ||L\sigma(t)||_{H}
\end{equation}
for all $\sigma\in\mathcal{L}_{\rm loc}(e)$ and $t\in\R_{+}$.

We denote by $H^{1}(e)$ the closure of $\mathcal{L}_{c}(e)$
under the norm
\begin{equation}\label{cc:sobar}
   \| \sigma \|_{H^{1}(e)}^2
   := \| \sigma \|_{L^2(\R_+,H)}^2
      + \| \partial_t\sigma \|_{L^2(\R_+,H)}^2
      + \| A\sigma \|_{L^2(\R_+,H)}^2,
\end{equation}
which is naturally associated to the data
defining the evolution system\footnote{
The notation $H_{1}(\R_+,A)$ is also common
and was used e.g. in \cite{BL2}.}.
We will also use the space
\begin{multline}\label{cc:sobarloc}
  H^1_{\loc}(e) := \\
  \{ \text{$\sigma:\R_{+}\to H$ measurable} \,:\,
  \text{$\psi\sigma\in H^{1}(e)$
  for all $\psi\in\Lip_{c}(\R_{+})$} \} .
\end{multline}
Note that the norm of $H^{1}(e)$ is stronger
than the graph norm of $D$.
In particular, we have a continuous extension
\begin{equation}\label{donh1}
  D: H^{1}(e)\to L^2(\R_+,H) .
\end{equation}
Moreover, if the domain $\bar{\mathcal{D}}$
of some closed extension $\bar{D}$ of $D_{0,c}$
is contained in $H^{1}(e)$,
then the $H^1(e)$-norm and the graph norm of $\bar D$
are equivalent on $\bar{\mathcal{D}}$,
by the Closed Graph Theorem.
This fact will be used repeatedly.

Spectral projections of $A$ will play a specific role;
we reserve the letter $Q$ for these.
For a subset $J\subset\R$, $Q_J = Q_J^{*}$ denotes
the associated spectral projection of $A$ in $H$.
As shorthand, we use, for $\Lambda\in\R$,
\begin{equation}\label{cc:q}
\begin{split}
   Q_{>\Lambda} &:= Q_{(\Lambda,\infty)} ,\quad
   Q_{\ge\Lambda} := Q_{[\Lambda,\infty)} ,\\
   Q_{<\Lambda} &:= Q_{(-\infty,\Lambda)} ,\quad
   Q_{\le\Lambda} := Q_{(-\infty,\Lambda]} ,\\
\end{split}
\end{equation}
We also use $Q_0:=Q_{\{0\}}$ and
\begin{equation}\label{cc:q0}
\begin{split}
   Q_> &:= Q_{>0} ,\quad
   Q_\ge := Q_{\ge0} ,\quad
   Q_< := Q_{<0} ,\quad
   Q_\le := Q_{\le0} ,\\
   Q_{\ne} &= Q_< + Q_> = I-Q_0 .
\end{split}
\end{equation}
Since $\gamma$ anticommutes with $A$,
we have $\gamma^{*}Q_{J}\gamma = Q_{-J}$.
In particular,
\begin{equation}\label{cc:q2}
   \gamma^* Q_>\gamma = Q_< ,
   \quad
   \gamma^* Q_\le\gamma = Q_\ge ,
   \quad\mbox{and}\quad
   \gamma^* Q_0\gamma = Q_0 .
\end{equation}
Furthermore, we set $H_<:=Q_<(H)$ and use a similar notation
in the other cases
and for the Sobolev spaces associated to $A$ below.

\subsection{Sobolev spaces associated to $A$}\label{soba}
Let $H$ and $A$ be as above.
For $s\ge0$, let $H^s\subset H$ be the domain of $|A|^s$.
Then $H=H^0$ and $H_A=H^1$.
We set $H^\infty=\cap_{s\ge0}H^s$,
which is a dense subspace of $H$.

For $s\in\R$, we define a scalar product
$\la\cdot,\cdot\ra_s$ on $H^\infty$,
\begin{equation}\label{cc:sobscale}
   \la x,y \ra_s
   := \la (I+A^2)^{s/2}x,(I+A^2)^{s/2}y \ra .
\end{equation}
For $s\ge0$, the norm $||\cdot||_s$ 
associated to $\la\cdot,\cdot\ra_s$ 
is equivalent to the graph norm of $|A|^s$
and $H^s$ is equal to the completion of $H^\infty$
with respect to $||\cdot||_s$.
For $s<0$, we define $H^s$ to be the completion
of $H^\infty$ with respect to $||\cdot||_s$
and set $H^{-\infty}:=\cup_{s\in\R} H^s$.
The pairing
\begin{equation}\label{cc:sobpairing}
   B_s: H^s \times H^{-s} \to \C , \quad
   B_s(x,y) := \la (I+A^2)^{s/2}x,(I+A^2)^{-s/2}y \ra ,
\end{equation}
is perfect, that is,
it identifies $H^{-s}$ with the dual space of $H^s$.
In particular, any $S\in\mathcal L(H^s)$ admits a
dual operator $S'\in\mathcal L(H^{-s})$ with
\begin{equation}\label{cc:sobdual}
   B_s(Sx,y) = B_s(x,S'y) .
\end{equation}
This defines an algebra antimorphism
$\mathcal L(H^{s})\to\mathcal L(H^{-s})$.
More generally,
for $S\in\mathcal{L}(H^{s_{1}}, H^{s_{2}})$,
we obtain a dual operator
$S'\in\mathcal{L}(H^{-s_{2}}, H^{-s_{1}})$;
in particular, if $s = s_{1} = -s_{2}$,
then $S, S'\in\mathcal{L}(H^{s}, H^{-s})$.

Since $A$ is discrete,
the embedding $i_{t,s}: H^t\hookrightarrow H^s$ is compact
for $s<t$.
For $0\le\theta\le1$ and $r=\theta s+(1-\theta)t$,
$H^r$ is (isomorphic to) the interpolation space
$[H^s,H^t]_\theta$, see for example \cite[Chapter 4.2]{T}.
If $S\in\mathcal L(H^s)$ satisfies $S(H^t)\subset H^t$,
then $S:H^t\to H^t$ is continuous,
by the Closed Graph Theorem.
Moreover, $S(H^r)\subset H^r$ for any $r$ as above,
by interpolation.
Note also that $(i_{t,s})' = i_{-s,-t}$.

We say that an operator $S\in\mathcal L(H)$
is {\em $s$-smooth}, for $s\ge 0$,
if both $S$ and $S^{*}$ restrict to $H^{s}$;
this implies that $S,S^{*}$ restrict to $H^{t}$ and extend
(continuously) to $H^{-t}$ for $0\le t\le s$.
In fact, the dual of the restriction of $S$ and $S^*$ to $H^s$
extends $S^*$ and $S$ to $H^{-s}$, respectively.

An $s$-smooth operator $S$ is said to be {\em $(-s,s)$-smoothing}
if $S$ maps $H^{-s}$ into $H^{s}$;
if $S$ is $(-s,s)$-smoothing, then so is $S^*$.
In the special case $s = 1/2$ we simply speak of
{\em smoothing operators}.
Note that $S\in\mathcal{L}(H)$ is smoothing if $S$
extends to $H^{-1/2}$ with image in $H^{1/2}$.

We say that an operator $S\in\mathcal L(H)$ has {\em order} $0$,
if $S$ and $S^*$ restrict to $H^s$ for any $s>0$;
that is, $S$ is of order $0$ if and only if $S$
is $s$-smooth for all $s\ge 0$.
The space of operators of order $0$ is denoted $\mbox{Op}^0(A)$.
By definition,
all spectral projections of $A$ have order $0$.

We are primarily interested in the cases $s=-1/2,0,1/2$
and $s=1$.
If $S\in\mathcal L(H)$ extends continuously to $H^{-1/2}$,
then the extension is denoted by $\tilde S$,
\begin{equation}\label{cc:tilde}
   \tilde S: H^{-1/2} \to H^{-1/2};
\end{equation}
if $S\in\mathcal L(H)$ restricts to $H^{1/2}$,
then the restriction is denoted by $\hat S$,
\begin{equation}\label{cc:hat}
   \hat S: H^{1/2} \to H^{1/2}.
\end{equation}
If there is no danger of confusion,
we also write $S$ instead of $\hat S$ or $\tilde S$.

If the adjoint operator $S^*$ of $S\in\mathcal L(H)$
restricts to $H^{1/2}$,
then $S$ extends continuously to $H^{-1/2}$,
\begin{equation}\label{cc:hattil}
  \tilde S = (\widehat{S^*})' .
\end{equation}
In particular, if $S=S^*$ and $S(H^{1/2})\subset H^{1/2}$,
then $\tilde S=\hat S'$.
If $Q$ is a spectral projection of $A$,
then $Q(H^s)\subset H^s$ for any $s\in\R$,
by the definition of $H^s$.
Since $Q^*=Q$,
we have $\tilde Q=\hat Q'$ for any such $Q$.

The following lemma and corollary will only be used
in the discussion of regular pairs of projections
in Section \ref{pop}.

\begin{lem}\label{lem:preim}
Let $S\in\mathcal{L}(H)$ be $1/2$-smooth.
Then the following conditions are equivalent:
\begin{enumerate}
\item\label{eq:smprim}
Let $x\in H^{-1/2}$.
If $\tilde Sx\in H^{1/2}$ or $\widetilde{S^*}x\in H^{1/2}$, 
then $x\in H^{1/2}$.
\item\label{eq:sapriori}
$\ker\tilde S = \ker\hat S$,
$\ker\widetilde{S^{*}} = \ker\widehat{S^{*}}$,
and there is a constant $C$ with
\begin{equation*}
\begin{matrix}
   ||x||_{1/2}\le C(||\hat Sx||_{1/2} + ||x||_{-1/2}) \\
   ||x||_{1/2}\le C(||\widehat{S^*}x||_{1/2} + ||x||_{-1/2})
\end{matrix}   \quad \text{for all $x\in H^{1/2}$} .
\end{equation*}
\item\label{eq:sind}
$\hat S$ and $\widehat{S^{*}}$ are Fredholm operators
with $\ind \hat S   + \ind\widehat{S^{*}}  = 0$.
\item\label{eq:lrparam}
There are a $1/2$-smooth operator $U$ and smoothing operators
$K_{r}, K_{l}$ in $\mathcal{L}(H)$ such that
\begin{equation*}
    \tilde S\tilde U = \widetilde{U^{*}}\widetilde{S^{*}} =
    I - \widetilde{K_{l}}
    \quad\text{and}\quad
    \tilde U\tilde S = \widetilde{S^{*}}\widetilde{U^{*}} =
    I - \widetilde{K_{r}} .
\end{equation*}
\end{enumerate}
\end{lem}

\begin{proof}
$(1)\Rightarrow (2)$.
The assertion on the kernels is an obvious consequence of
\eqref{eq:smprim}.
Consider next $\widehat{S^{(*)}}$ as an unbounded operator
in $H^{-1/2}$ with domain and target space $H^{1/2}$.
Then it follows from \eqref{eq:smprim}
that $\widehat{S^{(*)}}$ is closed.
The projection $\pi_{1}: H^{-1/2}\times H^{1/2}\to H^{-1/2}$
takes values in  $H^{1/2}$
when restricted to the graph of $\widehat{S^{(*)}}$.
Applying the Closed Graph Theorem to this map
we derive the asserted inequalities in \eqref{eq:sapriori}.

\noindent $(2)\Rightarrow (3)$.
By \lref{trfred} in Appendix \ref{fredpair},
the a priori estimate in \eqref{eq:sapriori} implies
that $\hat S$ and $\widehat{S^{*}}$ have finite-dimensional
kernels and closed images in $H^{1/2}$.
 From the assumption on the kernels and duality
we deduce that
\begin{align*}
    \codim\hat S &= \dim\, (\im\hat {S})^{0}
    = \dim\ker\widetilde{S^{*}} = \dim\ker\hat{S^{*}}, \\
    \codim\widehat{S^{*}} &= \dim\,\ker\hat S ,
\end{align*}
where the superscript $0$ indicates 
the annihilator (or polar set) in $H^{-1/2}$.
This establishes \eqref{eq:sind}.

\noindent $(3)\Rightarrow (4)$.
It is immediate from the assumptions that
$\ker\tilde S = \ker\hat S$ and
$\ker\widetilde{S^{*}} = \ker\widehat{S^{*}}$.
Choose a basis $(e_{j}^{(*)})\subset H^{1/2}$ of
$\ker\widehat{S^{(*)}}$ which is orthonormal in $H$ 
and set
\[
  K_{r(l)} x 
  := \sum B_{-1/2}(x,e_{j}^{(*)})\, e_{j}^{(*)} ,
  \quad x\in H^{-1/2} .
\]
Then $K_{r(l)}\in\mathcal{L}(H^{-1/2}, H^{1/2})$
is a projection in $H^{-1/2}$ onto $\ker\widehat{S^{(*)}}$
and $\tilde S:\ker K_r\to\ker K_l$ is an isomorphism.
It follows that there is a $1/2$-smooth operator
$U\in\mathcal{L}(H)$ with
\[
     \tilde S\tilde U  = I - \widetilde{K_{l}}
     \quad\text{and}\quad
     \tilde U\tilde S  = I - \widetilde{K_{r}}.
\]
Restricting to $H^{1/2}$ and computing the dual operators
gives the remaining identities in \eqref{eq:lrparam}.

\noindent $(4)\Rightarrow (1)$.
Consider $x\in H^{-1/2}$ with $y :=\tilde{S}x\in H^{1/2}$.
Then we obtain from \eqref{eq:lrparam}
\[
     x = \hat{U}y + K_{r}x \in H^{1/2},
\]
since $K_{r}$ is smoothing;
a similar argument works for $\widetilde{S^{*}}$.
\end{proof}

Since, by complex interpolation, both $\tilde S$ and $\tilde U$
restrict to $H^{s}$, for $|s|\le 1/2$, we have the following
consequence.

\begin{cor}\label{cor:sfred}
Under the conditions of \lref{lem:preim}, $\widetilde{S^{(*)}}$
restricts respectively extends to a Fredholm operator on each
$H^{s},\,|s|\le 1/2$, with kernel and index independent of $s$.
\end{cor}

\subsection{The domain of the maximal extension}
In our approach, boundary conditions at $0$
will play a prominent role; 
for that reason,
the existence of restriction maps is of interest.
We begin with the following regularity lemma;
its third part reflects the usual trace properties
of Sobolev spaces.

\begin{lem}[Regularity I]\label{lem:h1restr}
We have 
\begin{enumerate}
\item
$\mathcal L_{\loc}(e)\subset C^{0+1/2}(\R_+,H^{1/2})$.
\item
$H^1(e)\subset C(\R_+,H^{1/2})$.
\item
  $\mathcal{R}: H^{1}(e)\to H^{1/2}$, 
  $\mathcal{R}\sigma := \sigma(0)$, is continuous.
\end{enumerate}
\end{lem}

\begin{proof}
By the Cauchy-Schwarz inequality, we have, 
for any $x\in H_A$,
\[
  ||x||_{H^{1/2}}^2
  \le ||x||_{H_A} ||x||_H .
\]
Hence if $\sigma\in\mathcal L_{\loc}(e)$
with $||\sigma||_{H_A}\le K$ and $||\sigma'||_H\le L$
on $[0,T]$,
then 
\begin{align*}
  ||\sigma(s)-\sigma(t)||_{H^{1/2}}^2
  &\le ||\sigma(s)-\sigma(t)||_{H_A}
      ||\sigma(s)-\sigma(t)||_H \\
  &\le 2KL |s-t| ,   
\end{align*}
for all $s,t\in[0,T]$.
This shows the first claim.

As for the proof of the second and third claim,
we choose an orthonormal basis, $(e_{n})_{n\in\N}$, for $H$,
consisting of eigenvectors of $A$, i.e. $Ae_{n} = a_{n}e_{n}$ for
some $a_{n}\in\R$.
For $\sigma\in\mathcal{L}_{c}(e)$
we set $\sigma_{n}(t) := \langle \sigma(t),e_{n}\rangle$.
Then $\sigma_{n}\in\Lip_{c}(\R_{+})$ and hence,
by \eqref{ineq6},
\[
  |a_{n}||\sigma_{n}(t)-\sigma_n(s)|^{2}
  \le 2||\sigma_n'||_{L^2([s,t])}^2 
      + 2 a_n^2 ||\sigma_n||_{L^2([s,t])}^2 .
\]
Therefore
\begin{multline}\label{h1diff}
  ||\sigma(t)-\sigma(s)||_{H^{1/2}}^2 \\
  \le C  \big( ||\sigma||_{L^2([s,t],H)}^2 
      +||\sigma'||_{L^2([s,t],H)}^2 
      + ||A\sigma||_{L^2([s,t],H)}^2 \big) .
\end{multline}
Since $\sigma$ has compact support,
\begin{multline}\label{h1restr1}
   ||\sigma(s)||_{H^{1/2}} \\
   \le C  \big( ||\sigma||_{L^2([s,\infty],H)}^2 
      +||\sigma'||_{L^2([s,\infty],H)}^2 
      + ||A\sigma||_{L^2([s,\infty],H)}^2 \big) .
\end{multline}
In particular, 
\[
  ||\sigma(0)||_{H^{1/2}}\le C||\sigma||_{H^{1}(e)} .
\]
Since $H^1(e)$ is the closure of $\mathcal{L}_{c}(e)$
in the $H^1(e)$-norm, \eqref{h1diff} and \eqref{h1restr1}
hold for all $\sigma\in H^1(e)$.
Claims (2) and (3) follow.
\end{proof}

To get a satisfactory description of the domain
$\mathcal D_{\max}\subset L^2(\R,H)$ of the maximal
extension $D_{\max}$ of $D_{0,c}$,
we employ the solution theory of the evolution operator $L$.
For $\sigma\in L^2(\R_+,H)$ we set
\begin{equation}\label{cc:sobsa}
   (S_L\sigma) (t)
   := \int_0^t e^{(s-t)A_>} \sigma(s) ds
   - \int_t^\infty e^{(s-t)A_<}\sigma(s) ds ,
\end{equation}
where we have written $A_>:=AQ_>$ and $A_<:=AQ_<$.
The {\em solution operator} $S_L$ has been studied
in \cite[Proposition (2.5)]{APS} via the corresponding
ordinary differential equations in the eigenspaces of $A$.
The result is that
\[
  S_L: L^{2}(\R_{+},H_{\ne}) \to
  \{ \sigma\in Q_{\ne}(H^{1}(e)) \,:\, \sigma(0)\in H^{1/2}_< \}
\]
is continuous and bijective with inverse $L$.
We conclude:

\begin{lem}\label{lem:Sprops}
The solution operator
$S_D:=S_L\gamma^*:L^2(\R_+,H)\to H^{1}(e)$
of $D$ is continuous with $(Q_>S_D\sigma)(0)=0$ and
\begin{equation}\label{cc:sobsa1}
   DS_D\sigma = Q_{\ne}\sigma
\end{equation}
for all $\sigma\in L^2(\R_+,H)$.
Moreover,
\begin{equation}
     S_DD\sigma = Q_{\ne}\sigma
     \label{eq:sobsa1a}
\end{equation}
for all $\sigma\in H^{1}(e)$ with $Q_{>}\sigma(0)=0$.
In particular, the map
\begin{equation}\label{cc:sobsa2}
   \mathcal{R}S_D:L^2(\R_+,H) \to H^{1/2}_< ,
   \quad \sigma\mapsto (S_D\sigma)(0) ,
\end{equation}
is surjective. \qed
\end{lem}

We also use the {\em extension operator}
\begin{equation}\label{eq:extdef}
     \mathcal{E}x(t) := e^{-t(|A| +Q_{0})}x ,
\end{equation}
which is defined for $x\in H^{-\infty}$ and $t\ge0$.
We note that $\mathcal{E}x(t)\in H^{\infty}$ for all $t>0$.
The following assertions are readily verified by studying
the respective ordinary differential equations
in the eigenspaces of $A$.

\begin{lem}\label{cc:ne}
For any $s\in\R$ and $x\in H^s$,
\begin{enumerate}
\item
$\mathcal{E}x\in C(\R_+,H^s)$
and $||(\mathcal{E}x)(t)||_s\le||x||_s$ for all $t\ge0$.
\item
$\mathcal{E}x\in C^1(\R_+,H^{s-1})$
with $(\mathcal{E}x)'=-(|A|+Q_0)\mathcal{E}x$.
\item
$C_{s}^{-1}||x||_{s-\frac{1}{2}}
   \le|| (|A|+Q_{0})^{s}\mathcal{E}x||_{L^{2}(\R_{+},H)}
   \le C_{s}||x||_{s-\frac{1}{2}}$ . \qed
\end{enumerate}
\end{lem}

Since $||(|A|+Q_{0})x||_{s-1}\le||x||_s$,
the second equation implies that,
for any $x\in H_A=H^1$,
the extension $\mathcal{E}x:\R_+\to H$
is Lipschitz continuous with Lipschitz constant $1$.
In particular,
$\mathcal{E}x\in\mathcal L_{\loc}(e)$ for any $x\in H_A$.

\begin{prop}\label{eq:kermap}
The map
\[
  \mathcal{E}_{>}: H^{-1/2}_> \to \ker D_{\max} ,
  \quad \mathcal{E}_>x:=\mathcal{E}x ,
\]
is a continuous isomorphism.
The restriction map $\mathcal R$ extends to a continuous
map $\mathcal R$ on $\ker D_{\max}$
with $\mathcal R\mathcal{E}_>x=x$.
\end{prop}

\begin{proof}
It follows from Lemma \ref{cc:ne}.2
that $\mathcal{E}_{>}$ maps $H^{1}_>$
to $\ker D_{\max}$.
Lemma \ref{cc:ne}.3 implies that it extends
to $H^{-1/2}_>$ as a continuous and injective map,
where we recall that $\ker D_{\max}\subset L^2(\R_+,H)$
is closed.

To prove surjectivity,
choose a unitary basis $(e_n)$ of $H$
of eigenvectors of $A$, $Ae_n=a_ne_n$.
Let $\sigma\in\ker D_{\max}$ and set
$\sigma_n(t):=\la\sigma(t),e_n\ra$.
Then $\sigma_n$ solves the ordinary differential equation
$\sigma_n'+a_n\sigma_n=0$ weakly,
and hence $\sigma_n(t)=e^{-ta_n}x_n$,
where $x_{n} =\sigma_{n}(0)$.
Since $\sigma$ is square integrable,
$x_n=0$ for $a_n\le 0$
and $x=\sum_{a_{n}>0}x_{n}e_{n}\in H^{-1/2}_>$.
Hence $\sigma = \mathcal{E}_{>}x$.

The assertion about $\mathcal R$ is clear.
\end{proof}

We note that the Dirac operator $D$ commutes with $Q_0$
and $Q_{\ne}$,
hence
\begin{equation}\label{dmaxdecom}
  \mathcal D_{\max}
  = Q_{\ne}\mathcal D_{\max}\oplus Q_0\mathcal D_{\max} .
\end{equation}
Moreover $Q_0\mathcal D_{\max}=H^1(\R_+,Q_0H)$,
the standard Sobolev space.

\begin{cor}[Representation Formula]\label{cc:rf}
The map
\begin{align*}
  H^{-1/2}_>\oplus L^2(\R_+,H_{\ne})&\oplus H^1(\R_+,Q_0H)
  \to \mathcal D_{\max} , \\
  (x,\tau,\sigma_0) &\mapsto
  \sigma = \mathcal E_>x + S_D\tau + \sigma_0 ,
\end{align*}
is a continuous isomorphism
with $D_{\max}\sigma=\tau+\gamma\sigma_0'$.
\end{cor}

\begin{proof}
Clearly $\mathcal E_>x+S_D\tau+\sigma_0\in\mathcal D_{\max}$
for all $x\in H^{-1/2}_>$, $\tau\in L^2(\R_+,H_{\ne})$,
and $\sigma_0\in H_1(\R_+,Q_0H)$.
Vice versa, let $\sigma\in\mathcal D_{\max}$
and set $\tau=D_{\max}\sigma$ and $\sigma_0=Q_0\sigma$.
Then $\sigma-S_D\tau-\sigma_0\in\ker D_{\max}$,
by \lref{lem:Sprops}.
Hence our map is a continuous isomorphism,
by the continuity of $S_D$ and \pref{eq:kermap}.
\end{proof}

\begin{prop}[Boundary Values]\label{dmaxbou}
Let
\[
  \check H :=
  H^{-1/2}_> \oplus Q_{0}H \oplus H^{1/2}_< .
\]
Then $\mathcal{R}$ and $\mathcal{E}$ extend to
respectively define continuous operators
\[
  \mathcal{R}: \mathcal D_{\max} \to \check H
  \quad\text{and}\quad
  \mathcal{E}: \check H\to\mathcal{D}_{\max}
\]
with $\mathcal{R}\mathcal{E}=I$ on $\check H$.
In particular, $\mathcal{R}$ is surjective.
\qed
\end{prop}

Now we can derive the precise regularity properties of
elements in $\mathcal D_{\max}$ which will make
the special role of $0$ even more apparent.
For ease of notation, we set $\mathcal R\sigma=:\sigma(0)$.

\begin{lem}[Regularity II]\label{lem:regprop}
The maximal domain $\mathcal D_{\max}$ has the following
properties:
\begin{enumerate}
\item
$\mathcal{L}_{c}(e)$ is dense in $\mathcal D_{\max}$.
\item
$H^{1}(e) = \{\sigma\in\mathcal D_{\max}
 \,:\, \sigma(0)\in H^{1/2} \}\subset\mathcal D_{\max}$.
\item
$\mathcal D_{\max}\subset
  C(\R_{+},\check H) \cap C((0,\infty),H^{1/2})$.
\item
$\lim_{t\to\infty}\sigma(t)=0$ in $H^{1/2}$
for any $\sigma\in\mathcal D_{\max}$.
\item
If $\phi\in\Lip(\R_{+})$ is bounded
and $\sigma\in\mathcal D_{\max}$, \\
then $\phi\sigma\in\mathcal D_{\max}$
and $(\phi\sigma)(0)=\phi(0)\sigma(0)$.
\end{enumerate}
\end{lem}

\begin{proof}
(1) By definition,
$\mathcal{L}_{c}(e)$ is dense in $H^{1}(e)$.
Hence it suffices to consider $\sigma\in\ker D_{\max}$,
by \cref{cc:rf}.
Write $\sigma=\mathcal{E}_>x$ with $x\in H^{-1/2}_>$.
Choose a sequence $(x_n)$ in $H^1_>$
with $x_n\to x$ in $H^{-1/2}$
and $\phi\in\Lip_c(\R_+)$ with $\phi=1$ near $0$.
Set $\phi_n(t):=\phi(t/n)$, then by \lref{cc:ne}, 
$\phi_n\mathcal{E}_>x_n\in\mathcal{L}_c(e)$.
It is easy to see
that $\phi_n\mathcal{E}_>x_n\to\mathcal{E}_>x$
in $\mathcal{D}_{\max}$.

(2) Clearly $H^1(e)\subset\mathcal{D}_{\max}$.
Since the image of $S_D$ is contained in $H^1(e)$,
the asserted characterization of $H^1(e)$
is immediate from \lref{cc:ne}.3 and \cref{cc:rf}.

(3) $\mathcal D_{\max}\subset C(\R_{+},\check H)$
is clear from \lref{cc:ne}.1.
By \lref{lem:h1restr}.2, $H^1(e)$ is contained 
in $C(\R_+,H^{1/2})$, thus in $C(\R_+,\check H)$.
By \cref{cc:rf}, it is hence sufficient to
consider the image of $\mathcal E_>$.
Now $\mathcal Ex(t)\in H^{\infty}$
and $\mathcal Ex(t+t')=\mathcal E(\mathcal Ex(t))(t')$
for all $x\in H^{-1/2}$ and $t,t'>0$.
Hence $\mathcal E_>x\in C((0,\infty),H^{1/2})$
for all $x\in H^{-1/2}$, by \lref{cc:ne}.1.

(4) Let $\sigma\in\mathcal D_{\max}$.
It follows from (2) and (3) that $\sigma$ shifted by $t>0$,
$\tau_t\sigma(t'):=\sigma(t+t')$, is in $H^1(e)$.
Hence by \eqref{h1restr1},
\begin{align*}
  ||\sigma(t)||_{H^{1/2}}^2
  &= ||\tau_{t}\sigma(0)||_{H^{1/2}}^2
  \le C||\tau_{t}\sigma||_{H^{1}(e)}^2 \\
  &= C \int_{t}^{\infty}(||\sigma'||^{2}
  + ||A\sigma||^{2} + ||\sigma||^{2}) .
\end{align*}
Hence $\sigma(t)\to0$ in $H^{1/2}$ as $t\to\infty$.

(5) Let $\sigma\in\mathcal{D}_{\max}$
and $\phi\in\Lip(\R_+)$ be bounded.
Choose a sequence $(\sigma_n)$ in $\mathcal{L}_c(e)$
which converges to $\sigma$ in $\mathcal{D}_{\max}$.
Then $\phi\sigma_n\in\mathcal{L}_c(e)$
and $\phi\sigma_n\to\phi\sigma$ in $\mathcal{D}_{\max}$,
hence the claim.
\end{proof}

Now we can extend \eqref{eq:intpart1}
(cf. \cite[Lemma 2.15]{BL2}) to $\mathcal{D}_{\max}$.
We only have to use Part 1 of \lref{lem:regprop}
and to note that the skew-Hermitian form $\omega$
defined in \eqref{eq:intpart1} extends naturally
to $(x,y)\in\check H\times\check H$ by
\begin{equation}\label{defomegagen}
  \omega(x,y) :=
  B_{-1/2}(Q_{>\Lambda}x, \gamma Q_{<-\Lambda}y)
  + B_{1/2}(Q_{\le\Lambda}x, \gamma Q_{\ge-\Lambda}y) ,
\end{equation}
where $\Lambda\in\R$ is arbitrary.

\begin{cor}\label{intpart2}
For $\sigma_{1}, \sigma_{2}\in\mathcal D_{\max}$ we have
\begin{equation*}
  (D_{\max}\sigma_{1}, \sigma_{2})
  - (\sigma_{1}, D_{\max}\sigma_{2})
  = \omega(\sigma_{1}(0), \sigma_{2}(0)) .
  \qed
\end{equation*}
\end{cor}

We note that $\omega$ is non-degenerate on $\check H$.
For a linear subspace $B\subset\check H$,
the annihilator of $B$ with respect to $\omega$ is
\begin{equation}\label{annih}
  B^a := \{ y \in \check H \,:\,
  \text{ $\omega(x,y)=0$ for all $x\in B$} \} ;
\end{equation}
$B^a\subset\check H$ is closed,
and $B^{aa}$ is the closure of $B$ in $\check H$.
The description of $B^a$ is easy in the case
where $B$ is contained in $H^{1/2}$.

\begin{lem}\label{bcbann}
If $B\subset H^{1/2}\subset\check H$,
then $B^a=(\gamma B^0)\cap\check H$,
where
\[
  B^0 = \{y\in H^{-1/2} \,:\,
  \text{$B_{1/2}(x,y)=0$
  for all $x\in B$} \} .
\]
In particular,
$B^a\cap H^{1/2}=\gamma(B^\perp\cap H^{1/2})$,
where $B^\perp$ is the orthogonal complement
of $B\subset H$ in $H$.
\end{lem}

\begin{proof}
For $x,y\in\check H$ with $x\in H^{1/2}$,
we have $\omega(x,y)=B_{1/2}(x,\gamma y)$.
\end{proof}

\subsection{Boundary conditions and Fredholm properties}\label{boco}
With any linear subspace, $B\subset\check H$,
we now associate various extensions of $D_{0,c}$.
We define:
\begin{align}
  \mathcal L_{B,c}(e)
  :&= \{\sigma\in\mathcal L_c(e) \,:\, \sigma(0)\in B \} ,
  \label{bcbc} \\
  \quad D_{B,c}
  :&= D|\mathcal L_{B,c}(e) ;
  \notag\\
 \mathcal{D}_{B}
  :&= \{\sigma\in\mathcal D_{\max} \,:\,
       \sigma(0) \in B\cap H^{1/2} \}
  \label{bcb1} \\
  &= \{\sigma\in H^{1}(e) \,:\, \sigma(0) \in B \} ,
  \notag \\
  \quad D_{B}
  :&= D | \mathcal{D}_{B} ;
  \notag \\
 \mathcal{D}_{B,\max}
  :&= \{\sigma\in\mathcal D_{\max} \,:\, \sigma(0) \in B \} ,
  \label{bcbmax} \\
  \quad D_{B,\max}
  :&= D_{\max} |\mathcal{D}_{B,\max} .
  \notag
\end{align}
Since the restriction map
$\mathcal R:\mathcal D_{\max}\to\check H$ is continuous,
$D_{B,\max}$ is a closed operator if $B$
is a closed subspace of $\check H$.
Vice versa, we have:

\begin{prop}\label{bcbclex}
Let $\bar D\subset D_{\max}$ be a closed extension
of $D_{0,c}$ and $\bar{\mathcal D}$ be the domain
of $\bar D$.
Then $\bar D=D_{B,\max}$,
where $B=\{\sigma(0) \,:\, \sigma\in\bar{\mathcal D}\}$
is a closed subspace of $\check H$.
\end{prop}

\begin{proof}
Since $\bar D$ is a closed extension of $D_{0,c}$,
the closure of $\mathcal L_{0,c}(e)$ in the $H^1(e)$-norm
is contained in $\bar{\mathcal D}$,
\[
  H^1_0(e)
  := \{\sigma\in H^1(e) \,:\, \sigma(0)=0\}
  \subset \bar{\mathcal D} .
\]
Since the difference of any two elements from
$\mathcal D_{\max}$ with the same value at $0$
is in $H^1_0(e)$, by \lref{lem:regprop}.2,
we conclude that
\[
  \bar{\mathcal D}
  = \{\sigma\in\mathcal D_{\max}\,:\,\sigma(0)\in B\} ,
\]
hence that $\bar D=D_{B,\max}$.
Suppose now that $(x_n)$ is a sequence in $B$
such that $x_n\to x$ in $\check H$.
Then, by what we just said,
$(\mathcal Ex_n)$ is a sequence in $\bar{\mathcal D}$
and $\mathcal Ex_n\to\mathcal Ex$ in $\mathcal D_{\max}$,
by \pref{dmaxbou}.
Since $\bar D$ is a closed operator
and $\mathcal R$ is continuous,
we get that $x\in B$.
\end{proof}

\begin{dfn}\label{bcbdef}
A (linear) {\em boundary condition} for a Dirac system
is a closed linear subspace $B\subset\check H$.
\end{dfn}

\begin{rem}\label{bcbaps}
Since the seminal article \cite{APS}
of Atiyah, Patodi, and Singer,
it is customary to state boundary conditions
for Dirac systems in terms of projections in $H$.
In our setup, the boundary condition introduced
by Atiyah, Patodi, and Singer
is given by the subspace $B_{APS}:=\check H_\le$
of $\check H$.
We will discuss boundary conditions given by projections
in Section \ref{pop}.
Our approach to the description of boundary conditions
for Dirac systems, however, does not only seem to be
more general but will also lead to a more satisfying
analysis of the corresponding operators,
as we are going to explain.
\end{rem}

For any $\sigma\in\mathcal L_c(e)$, $\sigma(0)\in H_A=H^1$.
Vice versa, for any $x\in H_A$
there is $\sigma\in\mathcal L_c(e)$ with $\sigma(0)=x$.
Similarly, for any $x\in H^{1/2}$
there is $\sigma\in H^1(e)$ with $\sigma(0)=x$.
Let $B\subset\check H$ be a boundary condition.
We conclude, using \eqref{intpart2},
that the adjoint operators of the above operators are
\begin{align}
  (D_{B,c})^* &= D_{B_1,\max}
  \quad\text{with $B_1=(B\cap H_A)^a$} ,
  \label{bcbcad} \\
  (D_{B})^* &= D_{B_2,\max}
  \quad\text{with $B_2=(B\cap H^{1/2})^a$} ,
  \label{bcb1ad} \\
  (D_{B,\max})^* &= D_{B^a,\max} .
  \label{bcbmaxad}
\end{align}
Since the closure of a linear subspace of $\check H$
is the annihilator of its annihilator,
the closures of the above operators are
\begin{align}
  D_{B,\min}
  = (D_{B,c})^{**}
  &= D_{C_1,\max} ,
  \label{bcbcclos} \\
  (D_{B})^{**}
  &= D_{C_2,\max} ,
  \label{bcb1clos} \\
 (D_{B,\max})^{**}
  &= D_{B,\max} .
  \label{bcbmaxclos}
\end{align}
where $C_1$ is the closure of $B\cap H_A$
in $\check H$ in \eqref{bcbcclos}
and $C_2$ is the closure of $B\cap H^{1/2}$ 
in $\check H$ in \eqref{bcb1clos}.
In particular,
\begin{equation}\label{bcbminmax}
  D_{B,\min} = D_{B,\max}
  \Longleftrightarrow
  \text{$B\cap H_A$ is dense in $B$} .
\end{equation}

\begin{dfn}\label{bcbreg}
We say that a boundary condition $B\subset\check H$
is {\em regular} if $D_{B,\max}=D_B$.
We say that a boundary condition $B$ is
{\em elliptic} if $B$ and $B^a$ are regular.
\end{dfn}

By the representation formula \ref{cc:rf},
the boundary condition $B_{APS}=\check H_{\le}$
of Atiyah, Patodi, and Singer is the most natural
regular boundary condition.
The following reformulations of regularity are
immediate from the definition of regularity
and the properties of the maximal domain.

\begin{prop}\label{bcbregeq}
A closed linear subspace $B$ of $\check H$ is a
regular boundary condition iff any of the following
equivalent conditions is satisfied:
\begin{enumerate}
\item
$D_{B,\max}=D_B$.
\item
$\mathcal D_{B,\max}\subset H^1(e)$.
\item
$B\subset H^{1/2}\subset\check H$.
\end{enumerate}
A closed linear subspace $B$ of $H^{1/2}$
is a regular boundary condition iff one of the following
two equivalent conditions is satisfied:
\begin{enumerate}
\item[(4)]
The $H^{1/2}$ and $\check H$-norms
are equivalent on $B$.
\item[(5)]
For some or, equivalently, any $\Lambda\in\R$, 
there is a constant $C$ such that, for all $x\in B$,
\begin{equation*}
  ||Q_{>\Lambda}x||_{1/2}
  \le C ( ||Q_{>\Lambda}x||_{-1/2} 
      + ||Q_{\le\Lambda} x||_{1/2} ) .
  \qed
\end{equation*}
\end{enumerate}
\end{prop}

\begin{lem}\label{closeext}
Let $B\subset\check H$ be a regular boundary condition
and $\Lambda\in\R$.
Then the map $Q_{\le\Lambda}:B \to H^{1/2}_{\le\Lambda}$
is a left-Fredholm operator, that is,
has finite-dimensional kernel and closed image.
Moreover,
$(H_{>\Lambda}^{1/2},B)$ is a left-Fredholm pair
in $H^{1/2}$ with
\begin{align*}
  \nuli(H_{>\Lambda}^{1/2},B)
  &= \dim\ker (Q_{\le\Lambda}:B 
     \to H^{1/2}_{\le\Lambda}) 
  = \dim(H_{>\Lambda}^{1/2} \cap B) , \\
  \defi(H_{>\Lambda}^{1/2},B)
  &= \dim\coker (Q_{\le\Lambda}:B \to H^{1/2}_{\le\Lambda})
  = \dim (\check H_{\ge-\Lambda} \cap B^a) .
\end{align*}
\end{lem}

\begin{proof}
We use H\"ormander's Criterion, see \lref{trfred}.
Suppose that $(x_n)$ is a bounded sequence in $B$
such that $Q_{\le\Lambda}(x_n)$ converges 
in $H^{1/2}_{\le\Lambda}$.
Since the inclusion $H^{1/2}\to H^{-1/2}$ is compact
and $(x_n)$ is bounded in $H^{1/2}$,
we may assume,
by passing to a subsequence if necessary,
that $(x_n)$ converges in $H^{-1/2}$.
But then $(Q_{>\Lambda}x_n)$ is a Cauchy sequence 
in $H^{1/2}$, by \pref{bcbregeq}.5.
It follows that 
$Q_{\le\Lambda}:B \to H^{1/2}_{\le\Lambda}$
is a left-Fredholm operator
and hence, by \pref{bp2},
that $(H_{>\Lambda}^{1/2},B)$ is a left-Fredholm pair.
The formulas for the nullity and the first formula 
for the deficiency of the pair $(H_{>\Lambda}^{1/2},B)$
are clear. 
As for the last equality, we have, using \eqref{fgd},
\begin{align*}
  (H_{>\Lambda}^{1/2}+B)^0
  &= (H_{>\Lambda}^{1/2})^0 \cap B^0
  = H_{\le\Lambda}^{-1/2} \cap B^0 \\
  &= \{ x\in H_{\le\Lambda}^{-1/2} \,:\,
    \text{$B_{-1/2}(x,y)=0$ for all $y\in B$} \} \\
  &= \gamma(\{ x\in\check H_{\ge-\Lambda} \,:\,
    \text{$\omega(x,y)=0$ for all $y\in B$} \} ) \\
  &= \gamma(\check H_{\ge-\Lambda} \cap B^a) .
  \qedhere
\end{align*}
\end{proof}

\begin{prop}\label{bcbrdata}
Let  $\Lambda$ be a real number,
$\hat U\subset H_{\le\Lambda}^{1/2}$ be a closed subspace, 
$F\subset H^{1/2}_{<-\Lambda}$ be a finite-dimensional subspace, 
$\hat V:=F^0\cap H^{1/2}_{<-\Lambda}$,
and let $g:\hat U\to\hat V$ be a continuous linear map.
Then
\[
  B = \gamma F \oplus \{ u+\gamma gu \,:\, u\in\hat U \}
\]
is a regular boundary condition,
and all regular boundary conditions arise in this way.
\end{prop}

\begin{proof}
It is clear that any boundary condition $B$ of the
given form is regular.
Conversely, let $B\subset H^{1/2}$ be a regular
boundary condition.
By \lref{closeext}, 
\[
  \hat U := \im(Q_{\le\Lambda}:B\to H^{1/2}_{\le\Lambda})
\]
is a closed subspace of $H_{\le\Lambda}^{1/2}$ and
\[
  F := \gamma (B\cap H^{1/2}_{>\Lambda})
  = \gamma (\ker(Q_{\le\Lambda}:B\to H^{1/2}_{\le\Lambda}))
\]
is a finite-dimensional subspace of $H^{1/2}_{<-\Lambda}$.
It follows that $G=(\gamma F^\perp)\cap B$ is a complement
of $\gamma F$ in $B$ and 
that $Q_{\le\Lambda}:G\to\hat U$ is an isomorphism.
Hence there is a continuous linear map
$g:\hat U\to H^{1/2}_{<-\Lambda}$ 
such that
\[
  G = \{ u+\gamma gu \,:\, u\in\hat U \} .
\]
Since $G\subset\gamma F^\perp$, 
$g$ takes values in $\hat V$.
\end{proof}

\begin{rem}\label{bcbregrem}
In \pref{bcbrdata} above and \pref{bcbsrdata} below,
the roles of weak and strong inequalities 
can be interchanged.
\end{rem}

\begin{prop}\label{bcbsrdata}
Let $\Lambda$ be a real number and let
\[
  H_{\le\Lambda} = E \oplus U
  \quad\text{and}\quad
  H_{<-\Lambda} = F\oplus V
\]
be orthogonal decompositions,
where $E,F\subset H^{1/2}_{<-\Lambda}$ are 
finite-dimensional subspaces, 
and $g: U \to V$ be a $1/2$-smooth linear map.
Then
\[
  B = \gamma F \oplus 
      \{u+\gamma gu \,:\, u\in U\cap H^{1/2} \}
\]
is an elliptic boundary condition with
\[
  B^a = \gamma E \oplus
        \{v+\gamma g^*v \,:\; v\in V \cap H^{1/2} \} .
\]
All elliptic boundary conditions arise in this way.
\end{prop}

\begin{rem}\label{baeb}
In previous work, but in a different context,
the first author and Christian B\"ar
considered boundary conditions of precisely this form.
For details see the forthcoming article \cite{BB}.
\end{rem}

\begin{proof}[Proof of \pref{bcbsrdata}]
With data as in \pref{bcbrdata}, write
\[
  B = \gamma F \oplus \{u+\gamma\hat gu \,:\, u\in\hat U \} ,
\]
where the map $g$ there is decorated with a hat here.
Since $F\subset H^{1/2}_{<-\Lambda} $ 
is of finite dimension,
\begin{equation*}
  F\oplus\tilde V = H^{-1/2}_{<-\Lambda} ,\quad
  F\oplus V = H_{<-\Lambda} , \quad
  F\oplus\hat V = H^{1/2}_{<-\Lambda} , 
\end{equation*}
where 
\[
  \tilde V = F^0\cap H^{-1/2}_{<-\Lambda} , \quad
  V = F^0\cap H_{<-\Lambda} , \quad
  \hat V = F^0\cap H^{1/2}_{<-\Lambda} .
\]
Let $x\in\gamma B^0\subset H^{-1/2}$.
Then there exist $f\in F$ and $v\in\tilde V$
with $Q_{<-\Lambda}x=f+v$.
We compute $B_{-1/2}(x,f)=|f|^2$.
Since $f\in\gamma B$, we conclude that $f=0$
and hence that 
\[
  Q_{<-\Lambda}(\gamma B^0)\subset\tilde V .
\]
Conversely, let $v\in\tilde V$.
Then $B_{-1/2}(v+\gamma w,f)=0$ 
for all $w\in H_{\le\Lambda}^{-1/2}$
and $f\in F$, by the definition of $\tilde V$ 
and since $F\subset H^{1/2}_{<-\Lambda}$.
With $u\in\hat U$, we compute
\begin{align}
  B_{-1/2}(v+\gamma w,\gamma u-\hat gu)
  &= B_{-1/2}(\gamma w,\gamma u) - B_{-1/2}(v,\hat gu) 
  \notag \\
  &= B_{-1/2}(w,u) - B_{-1/2}(v,\hat gu) 
  \notag \\
  &= B_{-1/2}(w,u) - B_{-1/2}(u',u) \notag
\end{align}
for some appropriate $u'\in H^{-1/2}_{\le\Lambda}$,
by the duality 
$(H^{1/2}_{\le\Lambda})'=H^{-1/2}_{\le\Lambda}$.
We conclude that $v+\gamma u'\in\gamma B^0$.
In particular,
\[
  \tilde V = Q_{<-\Lambda}(\gamma B^0) .
\]
Since
$\check H=H^{1/2}_{<-\Lambda}\oplus H^{-1/2}_{\ge\Lambda}$,
we have $v+\gamma u'\in\check H$ if and only 
if $v\in H^{1/2}_{<-\Lambda}$.

We now use that $B$ is elliptic.
Then $B^a=(\gamma B^0)\cap\check H$ is regular
and hence 
$(\gamma B^0)\cap\check H=(\gamma B^0)\cap H^{1/2}$.
It follows that $v+\gamma u'\in\gamma B^0$ as above
belongs to $H^{1/2}$ if and only if $v\in H^{1/2}$,
and therefore
\[
  \hat V = Q_{<-\Lambda}(B^a) .
\]
By the symmetry of the roles of $B=(B^a)^a$ and $B^a$ 
and switching the roles of weak and strong inequalities,
see Remark \ref{bcbregrem}, we get
\[
  \hat U =  Q_{\le\Lambda}(B) 
  = E^0 \cap H^{1/2}_{\le\Lambda} ,
\]
where $E=\gamma(B^a\cap H^{1/2}_{\ge\Lambda})$.
By \lref{closeext}, $E$ is finite-dimensional.
Hence the sesquilinear form $B_{-1/2}$ identifies
$\tilde U = E^0 \cap H^{-1/2}_{\le\Lambda}$
with the dual space of $\hat U$.
In particular, in the above $v+u'$,
we may take $u'=\hat g'v$, 
where $\hat g':\tilde V\to\tilde U$ is the
dual map of $\hat g$.

We now recall that $u'=\hat g'v$ is in $H^{1/2}$
if $v\in H^{1/2}$, by the regularity of $B^a$.
By interpolation we get that $\hat g'$ is the extension 
of a $1/2$-smooth linear map $g^*:V\to U$.
By symmetry, $g^*$ is the adjoint of a $1/2$-smooth
map $g:U\to V$ and $\hat g$ is the restriction of $g$
to $\hat U$.
\end{proof}

\begin{cor}\label{bcbsrind}
Let $B\subset\check H$ be an elliptic boundary 
condition and $\Lambda\in\R$.
Then $\gamma B^\perp$ is the closure of $B^a$ in $H$ and
\begin{equation}
 \bar B\cap H_{\ge\Lambda}
 = B\cap H_{\ge\Lambda}^{1/2} , \quad
 B^\perp\cap H_{<\Lambda}
 = \gamma (B^a\cap H_{>-\Lambda}^{1/2} ) , \tag{1}
\end{equation}
where $\bar B$ denotes the closure of $B$ in $H$.
Moreover, $(\bar B,H_{\ge\Lambda})$ is a Fredholm pair
in $H$ with index
\begin{align}
  \ind(\bar B,H_{\ge\Lambda})
  &= \dim (\bar B\cap H_{\ge\Lambda})
    - \dim (B^\perp\cap H_{<\Lambda}) \tag{2} \\
  &= \dim (B\cap H_{\ge\Lambda}^{1/2})
    - \dim (B^a\cap H_{>-\Lambda}^{1/2}) . \notag
\qed
\end{align}
\end{cor}

It is natural to ask whether the index formula
in (\ref{bcbsrind}.2) gives the index of $D_{B,\max}$
for suitable $\Lambda$;
this is in fact true for $\Lambda=0$ if $\ker A=0$.

\begin{prop}\label{bcbreg0}
Let $\ker A=0$.
If $B\subset\check H$ is a regular boundary condition,
then $D_{B}=D_{B,\max}$ is a left-Fredholm operator
with
\begin{equation}
  (\im D_{B})^\perp = \ker D_{B^a,\max} . \tag{1}
\end{equation}
If $B$ is elliptic, then $D_{B}$ is a Fredholm 
operator with
\begin{equation}
  \ind D_B = 
  \dim \bar B \cap H_\ge - \dim B^\perp \cap H_< .
  \tag{2}
\end{equation}
\end{prop}

\begin{proof}
We again use H\"ormander's Criterion
from \lref{trfred}.
Since the kernel of $A$ vanishes,
we have the representation formula
\[
  \sigma
  = \mathcal EQ_>\sigma(0) + S_DD_{\max}\sigma ,
\]
characterizing elements
$\sigma\in\mathcal D_{\max}$.
Furthermore,
\[
  S_D:L^2(\R_+,H) \to
  \{ \sigma \in H^1(e) \,:\, Q_>\sigma(0)=0 \}
\]
is an isomorphism, by \lref{lem:Sprops}.
Let $(\sigma_n)$ be a bounded sequence
in $\mathcal D_{B,\max}$
such that $D_{\max}\sigma_n$ converges in $L^2(\R_+,H)$.
Then $(\sigma_n(0))$ is a bounded sequence in $B$
and $(S_DD_{\max}\sigma_n)$ converges in $H^1(e)$.
It follows that the sequence
$(Q_\le\sigma_n(0)=(S_DD_{\max}\sigma_n)(0))$
converges in $H^{1/2}_\le$.
By \lref{closeext} and H\"ormander's criterion again, 
$(\sigma_n(0))$ has a convergent subsequence in $B$.
Hence
$(\sigma_n=\mathcal EQ_>\sigma_n(0)+S_DD_{\max}\sigma_n)$
has a convergent subsequence in $\mathcal D_{B,\max}$.
This shows that $D_B$ is a left-Fredholm operator.
Now $D_{B}^*=D_{B^a,\max}$, see \eqref{bcb1ad},
therefore $(\im D_{B})^\perp = \ker D_{B^a,\max}$
as claimed.
\end{proof}

We note that the image of $D_{B,\max}$ is not closed
if $\ker A\ne0$ while the index formula
in (\ref{bcbsrind}.2) holds in general.
This suggests a possible extension of \pref{bcbreg0}
which we achieve by conveniently enlarging the domain 
of $D_{\max}$.
We recall that $Q_0$ and $Q_{\ne}$ commute with $D_{\max}$
and that $\mathcal D_{\max}$ splits perpendicularly
with components $H^1(\R_+,Q_0H)$ and $Q_{\ne}\mathcal D_{\max}$.
As is well known, the source of trouble is the part
\[
  D_{\max}: H^1(\R_+,Q_0H) \to L^2(\R_+,Q_0H) .
\]
of $D_{\max}$.
We restore Fredholm properties of $D$
by enlarging $H^1(\R_+,Q_0H)$.
Our discussion is motivated by the work of the third
author on non-parabolic Dirac operators,
compare \cite{Ca2} and Section \ref{bocogen} below.

By \cref{cc:rf},
we have equivalences of norms on $\mathcal D_{\max}$,
\begin{align}
  ||\sigma||_{D_{\max}}^2
  &\approx
  ||Q_>\sigma(0)||_{-1/2}^2
  + ||\tau||_{L^2(\R_+,H_{\ne})}^2
  + ||\sigma_0||_{H^1(\R_+,Q_0H)}^2 \notag \\
  &=
  ||Q_>\sigma(0)||_{-1/2}^2
  + ||D_{\max}\sigma||_{L^2(\R_+,H)}^2
  + ||\sigma_0||_{L^2(\R_+,Q_0H)}^2 ,
  \label{equivnorm} \\
  &\approx
  ||\sigma(0)||_{\check H}^2
  + ||D_{\max}\sigma||_{L^2(\R_+,H)}^2
  + ||\sigma_0||_{L^2(\R_+,Q_0H)}^2 ,
  \label{equivnorm2}
\end{align}
where $\tau=D_{\max}Q_{\ne}\sigma$ and $\sigma_0=Q_0\sigma$
and where we note, for the last equivalence,
that $\mathcal R$ is continuous
on $\mathcal D_{\max}$.
We now introduce a continuous seminorm $||\cdot||_W$
on $\mathcal D_{\max}$,
\begin{equation}\label{eq:west}
  ||\sigma||_W^2
  := ||\sigma(0)||_{\check H}^2 +
  ||D_{\max}\sigma||_{L^2(\R_+,H)}^2
  \le C\cdot ||\sigma||^2_{D_{\max}} .
\end{equation}
\cref{cc:rf} implies that $||\cdot||_W$
is actually a norm on $\mathcal D_{\max}$.
Clearly, $||\cdot||_W$ and the graph norm of $D_{\max}$
are equivalent if $\ker A=0$.
On the other hand, if $\ker A\ne 0$,
then $||\cdot||_W$ is strictly weaker
than the graph norm of $D_{\max}$.
However, one easily verifies that for any $T>0$
there is a constant $C_T$ such that
\begin{equation}\label{eq:west0}
   ||\sigma||_{L^{2}([0,T],H)}
   \le C_{T}||\sigma||_{W} ,
\end{equation}
for all $\sigma\in\mathcal D_{\max}$.

We now let $W$ be the closure of $\mathcal D_{\max}$
under the norm $||\cdot ||_{W}$.
By \lref{lem:regprop}.1,
$\mathcal{L}_{c}(e)$ is dense in $W$.
By definition, $D_{\max}$ extends to a continuous
operator
\begin{equation}\label{dext}
   D_{\ext}: W \to L^2(\R_+,H) .
\end{equation}
We observe now that
\begin{equation}\label{wsplit}
  W = Q_{\ne}W \oplus Q_0W
  = Q_{\ne}\mathcal D_{\max} \oplus Q_0W .
\end{equation}
The linear map
$S_0: Q_0H \oplus L^2(\R_+,Q_0H) \to Q_0W$
defined by
\begin{equation}
  S_0(x,\tau)(t) := x + \gamma^*\int_0^t\tau(s) ds ,
\end{equation}
is an isomorphism with $D_{\ext}S_0(z,\tau)=\tau$.
In particular,
\begin{equation}\label{w0h1}
  Q_0W \subset H^1_{\loc}(\R_+,Q_0H) .
\end{equation}
With $\mathcal R(S_0(x,\tau)):=x$
we obtain a continuous extension
\begin{equation}\label{rext}
  \mathcal R: W \to\check H ,
  \quad
  \mathcal R\sigma =: \sigma(0) ,
\end{equation}
of $\mathcal R$ to $W$.
For a boundary condition $B\subset\check H$, we set
\begin{equation}\label{wbc}
  W_B := \{ \sigma\in W \,:\, \sigma(0)\in B \}
  \quad\text{and}\quad
  D_{B,\ext} := D_{\ext}|W_B .
\end{equation}
We see from the above
that $L^2(\R_+,Q_0H)\subset\im D_{B,\ext}$,
irrespective of the boundary condition $B$.

\begin{thm}\label{wfred}
If $B$ is regular,
then $D_{B,\ext}$ is a left-Fredholm operator
with $(\im D_{B,\ext})^{\perp} = \ker D_{B^a,\max}$.
\end{thm}

\begin{proof}
Use the representation
$\mathcal Ex + S_D\tau + S_0(y,\rho)$
of elements of $W$,
where $x\in H^{-1/2}_>$, $\tau\in L^2(\R,Q_{\ne}H)$,
$y\in Q_0H$, and $\rho\in L^2(\R,Q_0H)$,
and adapt the argument from the proof of \pref{bcbreg0}.
\end{proof}

For any boundary condition $B\subset\check H$,
\begin{equation}\label{bcbker}
\begin{split}
  \ker D_{B,\max}
  &= \mathcal D_{B,\max} \cap \ker D_{\max} , \\
  \ker D_{B,\ext}
  &= W_B \cap \ker D_{\ext}
   = W_B \cap (\ker D_{\max}+Q_0H) .
\end{split}
\end{equation}
In particular, we have isomorphisms
\begin{equation}\label{bcbker2}
\begin{split}
  \mathcal R: &\ker D_{B,\max} \to B\cap \check H_> , \\
  \mathcal R: &\ker D_{B,\ext} \to B\cap \check H_\ge .
\end{split}
\end{equation}
Recall that a boundary condition $B$ is elliptic
if $B$ and $B^a$ are regular.
As above, we let $\bar B$ denote the closure of $B$ in $H$.

\begin{nameit}{Corollary and Definition}
\label{bcbsr}
If $B$ is elliptic,
then $D_{B,\ext}$ is a Fredholm operator
with index
\begin{equation*}
  \ind D_{B,\ext}
 = \dim (B\cap H_\ge)
  - \dim (B^\perp\cap H_<)
 = \ind (\bar B,H_\ge) ,
\end{equation*}
the {\em extended index} of $D_{B}$,
also denoted by $\ind_{\ext}D_B$
\end{nameit}

\begin{proof}
Immediate from \eqref{bcbker2}, \tref{wfred},
and \cref{bcbsrind}.
\end{proof}

\subsection{Self-adjoint boundary conditions}
We say that a boundary condition $B\subset\check H$
is {\em self-adjoint} if $B=B^a$.
By definition, a regular self-adjoint boundary condition
is elliptic.

We say that $(H_0,\omega)$ is a 
{\em Hermitian symplectic vector space}
if the $\pm1$-eigenspaces of the involution $i\gamma$
of $H_0$ have equal dimension.
Then a subspace $L\subset H_0$ is {\em Lagrangian}
if $L\perp\gamma L$ and $L\oplus\gamma L=H_0$.

\begin{thm}\label{bcsrsa}
Regular self-adjoint boundary conditions exist
if and only if $(H_0,\omega)$ is a Hermitian symplectic vector 
space (where $H_0=0$ is not excluded).
Then regular self-adjoint boundary conditions are given
by the following data: a Lagrangian subspace $L\subset H_0$,
an orthogonal decomposition $H_<=F\oplus V$,
where $F\subset H_<^{1/2}$ is of finite dimension,
and a $1/2$-smooth map $g:V\oplus L\to V\oplus L$ with $g^*=g$.
The regular self-adjoint boundary condition $B$ 
given by such data is
\[
  B = \gamma F \oplus
  \{ w + \gamma gw \,:\, w\in (V\oplus L)\cap H^{1/2} \} .
\]
\end{thm}

Write $H=H^+\oplus H^-$, where $H^\pm$ is the $\pm1$
eigenspace of $i\gamma$.
Since $A$ anti-commutes with $\gamma$,
$A$ maps $H^\pm$ to $H^\mp$
so that the restriction of $A$ to $H^+$ is
a Fredholm operator (in general unbounded) to $H^-$.
Since $\gamma$ intertwines eigenspaces of $A$
with opposite eigenvalues,
it follows easily that $(H_0,\omega)$ 
is a Hermitian symplectic vector space if and only if 
the Fredholm operator $A^+$ has index $0$.

\begin{cor}\label{bcsrsacor}
With $H^\pm$ and $A^+$ as above,
$\check H$ contains elliptic self-adjoint boundary
conditions if and only if $\ind A^+=0$. 
\qed
\end{cor}

\begin{proof}[Proof of \tref{bcsrsa}]
Any data as in the assertion give rise to
a regular self-adjoint boundary condition.
As for the existence, 
if $L\subset H_0$ is a Lagrangian subspace,
then $L\oplus H_<$ is a regular self-adjoint
boundary condition.

To prove the asserted characterization,
we first observe that regular self-adjoint boundary
conditions are elliptic,
so that we can use the description of elliptic
boundary conditions given in \pref{bcbsrdata}.

Let $B$ be an elliptic boundary condition.
By \pref{bcbsrdata},
there are orthogonal decompositions
\[
  H_{\le} = E \oplus U 
  \quad\text{and}\quad
  H_< = F \oplus V ,
\]
where $E,F\subset H^{1/2}$ are of finite dimension,
and a $1/2$-smooth linear map $b:U\to V$
such that
\begin{align*}
  B &= \gamma F \oplus 
        \{ u+\gamma bu \,:\, u \in U\cap H^{1/2} \} , \\
  B^a &= \gamma E \oplus
          \{ v+\gamma b^*v \,:\, v \in V\cap H^{1/2} \} .
\end{align*}
 From now on we assume that $B=B^a$.
Then the $H$-closure $\bar B=\gamma B^\perp$,
and hence any element in $\bar B$  can be written
in any of the following two ways:
\begin{equation*}
  \gamma f + u + \gamma bu
  = \gamma f + u_< + u_0 + \gamma bu ,
\end{equation*}
where $u_<=Q_<u$ and $u_0=Q_0u$, and
\begin{equation*}
  \gamma e + v + \gamma b^*v
  = \gamma e_< + \gamma e_0 + v + \gamma b^*v ,
\end{equation*}
where $e_<=Q_<e$ and $e_0=Q_0e$.
We are going to compare the $H_<$, $H_0$, and $H_>$ components 
of elements of $\bar B$ in the above two representation: 

We observe first that $V=Q_<(U)=\{u_< \,:\, u\in U \}$.
Since $E$ and $F$ are the orthogonal complements of $U$ 
in $H_\le$ and $V$ in $H_<$,
it follows that 
\[
  F = E\cap H_< \subset E .
\]
Let $L:=U\cap H_0$ 
and 
\[
  B_L := \{u+\gamma bu \,:\, u\in L\}
  \subset\bar B .
\]
Let $u\in L$.
Then $u+\gamma bu\in\bar B$ and hence there exist $e\in E$ 
and $v\in V$ such that
\[
  u + \gamma bu =
  \gamma e + v + \gamma b^*v .
\] 
Clearly $v=0$, hence $b^*v=0$,
and hence $u=\gamma e_0$ and $bu=e_<$.
We get $\gamma u-bu=-e$ and hence 
\[
  \gamma B_L
  = \{\gamma u-bu \,:\, u\in L\} \subset E .
\]
Let $e\in E$.
Then $\gamma e\in\bar B$ and hence there exist $f\in F$ 
and $u\in U$ such that
\[
  \gamma e = \gamma f + u_< + u_0 + \gamma bu .
\]
We obtain $u_<=0$, hence $\gamma e_0=u_0=u\in L$
and $e_<=f+bu$. 
Since $F\subset E$, we get
\[
  E = F \oplus \{\gamma u-bu \,:\, u\in L\} 
    = F \oplus \gamma B_L .
\]
Since $U$ is the orthogonal complement of $E$ in $H_\le$
and $Q_0(E)=\gamma L$,
the orthogonal complement of $\gamma L$ in $H_0$
belongs to $U$, that is, to $L$,
by the definition of $L$.
We conclude that we have an orthogonal decomposition
\[
  H_0 = L \oplus \gamma L .
\]
It follows that $(H_0,\omega)$ is a Hermitian symplectic vector
space and that $L$ is a Lagrangian subspace of $(H_0,\omega)$.

Since $\bar B=\gamma B^\perp$ and $E=F\oplus\gamma B_L$, 
we have orthogonal sums
\[
  \bar B = \gamma E \oplus
     \{ v+\gamma b^*v \,:\, v \in V \}
  =: \gamma F \oplus B_L \oplus B_V .
\]
Let $W:=V\oplus L$.
Then $H$ decomposes orthogonally as
\[
  H = F \oplus W \oplus \gamma F \oplus \gamma W .
\]
For a subspace $K\subset H$, 
let $Q_K$ be the orthogonal projection in $H$ onto $K$.
Then $Q_W=Q_V+Q_L$.

Let $x\in B_L\oplus B_V$
and write $x=u+\gamma bu+v+\gamma b^*v$
with $u\in L$ and $v\in V$.
Since $Q_F(b^*v)=Q_F(\gamma b^*v)=0$, 
we have 
\begin{align*}
  \gamma Q_Wb^*v 
  = \gamma Q_{F\oplus W} b^*v
  = (I-Q_{F\oplus W})\gamma b^*v 
  = (I-Q_W)\gamma b^*v .
\end{align*}
Therefore
\begin{align*}
  x &= u + \gamma bu + v + \gamma b^*v \\
  &= u + Q_W\gamma b^*v + v 
     + \gamma bu + (I-Q_W)\gamma b^*v \\
  &= u + Q_W\gamma b^*v + v
     + \gamma(bu + Q_Wb^*v) .
\end{align*}
Since $\gamma b^*v\in H_\ge$, 
we have $Q_W\gamma b^*v=Q_L\gamma b^*v$.
Hence
\begin{align*}
  x &= (u + Q_L\gamma b^*v) + v
   + \gamma\big( b(u+Q_L\gamma b^*v)
   + (Q_Wb^*- bQ_L\gamma b^*) v \big) \\
  &= (Q_L+Q_V)x + \gamma g (Q_Wx) 
  = Q_Wx + \gamma g(Q_Wx) ,
\end{align*}
where $g:W\to W$ is the $1/2$-smooth linear map given by
\[
  gw = bQ_Lw + (Q_Wb^*-bQ_L\gamma b^*)Q_Vw .
\]
We conclude that
\[
  \bar B = \gamma F \oplus
  \{ w + \gamma gw \,:\, w\in W \} .
\]
Now
\[
  \gamma \bar B = B^\perp 
  = F \oplus \{ \gamma w - g^*w \,:\, w\in W \} ,
\]
hence $g=g^*$.
\end{proof}

\begin{exa}\label{bcsrsaexa}
Let $\beta:H\to H$ be $1/2$-smooth with
\begin{align*}
  \beta^* = \beta^{-1} &= \beta , \tag{1} \\
  \gamma\beta + \beta\gamma &= 0 , \tag{2} \\
  A\beta + \beta A &= 0 . \tag{3}
\end{align*}
Then $B=\{x\in H^{1/2} \,:\, \beta x=x\}$
is a regular self-adjoint boundary condition.

For example, given a Dirac system $d=(H,A,\gamma)$,
consider the Dirac system 
\[
  \tilde d = (H\oplus H,(A,-A),(\gamma,-\gamma)) .
\]
Then $\beta:H\oplus H\to H\oplus H$, $\beta(x,y)=(y,x)$,
satisfies (1)--(3).
The corresponding boundary condition 
$B=\{(x,x) \,:\, x\in H^{1/2} \}$
is regular and self-adjoint.
It arises as the transmission boundary condition
when cutting a manifold along a hypersurface.
\end{exa}

\subsection{Regular pairs of projections}\label{pop}
Let $P$ and $Q$ be $1/2$-smooth projections in $H$.
We say that the ordered pair $(P, Q)$ is {\em regular} if
\begin{equation}\label{cc:regular2}
   x\in H^{-1/2} ,\, \tilde Px=0 ,\,
   \tilde Qx\in H^{1/2}
   \Longrightarrow x\in H^{1/2} .
\end{equation}
Roughly speaking,
this means that $Q$ is close to $I-P$;
compare \pref{cc:cr} below.

\begin{lem}\label{lem:reg2}
Let $(P,Q)$ be a pair of 1/2-smooth projections in $H$.
Then $(P,Q)$ is regular if and only if
\begin{equation*}
  x\in H^{-1/2} ,\, \tilde Px\in H^{1/2} ,\,
  \tilde Qx\in H^{1/2}
  \Longrightarrow x\in H^{1/2} .
\end{equation*}
\end{lem}

\begin{proof}
Assume that $(P,Q)$ is regular.
Consider $x\in H^{-1/2}$ with $\tilde{P}x$
and $\tilde{Q}x$ in $H^{1/2}$.
Set $y:=(I - \tilde P)x\in H^{-1/2}$.
Then $\tilde{P} y=0$ and
\begin{equation*}
    \tilde{Q} y
    = \tilde{Q} x - \tilde{Q}\tilde{P} x
    = \tilde{Q} x - \hat{Q}\tilde{P} x\in H^{1/2} .
\end{equation*}
By regularity, $y\in H^{1/2}$
and hence $x=\tilde{P}x+y\in H^{1/2}$.
\end{proof}

\begin{cor}[Symmetry and Stability]\label{cc:stabsym}
$\phantom{x}$
\begin{enumerate}
\item
The regularity relation on pairs
of $1/2$-smooth projections is symmetric.
\item
The regularity relation is stable under smoothing perturbations,
i.e. if $P_1,P_2,Q_1,Q_2$ are $1/2$-smooth projections in $H$
with $P_1-P_2$ and $Q_1-Q_2$ smoothing,
then $(P_1, Q_1)$ is regular if and only if $(P_2, Q_2)$ is regular.
\qed
\end{enumerate}
\end{cor}

We need stronger regularity conditions:
The pair $(P, Q)$ is called {\em strongly regular}
if both $(P, Q)$ and $(I-P, I-Q)$ are regular.

\begin{thm}\label{cc:mt}
Let $P$ and $Q$ be 1/2-smooth projections in $H$.
Then the following conditions are equivalent.
\begin{enumerate}
\item
The pair $(P, Q)$ is strongly regular.
\item
The operator
\begin{equation*}
   T = T(P,Q) := P - Q = P(I-Q) - (I-P)Q
\end{equation*}
satisfies half of the condition
\ref{lem:preim}.\ref{eq:smprim}, i.e.,
\begin{equation*}
   x\in H^{-1/2} ,\, \tilde Tx \in H^{1/2}
    \Longrightarrow x\in H^{1/2} .
\end{equation*}
\end{enumerate}
\end{thm}

\begin{proof}
Assume that the pair $(P,Q)$ is strongly regular.
Let $x\in H^{-1/2}$ with $\tilde{T}x\in H^{1/2}$.
Then $(I -\tilde{P})\tilde Px=0$
and $(I-\tilde{Q})\tilde Px=(I-\hat{Q})\tilde{T}x$
is in $H^{1/2}$.
Hence $\tilde Px\in H^{1/2}$,
by the regularity of $(I-P,I-Q)$.
A similar argument shows that $\tilde{Q}x\in H^{1/2}$,
Hence $x\in H^{1/2}$, by the regularity of $(P,Q)$.
The other direction is obvious.
\end{proof}

In order to link strong regularity to Fredholm properties
of suitable operators, as in \cite{BL2},
we have to require regularity of the adjoint projections,
too.

\begin{thm}\label{thm:fredproj}
Let $P$ and $Q$ be $1/2$-smooth projections in $H$.
Then the following conditions are equivalent:
\begin{enumerate}
\item
The pairs $(P, Q)$ and $(P^{*},Q^{*})$ are strongly regular.
\item
With $T = T(P,Q) = P -Q$ as before,
the operators $\hat T$ and $\widehat{T^{*}}$ are Fredholm
in $H^{1/2}$ with $\ind\hat T+\ind\widehat{T^{*}}=0$.
\end{enumerate}
If any of these conditions holds then both $\tilde T$
and $\widetilde{T^{*}}$ restrict to Fredholm operators
in each $H^{s}$, $|s|\le 1/2$,
with kernels independent of $s$.
\end{thm}

\begin{proof}
 From \tref{cc:mt} we know that the strong regularity
of the pairs $(P,Q)$ and $(P^{*}, Q^{*})$ is equivalent
to the condition \ref{lem:preim}.\ref{eq:smprim}
for $T$ and $T^*$.
By \lref{lem:preim}, this condition is equivalent
to Condition (2) of the theorem.
The Fredholm property of the restrictions
and the constancy of their kernels
follows from \cref{cor:sfred}.
\end{proof}

\begin{rems}
2) $(I-P^{*},I-Q^{*})$ is (strongly) regular
if and only if $(P_\gamma,Q_\gamma)$ 
is (strongly) regular.
\\ \noindent
1) If $P$ and $Q$ are orthogonal,
that is, $P=P^*$ and $Q=Q^*$,
then $(P,Q)$ is strongly regular if and only if $(P,Q)$
and $(P_\gamma,Q_\gamma)$ are regular.
\end{rems}

\begin{cor}\label{orthproj}
For any pair $P,Q$ of orthogonal $1/2$-smooth projections
in $H$, the following conditions are equivalent.
\begin{enumerate}
\item
The pairs $(P, Q)$ and $(P_{\gamma}, Q_{\gamma})$
are regular.
\item $\hat T$ is a Fredholm operator,
necessarily of index 0, in $H^{1/2}$.
\qed
\end{enumerate}
\end{cor}

With any projection $Q$ in $H$,
we associate the involution $J(Q) := I - 2Q$.

\begin{prop}\label{cc:cr}
If there is a representation $P=I-Q +R_1+R_2$
in $\mathcal L(H^{1/2})$,
where $R_2$ and $R_2^*$ are compact in $H^{1/2}$ and
\[
  \|J(Q)R_1\|_{H^{1/2}} ,
  \|R_1^{*}J(Q^{*})\|_{H^{1/2}} < 1 ,
\]
then $(P,Q)$ and $(P^*,Q^*)$ are strongly regular.
\end{prop}

\begin{proof}
We show that Condition 2 of \tref{thm:fredproj} holds.
We have
\begin{equation*}
   T = J(Q) + R_1+R_2
   = J(Q)(I + J(Q)R_1)+R_2
\end{equation*}
and, similarly,
\begin{equation*}
   T^* = (I + R_1^{*}J(Q^{*}))J(Q^{*})+R_2^{*} .
\end{equation*}
The bound on the norms now implies that both $\hat T$
and $\widehat{T^{*}}$ are Fredholm operators in $H^{1/2}$
of index $0$, hence the assertion.
\end{proof}

In \cite[Theorem 1.3]{BL2} a criterion for regularity is given
which uses only properties of $P$ and $Q$ in $H$,
without referring to other Sobolev spaces,
at the expense of introducing more conditions on $P$ and $Q$.
This result is a special case of our analysis
as we will show now.

\begin{lem}\label{lem:hfred}
Let $S$ be a $1/2$-smooth Fredholm operator in $H$
and denote by $K_{r(l)}$ the orthogonal projections
onto $\ker S$ and $\ker S^{*}$, respectively.
Then the following conditions are equivalent:
\begin{enumerate}
\item
$S$ admits a $1/2$-smooth parametrix $U\in\mathcal{L}(H)$
such that
\begin{equation*}
    US = I - K_{r}
    \quad\text{and}\quad
    SU = I - K_{l} .
\end{equation*}
\item
$S$ and $S^{*}$ restrict respectively extend to
Fredholm operators in each $H^{s},\,|s|\le 1/2$,
with index independent of $s$.
\end{enumerate}
\end{lem}

\begin{proof}
$(1)\Rightarrow (2)$.
If $U$ restricts to $H^{1/2}$
then $K_{r}=I-US$ and $K_{l}=I-SU$ as well.
Since both projections have finite rank
and since $H^{1/2}$ is dense in $H$,
it follows that both projections are actually smoothing.
Now (2) follows from \lref{lem:preim} and \cref{cor:sfred}.

\noindent $(2)\Rightarrow (1)$.
This follows from the explicit construction of $K_{r(l)}$
in the proof of \lref{lem:preim}.4.
\end{proof}

This lemma gives a useful criterion for linking
the regularity of a $1/2$-smooth projection $P$
to Fredholm properties of $T = P-Q_{>}$ in $H$,
provided that we can control the mapping properties
of parametrices.
To construct a parametrix $U$ satisfying
Condition (1) of \lref{lem:hfred},
we start with the polar decomposition $T = V|T|$ of $T$,
where
\[
  |T| 
  = (T^{*}T)^{1/2}, \quad
  V^{*}V
  = I-K_r, \quad
  VV^{*}
  = I-K_l .
\]
Now $0$ is an isolated point in $\spec(T^{*}T)$
if $T$ is a Fredholm operator,
hence $\{\Re z >0\}\cap\spec(T^{*}T)$
is a compact subset of $(0,\infty)$.
The function $f=f(z)=1/\sqrt z$ is holomorphic
in $\{\Re z> 0\}$.
Thus we can define the operator $|T|^{-1} := f(T^{*}T)$
by the Dunford-Taylor integral of $f$
along a simple closed curve in $\{\Re z> 0\}$
surrounding $\spec T^{*}T\setminus\{0\}$
(cf. \cite[p.225]{Yo}).
Then we have $|T| |T|^{-1} = I - K_{r}$,
which implies that
\[
  U := |T|^{-1}V^{*}
\]
satisfies $UT = I - K_{r}$ and $TU = I - K_{l}$.
Now it is apparent that this parametrix construction leads
to a $1/2$-smooth parametrix for all Fredholm operators
inside an operator algebra,
$\mathcal{A}\subset \mathcal{L}(H)$,
if $\mathcal{A}$ has the following properties:
\begin{enumerate}
\item
$\mathcal{A}$ is a *-algebra with identity,
\item
$\mathcal{A}$ admits holomorphic functional calculus,
i.e., is closed under forming Dunford-Taylor integrals,
\item
$\mathcal{A}$ is contained in the space
of $1/2$-smooth operators.
\end{enumerate}
We combine these facts in the following result
which generalizes Theorem 1.3 in \cite{BL2}.

\begin{thm}\label{thm:hreg}
Let $P$  and $Q$ be $1/2$-smooth projections in $H$
and assume that $P$ and $Q$ are contained in some
operator algebra $\mathcal{A}\subset\mathcal{L}(H)$
which satisfies the above properties.
Then the following conditions are equivalent:
\begin{enumerate}
\item
The pairs $(P,Q)$ and $(P^{*},Q^{*})$
are strongly regular.
\item
The operator $T:= P - Q$ is Fredholm in $H$.
\end{enumerate}
\end{thm}

The conditions imposed on the algebra $\mathcal{A}$
are not unnatural; e.g.,
they are satisfied for the algebra of pseudodifferential
operators of order zero on a compact manifold.

We now come back to Dirac systems and
study the more traditional boundary conditions
defined by projections in $H$.
Let $P$ be a $1/2$-smooth projection in $H$.
Then $P$ induces a continuous projection 
in $\check H$ iff $Q_\le\tilde PQ_>$ is smoothing.
In any case,
\begin{equation}\label{rpbp}
  B_P := \ker\tilde P\cap\check H
\end{equation}
is a closed subspace of $\check H$,
that is, a boundary condition in the sense of
Definition \ref{bcbdef}.
Furthermore, $\ker\hat P$ is a closed subspace of $\hat H$.
In their work,
Atiyah, Patodi, and Singer consider the boundary 
condition given by $P_{APS}:=Q_>$,
see (2.3) in \cite{APS}.

\begin{rem}\label{bcbpro}
Let  $P$ be a $1/2$-smooth projection in $H$
that induces a continuous projection $\check P$ in $\check H$.
Since $H^{1/2}$ is dense in $\check H$
and $\check P(H^{1/2})\subset H^{1/2}$,
$B_P\cap H^{1/2}$ is dense in $B_P=\ker\check P$.
Hence $B_P$ is equal to the closure
of $B_P\cap H^{1/2}$ in $\check H$.

Suppose there is an $x\in\im\check P\setminus H^{1/2}$
and set $B=\ker\check P\oplus\R x$,
a closed subspace of $\check H$.
If $z=y+\alpha x\in B$ is in $H^{1/2}$,
then also $Pz=\alpha x$, hence $\alpha=0$.
It follows that $H^{1/2}$ is not dense in $B$.
By what we just said,
$B$ is a boundary condition that is not realizable
as the boundary condition $B_R$ of a $1/2$-smooth 
projection $R$
that induces a continuous projection in $\check H$.
\end{rem}

The Dirac operators and domains
corresponding to the boundary condition $B_P$
posed by a $1/2$-smooth projection $P$ in $H$
will be denoted as above,
except that we substitute the subscript $P$ for $B_P$.

\begin{dfn}\label{rpdef}
We say that a projection $P:H\to H$ is {\em regular}
if it is $1/2$-smooth and $B_P$
is a regular boundary condition.
\end{dfn}

\begin{prop}\label{rpprop}
For a $1/2$-smooth projection $P$ in $H$,
the following are equivalent:
\begin{enumerate}
\item
$P$ is regular.
\item
$B_P=\ker\hat P$.
\item
For some or, equivalently, any $\Lambda\in\R$,
we have
\[
  x\in H^{-1/2} , \tilde Px=0 ,
  Q_{\le\Lambda}x\in H^{1/2}
  \Longrightarrow x\in H^{1/2} .
\]
\end{enumerate}
\end{prop}

\begin{proof}
The condition in (2) expresses
that $B_P\subset H^{1/2}$,
hence that $B_P$ is a regular boundary
condition, by \pref{bcbregeq}.3.
Since $\check H$ is equal to the direct sum
$H_{\le\Lambda}^{1/2}\oplus H_{>\Lambda}^{-1/2}$,
the condition in (3) is just another way of saying
that $B_P\subset H^{1/2}$.
\end{proof}

Part 3 of the preceding result is the regularity criterion
introduced in condition (4.6c) of \cite{BL2}.

We note that for a regular projection $P$ in $H$ 
with corresponding boundary condition 
$B_P=\ker\hat P$,
the adjoint boundary condition is given by
\begin{equation}\label{bcbap}
  (B_P)^a = \ker\tilde P_\gamma \cap \check H
  \quad\text{with}\quad
  P_\gamma := \gamma^*(I-P^*)\gamma .
\end{equation}
We say that $P$ is {\em elliptic}
if $P$ and $P_{\gamma}$ are regular.
Then 
\begin{equation}\label{bcbap2}
  (B_P)^a = \ker\hat P_\gamma = \gamma\im\widehat{P^*} . 
\end{equation}

\begin{cor}\label{bcpsr}
If $P$ is an elliptic orthogonal projection in $H$,
then $D_{P,\ext}$ is a Fredholm operator 
with extended index
\[
  \ind D_{P,\ext} 
  = \dim (\ker P\cap H_\ge) - \dim (\im P\cap H_<) .
  \qed
\]
\end{cor}

\newpage
\section{Dirac-Schr\"odinger systems}\label{sec:axioms}

\subsection{Dirac systems with Lipschitz coefficients}
In this section, we construct and describe a model
for the geometric operators we are interested in;
this model will be introduced axiomatically.

Let $H$ be  a separable complex Hilbert space.
For $t\in\R_{+}$,
let $\langle\cdot,\cdot\rangle_{t}$ be a family
of scalar products with norm $||\cdot||_{t}$
compatible with the Hilbert space structure of $H$.

\begin{axia}\label{ax:metric}
For all $T\in\R_+$, there is a constant $C_{T}$ such that
\[
   |\langle u, v\rangle_{r}- \langle u, v\rangle_{s}|
   \le C_{T}||u||_{t}||v||_{t} |r-s|
\]
for all $u,v\in H$ and $r,s,t\in [0,T]$.
\end{axia}

It would be equivalent to require the estimate
for $t=0$ only instead of requiring it
for arbitrary $t\in[0,T]$.

In the following we will write
$\la\sigma,\tau\ra$ for the function
$t\mapsto\la\sigma(t),\tau(t)\ra_t$,
and similarly for related expressions.

Our data define a Lipschitz Hilbert bundle $\mathcal{H}$
over $\R_+$ with fibers $H_t=(H,\la\cdot,\cdot\ra_t)$,
$t\in\R_+$.
Any bundle $\mathcal{H}=(H_{t})_{t\in\R_{+}}$
of Hilbert spaces which is (locally) Lipschitz
over $\R_{+}$ is isometric to such a model bundle.

For $t\in\R_{+}$, define a positive definite operator
$G_{t}\in\mathcal{L}(H)$ by
\begin{equation}
   \langle G_{t}u, v\rangle_{0} = \langle u, v\rangle_{t},
   \quad u,v\in H .
\end{equation}
The operators $G_t$ and $G_t^{-1}$
are locally Lipschitz functions of $t$ in $\mathcal L(H)$.
An easy application of \lref{lem:weakdiff}
gives the following result.

\begin{lem}\label{lem:liparg}
The operator function $G$ is  weakly differentiable almost
everywhere in $\R_{+}$ with symmetric derivative
$G'_{t}\in L^{\infty}_{\rm loc}(\R_{+},\mathcal{L}(H))$.
\end{lem}

More generally,
if $H_1$ and $H_2$ are separable Hilbert spaces,
then any function in $\Lip_{\loc}(\R_+,\mathcal{L}(H_1,H_2))$
is weakly differentiable almost everywhere, 
and the norm of the derivative is locally uniformly bounded.

Now we set
\begin{equation}\label{gamma}
      \Gamma := \frac12 G^{-1}_{t}G'_{t}
      \in L^{\infty}_{\rm loc}(\R_{+},\mathcal{L}(H)) .
\end{equation}
If $\partial_t$ denotes the derivative with respect to $t$,
$\partial_t\sigma=\sigma'$, then
\begin{equation}\label{connection}
      \partial := \big(\partial_{t} + \Gamma\big):
      \Lip_{\loc}(\R_{+},H) \to L^{\infty}_{\loc}(\R_{+},H)
\end{equation}
is  a continuous metric connection,
where {\em metric} means that
\begin{equation}\label{connectionmetric}
   \langle\sigma_1,\sigma_2\rangle'
   = \langle\partial\sigma_1,\sigma_2\rangle
   + \langle\sigma_1,\partial\sigma_2\rangle ,
\end{equation}
for all $\sigma_1,\sigma_2\in\Lip_{\loc}(\R_+,H)$.

\begin{rem}\label{connectionmet2}
Any other continuous metric connection
\[
   \tilde{\partial}:
   \Lip_{\loc}(\R_{+},H) \to L^{\infty}_{\rm loc}(\R_{+},H)
\]
is of the form
$\tilde{\partial} = \partial + \tilde{\Gamma}$,
where
$\tilde{\Gamma}\in L^{\infty}_{\loc}(\R_{+},\mathcal{L}(H))$
takes values in the space of skew-Hermitian operators.
\end{rem}

\begin{axia}\label{ax:A}
There is a family $\mathcal A$
of self-adjoint operators $A_{t}$ on $H_t$, $t\in\R_{+}$, 
with common domain $H_{A}$ and graph norm $||\cdot||_{A_{t}}$
such that
\begin{enumerate}
\item
with respect to the graph norm $||\cdot||_{A_{0}}$ on $H_{A}$,
\\ the embedding $H_{A}\to H$ is compact;
\item
for all $T\in\R_{+}$, there is a constant $C_{T}$ such that
\[
    |\langle A_{r}u, v\rangle_{r} 
    - \langle A_{s}u,v\rangle_{s}|
    \le C_{T}||u||_{A_{t}}||v||_{t}|r-s|
\]
for all $u\in H_{A},v\in H$, and $r,s,t\in[0,T]$.
\end{enumerate}
\end{axia}

As above in Axiom \ref{ax:metric},
it would be equivalent to require the estimate
for $t=0$ only instead of requiring it
for arbitrary $t\in[0,T]$.

\begin{rem}\label{problem}
It would be tempting to use the metric connection $\partial$
to identify $\mathcal{H}$ with $\R_{+}\times H_0$.
But this parallel transport may not preserve $H_{A}$
if $\Gamma$ does not,
and this happens indeed in important examples.
\end{rem}

A pair $e := (\mathcal{H},\mathcal A)$ satisfying
Axioms I and II will be called an {\em evolution system}.
To any evolution system $e$ we can naturally associate
a family of constant coefficient system $e^{t}$,
$t\in\R_+$, defined by
\begin{equation}\label{etoet}
  e^{t} := (H_t,A_{t}) .
\end{equation}
For any evolution system $e$, we introduce the Hilbert
space $L^{2}(\mathcal{H})$ as completion of the space
$\mathcal{L}_{c}(e^{0})$ under the norm
\begin{equation}\label{l2curlyh}
  ||\sigma||^{2}_{L^{2}(\mathcal{H})}
  := \int_{0}^{\infty}||\sigma||^{2}_t dt .
\end{equation}
Then we can form the linear operator
\begin{equation}\label{evolution}
   L := \partial + A:
   \mathcal{L}_{c}(e^{0}) \to L^{2}(\mathcal{H}),
\end{equation}
which we call the {\em evolution operator} associated to $e$.
Note that the domain of $L$ only depends
on the constant coefficient system $e^{0}$.

The evolution operator $L$ introduced above is not symmetric
on the dense subspace $\mathcal{L}_{0,c}(e^{0})$
of $L^{2}(\mathcal{H})$.
A modification as in the case of constant coefficients
leads to a symmetric operator.

\begin{axia}\label{ax:gamma}
There is a section
\[
  \gamma\in\Lip_{\loc}\big(\R_+,\mathcal{L}(H)\big)
  \cap L^{\infty}_{\loc}\big(\R_+,\mathcal{L}(H_{A})\big),
\]
such that the following relations hold:
\begin{alignat*}{2}
    -\gamma_t = \gamma_t^* &= \gamma_{t}^{-1} \quad
    &&\mbox{on $H_t$} ,
    \tag{1} \\
    A_t\gamma_t + \gamma_t A_t &= 0
    &&\mbox{on $H_{A}$} ,
    \tag{2} \\
    \lbrack\partial,\gamma\rbrack &= 0
    &&\mbox{on $\Lip_{\loc}(\R_{+},H)$} .
    \tag{3}
\end{alignat*}
\end{axia}

\noindent
Note that $\gamma\mathcal{L}_{c}(e^{0}) \subset
\mathcal{L}_{c}(e^{0})$, by assumption.

A triple $d := (\mathcal H,\mathcal A,\gamma)$
as above satisfying Axioms \ref{ax:metric}--\ref{ax:gamma},
is called a {\em Dirac system}.
Now we are ready to introduce our first model operator,
the {\em Dirac operator}
\begin{equation}\label{diracsys}
    D := \gamma (\partial + A):
      \mathcal L_{\loc}(e^0) \to L^{\infty}_{\rm loc}(\R_+,H),
\end{equation}
associated to the Dirac system $(\mathcal H,\mathcal A,\gamma)$.

For later purposes it is important to note that, pointwise,
\begin{equation}\label{normest}
   ||D\sigma|| = ||L\sigma||
\end{equation}
for all $\sigma\in\mathcal L_{\loc}(e^0)$, so that estimates
for the usual norms of $L\sigma$ also hold for $D\sigma$.

The restriction $D_{0,c}$ of $D$ to the domain
$\mathcal{D}_{0,c}:=\mathcal{L}_{0,c}(e^{0})$ is symmetric;
we denote by $D_{\min}$, with domain $\mathcal{D}_{\min}$,
the closure of $D_{0,c}$ in $L^2(\mathcal H)$,
and by $D_{\max}$, with domain $\mathcal{D}_{\max}$,
the adjoint operator.
In order to define self-adjoint extensions of $D_{\min}$,
we will introduce boundary conditions as in Chapter \ref{sec:cc}.
Again, this approach is based on integration by parts
and the boundary form $\omega$:
\eqref{eq:intpart0} and \eqref{eq:intpart1} translate literally
in view of the following computation, valid for all
$\sigma_{1},\sigma_{2}\in\mathcal{L}_{c}(e^{0})$,
\begin{equation}\label{eq:intpart0gen}
  \langle\gamma\sigma_{1},\sigma_{2}\rangle'
  = \langle\gamma\partial\sigma_{1},\sigma_{2}\rangle
     + \langle\gamma\sigma_{1},\partial\sigma_{2}\rangle
  = \langle D\sigma_{1},\sigma_{2}\rangle
     - \langle\sigma_{1},D\sigma_{2}\rangle ,
\end{equation}
which is an easy consequence of our axioms;
therefore, we also get
\begin{equation}
     (D\sigma_{1},\sigma_{2}) - (\sigma_{1},D\sigma_{2})
     = \omega (\sigma_{1}(0),\sigma_{2}(0)).
     \label{eq:intpart1gen}
\end{equation}
In particular, we have
$\mathcal{L}_{c}(e^{0})\subset\mathcal{D}_{\max}$.

\subsection{Comparison with constant coefficients}
Let $d$ be a Dirac system with Lipschitz coefficients.
Our strategy in dealing with $d$ aims at some kind of
comparison with constant Dirac systems,
where we have substantial control over the solution theory.
Any such attempt meets with two difficulties,
firstly that we lack any a priori control
on the domain of the maximal operator $D_{\max}=(D_{0,c})^*$
and secondly, that the domain of the adjoint operator
to $A_{t}$ in $H_{0}$ varies with $t$.

For any $t\ge0$,
we introduce the Dirac system
\begin{equation}\label{constantt}
  d^t = (H_t,A_t,\gamma_t)
\end{equation}
with constant coefficients and the Dirac system $d^{ct}$
with coefficients
\begin{align}\label{freezeget}
  \begin{matrix}
  H^{ct}_s = (H,\la.,.\ra_s) , &
  A^{ct}_s = A_s , &
  \gamma^{ct}_s = \gamma_s &
  \text{for $s\le t$} , \\
  H^{ct}_s = (H,\la.,.\ra_t) , &
  A^{ct}_s = A_t , &
  \gamma^{ct}_s = \gamma_t &
  \text{for $s\ge t$} .
  \end{matrix}
\end{align}
Objects associated to $d^t$ and $d^{ct}$
will be decorated with a superscript $t$ and $ct$,
respectively.
We think of $d^{ct}$ as a kind of interpolation
between $d^0=d^{c0}$ and $d$.

\begin{thm}\label{freezedom}
The Dirac systems $d^{ct}$ compare with $d^0$ as follows:
\begin{enumerate}
\item
For all $t\ge0$,
we have ${\mathcal D}^{ct}_{\min}=\mathcal D^0_{\min}$
and ${\mathcal D}^{ct}_{\max}=\mathcal D^0_{\max}$.
\item
For all $T\ge0$, there is a constant $C_T$
such that
\[
  C_T^{-1}||\cdot||_{{D}^{0}_{\max}}
  \le ||\cdot||_{{D}^{ct}_{\max}}
  \le C_T ||\cdot||_{{D}^{0}_{\max}}
\]
for all $t\in[0,T]$.
\end{enumerate}
\end{thm}

The proof of \tref{freezedom} will be given below.
In preparation, we will study the operators $G^{ct}D^{ct}$,
which are symmetric in $L^2(\mathcal H^0)$
with domain $\mathcal L_c(e^0)$.

We start with some estimates.
Axioms \ref{ax:metric} and \ref{ax:gamma} imply that,
for any $t\ge0$, there is a constant $C_T$ such that,
for all $r,s\in[0,T]$,
\begin{align}
  ||G_s\gamma_s\gamma_r^{-1}-G_r||_0
  &\le C_T |r-s| ,
  \label{estmetgain} \\
  ||G_r\gamma_r\Gamma_r||_0
  &\le C_T ,
  \label{estmetgaGa} \\
  ||G_r\gamma_r||_0
  &\le C_T .
  \label{estmetgam}
\end{align}
We will also need estimates on the operators 
$A_t$.
 From Axiom \ref{ax:A} we get,
for $0\le s,t \le T$ and $x\in H_A$,
\begin{equation*}
\begin{split}
  ||A_sx||_s^2
  &\le C_T ||x||_{A_t} ||A_sx||_t
        + \la A_tx,A_sx \ra_t \\
  &\le C_T ||x||_{A_t} ||A_sx||_t
        + ||A_tx||_t ||A_sx||_t \\
  &\le C_T ||x||_{A_t} ||A_sx||_s ,
\end{split}
\end{equation*}
where the constant $C_T$ may change from line to line.
Therefore
\begin{equation}\label{estas0}
  ||\cdot||_{A_s} \le C_T ||\cdot||_{A_t}
\end{equation}
for all $s,t\in[0,T]$.
In other words,
the graph norms $||\cdot||_{A_t}$ are locally uniformly
equivalent.
For all $r,s,t\in[0,T]$ and $x\in H_A$, we also have
\begin{align*}
  ||A_rx-A_sx||_t^2
  &\le C_T ||A_rx-A_sx||_r^2 \\
  &= C_T \cdot \la A_rx-A_sx,A_rx-A_sx \ra_r \\
  &\qquad + C_T \cdot
          \la A_sx-A_sx,A_rx-A_sx \ra_s \\
  &\le C_T |r-s| ||x||_{A_t}||A_rx-A_sx||_t \\
  &\qquad + C_T \big| \la A_sx,A_rx-A_sx\ra_r
   - \la A_sx,A_rx-A_sx\ra_s \big| \\
  &\le C_T |r-s|
    (||x||_{A_t}+||A_sx||_t) ||A_rx-A_sx||_t \\
  &\le C_T |r-s| \cdot ||x||_{A_t} ||A_rx-A_sx||_t ,
\end{align*}
where we use Axiom \ref{ax:metric} and \eqref{estas0}
in the last two inequalities.
Therefore
\begin{equation}\label{estars}
  ||A_rx-A_sx||_t
  \le C_T |r-s| \cdot ||x||_{A_t}
\end{equation}
for all $0\le r,s,t\le T$ and $x\in H_A$.

The main estimate we need is of Kato-Rellich type:

\begin{lem}\label{katorell}
Given $T\ge0$, there is a constant $C_T$ such that,
for all $r\le s$ in $[0,T]$
and $\sigma\in\mathcal{L}_{c}(e^{0})$,
\begin{multline*}
  ||G^{cr}D^{cr}\sigma
    - G^{cs}D^{cs}\sigma||_{L^{2}(\mathcal{H}^{0})} \\
  \le C_T ||\sigma||_{L^{2}(\mathcal{H}^{0})}
    +  C_T |r-s| \cdot
    ||G^{cr}D^{cr}\sigma||_{L^{2}(\mathcal{H}^{0})} .
\end{multline*}
\end{lem}

\begin{proof}
We start by comparing the coefficients of the two operators
$G^{cr}D^{cr}$ and $G^{cs}D^{cs}$.
On $[0,r]$, they coincide.
At $t\in(r,s]$, we have
\begin{align*}
  (G^{cs}D^{cs})|_t
  &= G_t\gamma_t(\partial_t+\Gamma_t+A_t)
  \notag \\
  &= G_t\gamma_t\gamma_r^{-1} D^{cr}
     + G_t\gamma_t\Gamma_t
     + G_t\gamma_t (A_t-A_r) .
\end{align*}
At $t\in[s,\infty)$, we have
\begin{equation*}
  (G^{cs}D^{cs})|_t
  = G_s\gamma_s(\partial_t+A_s)
  = G_s\gamma_s\gamma_r^{-1} D^{cr}
     + G_s\gamma_s (A_s-A_r) .
\end{equation*}
Let $\sigma\in\mathcal{L}_{c}(e^{0})$.
Then $G^{cr}D^{cr}\sigma$ and $G^{cs}D^{cs}\sigma$
coincide on $[0,r]$.
Using \eqref{estmetgain}, \eqref{estmetgaGa},
and \eqref{estmetgam}, we get
\begin{equation}\label{estgen}
\begin{split}
  ||G^{cr}D^{cr}\sigma
   &- G^{cs}D^{cs}\sigma||_{L^{2}(\mathcal{H}^{0})}
  \le C_T |r-s| \cdot
  ||D^{cr}\sigma||_{L^{2}(\mathcal{H}^{0})} \\
  &\quad
  + C_T ||\sigma||_{L^{2}(\mathcal{H}^{0})}
  + C_T
  ||(A^{cr}-A^{cs})\sigma||_{L^{2}(\mathcal{H}^{0})} .
\end{split}
\end{equation}
By Axiom \ref{ax:metric},
\[
  ||D^{cr}\sigma||_{L^{2}(\mathcal{H}^{0})} \le C_T
  ||G^{cr}D^{cr}\sigma||_{L^{2}(\mathcal{H}^{0})} ,
\]
hence the first two terms on the right
in \eqref{estgen} are under control as desired.
It remains to get a good upper bound for
$||(A^{cr}- A^{cs})\sigma||_{L^{2}(\mathcal{H}^{0})}$.
By \eqref{estas0} and \eqref{estars},
\begin{equation*}
  ||(A^{cr}- A^{cs})\sigma||_{L^{2}(\mathcal{H}^{0})}
  \le
  C_T|r-s| \cdot ||\sigma||_{L^2(\mathcal{H}^{0})}
   + C_T ||\varphi A_r\sigma||_{L^2(\R_+,H_r)} ,
\end{equation*}
where $\varphi(t)=\inf(t-r,s-r)$ for $t\ge r$
and $\varphi(t)=0$ for $t\le r$.
It remains to estimate
the second term on the right of this inequality.
We compute
\[
  ||(\varphi\sigma)'+\varphi A_r\sigma||_r^2
  = ||(\varphi\sigma)'||_r^2
    + ||\varphi A_r\sigma||_r^2
    + \la \varphi A_r\sigma,\varphi\sigma \ra_r' .
\]
Now $\varphi\sigma\in\mathcal L_c(e^0)$
vanishes at $0$, hence
\[
  ||D^r(\varphi\sigma)||_{L^2(\R_+,H_r)}^2
  = ||(\varphi\sigma)'||_{L^2(\R_+,H_r)}^2
    + ||\varphi A_r\sigma||_{L^2(\R_+,H_r)}^2 .
\]
Since $D^r(\varphi\sigma)=\varphi'\gamma^r\sigma
       + \varphi D^r\sigma$,
we conclude
\begin{align*}
  ||\varphi A_r\sigma||_{L^2(\R_+,H_r)}
  &\le C_T \cdot ||\sigma||_{L^2(\mathcal H^0)}
       + ||D^r(\varphi\sigma)||_{L^2(\R_+,H_r)} \\
  &\le
    C_T \cdot ||\sigma||_{L^2(\mathcal H^0)}
   + |s-r| ||D^{cr}\sigma||_{L^2(\R_+,H_r)} \\
  &\le C_T \cdot ||\sigma||_{L^2(\mathcal H^0)}
   + C_T |s-r| \cdot
   ||G^{cr}D^{cr}\sigma||_{L^2(\mathcal H_0)} .
  \qedhere
\end{align*}
\end{proof}

\begin{proof}[Proof of \tref{freezedom}]
We note first that the Hilbert spaces
$L^2(\mathcal H^{ct})$ and $L^2(\mathcal H^0)$
coincide as vector spaces of
(equivalence classes of) maps.
The operators $D^{ct}$ and $G^{ct}D^{ct}$ have 
the same minimal and maximal domains.
Hence we may as well consider the family
of operators $G^{ct}D^{ct}$ on $L^2(\mathcal H^0)$.
We introduce operators
\begin{equation}\label{minmax}
  S^t =
  \begin{pmatrix}
  0 & G^{ct}D^{ct} \\ G^{ct}D^{ct}_{0,c} & 0
  \end{pmatrix}
  \quad\text{and}\quad
  T^t =
  \begin{pmatrix}
  0 & G^{ct}D^{ct}_{\max} \\ G^{ct}D^{ct}_{\min} & 0
  \end{pmatrix}
\end{equation}
in $L^2(\mathcal H^0)\oplus L^2(\mathcal H^0)$
with domain
$\mathcal L_{0,c}(e^0)\oplus\mathcal L_{c}(e^0)$
and
$\mathcal D^{ct}_{\min}\oplus\mathcal D^{ct}_{\max}$,
respectively.
We note that $S^t$ is symmetric
and that $T^t$ is self-adjoint
with $S^t\subset T^t$.

Fix $T\ge0$ and assume that, for some $r\in[0,T]$,
the closure of $S^r$ is equal to $T^r$
with domain
$\mathcal D^0_{\min}\oplus\mathcal D^0_{\max}$.
By the results of the first section,
this holds for $r=0$.
By the Kato-Rellich Theorem,
see Theorem V.4.4 in \cite{Ka} and \lref{katorell},
we get that the closure of $S^s$ is self-adjoint
with domain
$\mathcal D^0_{\min}\oplus\mathcal D^0_{\max}$
for all $s\ge r$ in $[0,T]$ with $(s-r)C<1/2$,
where $C=C_T$ is the constant from \lref{katorell}.
Since $S^s\subset T^s$ and $T^s$ is self-adjoint,
we conclude that the closure of $S^s$ is equal
to $T^s$ for all such $s$.
By the connectedness of $[0,T]$, 
we get that the closure of $S^r$
is equal to $T^r$ with domain
$\mathcal D^0_{\min}\oplus\mathcal D^0_{\max}$
for all $r\in[0,T]$.
This proves the first assertion.

As for the proof of the second assertion,
we note that \lref{katorell} implies that $D^{cr}_{\max}$
and $D^{cs}_{\max}$ have equivalent graph norms
on their common domain $\mathcal D_{\max}^{0}$
as soon as $|r-s|C<1$.
Again by the connectedness of $[0,T]$,
the graph norm of $D^{ct}_{\max}$ is equivalent
to the one of $D^0_{\max}$.
Hence there is a constant as claimed.
\end{proof}

For applications it is useful to pass to a somewhat 
more general class of systems and operators.

\begin{dfn}\label{dsdef}
A {\em Dirac-Schr\"odinger system} is a pair $(d,V)$
consisting of a Dirac system $d$ with Lipschitz coefficients
and a {\em potential}
$V\in L^{\infty}_{\loc}(\mathcal L(\mathcal H))$
with $V=V^*$.
The associated {\em Dirac-Schr\"odinger operator}
is given by
\[
  D:= D^d+V : \mathcal L_{\loc}(e^0) 
      \to L^\infty_{\loc}(\mathcal H) ,
\]
where $D^d$ denotes the Dirac operator of $d$.
\end{dfn}

\begin{rem}
It is not really necessary to assume that the potential
is Hermitian, $V=V^*$.
However, assuming $V=V^*$ keeps the notation a bit simpler.
For most purposes, passing to the Dirac-Schr\"odinger system
with operator
\[
  \begin{pmatrix} 0 & D^d + V^* \\ D^d + V & 0 
  \end{pmatrix}
\]
reduces the general case to the case where $V$ is Hermitian.
\end{rem}

In what follows, $D$ is the Dirac-Schr\"odinger operator
associated to a Dirac-Schr\"odinger system $(d,V)$.
 From \eqref{eq:intpart1gen} we get
\begin{equation}\label{dsintpar}
  (D\sigma_1,\sigma_2) - (\sigma_1,D\sigma_2)
  = \omega (\sigma_1(0),\sigma_2(0)) ,
\end{equation}
where $\sigma_1,\sigma_2\in\mathcal L_c(e^0)$.
Therefore the restriction $D_{0,c}$ of $D$
to the domain $\mathcal D_{0,c}$ is symmetric.
We denote by $D_{\min}$, with domain $\mathcal D_{\min}$,
the closure of $D_{0,c}$ in $L^2(\mathcal H)$
and by $D_{\max}:=(D_{0,c})^*$, 
with domain $\mathcal D_{\max}$,
the adjoint operator of $D_{0,c}$ in $L^2(\mathcal H)$.

We let $D^0$ be the Dirac operator
associated to the constant coefficient Dirac system $d^0$
and $\mathcal D^0_{\max}$ be its domain.
The following result is crucial.

\begin{thm}\label{dcomp}
If $\sigma\in L^2(\mathcal H)$ has compact support,
then $\sigma\in\mathcal D_{\max}$
if and only if $\sigma\in\mathcal D_{\max}^{0}$.
\end{thm}

\begin{proof}
Suppose that $\sigma\in L^2(\mathcal H)$
has compact support in $[0,R]$.
Since $V\in L^{\infty}_{\loc}(\mathcal L(\mathcal H))$,
$V$ is uniformly bounded on $[0,R]$,
and hence we may assume that $V=0$.
Choose $T>R$.
For any $t\in (R,T)$, the coefficients of
$D$ and $D^{ct}$ coincide on $[0,R]\subset[0,t]$,
compare \eqref{freezeget}.
Hence $\sigma\in\mathcal D_{\max}$
if and only if $\sigma\in\mathcal D_{\max}^{ct}$,
and from \tref{freezedom},
$\mathcal D_{\max}^{ct}=\mathcal D_{\max}^{0}$.
\end{proof}

\begin{prop}[Regularity]\label{thm:diprops}
The maximal domain $\mathcal{D}_{\max}$ satisfies:
\begin{enumerate}
\item
$\mathcal{L}_{c}(e^{0})$
is dense in $\mathcal{D}_{\max}$.
\item
$\sigma\in\mathcal{D}_{\max}$ is in $H^1_{\loc}(e^{0})$
if and only if $\sigma(0)\in H^{1/2}$.
\item
$\mathcal D_{\max}\subset
 C(\R_{+},\check H)\cap C((0,\infty),H^{1/2})$.
\item
The restriction map on $\mathcal{L}_{c}(e^{0})$
extends to a continuous \\ surjective map
$\mathcal{R}:\mathcal{D}_{\max}\to\check H$
and $\mathcal D_{\min}=\mathcal R^{-1}(0)$.
\item
For $\sigma_{1},\sigma_2\in\mathcal{D}_{\max}$, we have
\begin{equation*}
  ( D_{\max}\sigma_{1},\sigma_{2} )_{L^2(\mathcal H)}
  - ( \sigma_{1},D_{\max}\sigma_2 )_{L^2(\mathcal H)}
  = \omega\big(\sigma_{1}(0),\sigma_{2}(0)\big) .
\end{equation*}
\end{enumerate}
\end{prop}

\begin{proof}
The first assertion follows from \lref{lem:regprop}.1
and \tref{dcomp}.
As for the proof of the second and third assertion,
multiply $\sigma\in\mathcal D_{\max}$ by
a Lipschitz cutoff function $\chi$ which is equal to $1$
on some interval $[0,R]$ and equal to $0$ after $2R$.
Then $\chi\sigma$ is in $\mathcal D^0_{\max}$,
by \tref{dcomp}, and $\chi\sigma$ has the asserted
regularity properties, by \lref{lem:regprop}.
By \tref{dcomp},
multiplication by $\chi$ defines a continuous operator
from $\mathcal{D}_{\max}$ to $\mathcal D^0_{\max}$,
hence the fourth assertion is immediate from \pref{dmaxbou}.
By (1) it is enough to check the last assertion
for $\sigma_{1},\sigma_{2}\in\mathcal{L}_{c}(e^{0})$.
This case was already observed in \eqref{dsintpar}.
\end{proof}

\subsection{Boundary conditions and Fredholm properties}
\label{bocogen}
We now turn to the description of closed extensions of $D$,
following closely the outline given in Section \ref{boco};
most proofs carry over easily via the link given by \tref{dcomp}.
In what follows, we fix a Dirac-Schr\"odinger system $(d,V)$
and define the Sobolev spaces $H^s$ and $\check H$
with respect to $A_0$ as in Section \ref{soba}.

As before, a boundary condition is a closed linear subspace
$B\subset\check H$.
Associated to a boundary condition $B$,
we consider extensions of $D_{0,c}$ as
in Section \ref{boco}:
\begin{align}
  \mathcal L_{B,c}
  &:= \{\sigma\in\mathcal L_c(e^0) \,:\, \sigma(0)\in B \} ,
  \label{bcbcgen} \\
  D_{B,c} &:= D|\mathcal L_{B,c} ;
  \notag \\
  \mathcal{D}_{B}
  &:= \{\sigma\in\mathcal D_{\max}\cap H^1_{\loc}(e^0)
  \,:\, \sigma(0)\in B \} ,
  \label{bcb1gen} \\
  D_{B} &:= D_{\max} | \mathcal{D}_{B} ;
  \notag \\
  \mathcal{D}_{B,\max}
  &:= \{\sigma\in\mathcal D_{\max} \,:\, \sigma(0) \in B \} ,
  \label{bcbmaxgen} \\
  D_{B,\max} &:= D_{\max} |\mathcal{D}_{B,\max} .
  \notag
\end{align}
As before, since the restriction map
$\mathcal R:\mathcal D_{\max}\to\check H$ is continuous
and $B$ is closed in $\check H$,
$D_{B,\max}$ is a closed operator.
Moreover, any closed extension of $D_{0,c}$ with domain
contained in $\mathcal D_{\max}$ is of this form.

\begin{rem}\label{bcbadclos}
The same formulas for the adjoint operators
and the closures as in \eqref{bcbcad}--\eqref{bcbmaxclos}
continue to hold and for the same reasons.
We do not repeat them here.
\end{rem}

As before,
we say that a boundary condition $B\subset\check H$
is {\em regular} if $D_{B,\max}=D_B$.
It is immediate from \pref{thm:diprops}.2. that
\begin{enumerate}
\item
in the case of constant coefficients with potential $V=0$,
the present definition coincides with the one 
in Section \ref{boco};
\item
a boundary condition $B$ is regular relative to $(d,V)$
if and only if it is regular relative to $d^{0}$.
\end{enumerate}
As in Section \ref{boco},
we say that a boundary condition $B$ is
{\em elliptic} if $B$ and $B^a$ are regular.

In the case of constant coefficients with potential $V=0$,
$D_B$ is not a Fredholm operator whenever $\ker A_0\ne0$,
even if $B$ is elliptic.
However, we may look for an analogue of the space $W$
which worked so nicely in the constant coefficient case.
 From the continuity of $\mathcal{R}$, established in
\tref{thm:diprops}.4 we get that there is a constant $C$
such that
\begin{equation}\label{eq:west2gen}
    ||\sigma||^{2}_{W}
    := ||\sigma(0)||_{\check H}^{2}
    + ||D_{\max}\sigma||_{L^{2}(\mathcal{H})}^{2}
    \le C||\sigma||_{\mathcal{D}_{\max}}^{2}.
\end{equation}
for all $\sigma\in\mathcal D_{\max}$.
The converse of \eqref{eq:west2gen} is not available
in general, as we know,
but a localized version may hold.
This requires the inequality \eqref{eq:west0}
which we now introduce as an additional axiom.

\begin{axia}\label{ax:nonpar}
For each $T>0$ there is a constant $C_{T}$ such that
\begin{equation*}
  ||\sigma||_{L^{2}([0,T],\mathcal H)}
  \le C_{T} ||\sigma||_{W}
  \quad\text{for all $\sigma\in\mathcal{L}_{c}(e^{0})$} .
\end{equation*}
\end{axia}

Following G. Carron (cf. the introduction to \cite{Ca2})
we will call a Dirac-Schr\"odinger system $(d,V)$ satisfying
Axiom \ref{ax:nonpar} {\em non-parabolic} (at infinity).
We say that a Dirac-Schr\"odinger system $(d,V)$ is 
{\em of Fredholm type}, if there is a constant $C$ such that
\begin{equation}\label{freddirac}
  ||\sigma||_{L^{2}(\mathcal H)}
  \le C ||\sigma||_{W}
  \quad\text{for all $\sigma\in\mathcal{L}_{c}(e^{0})$} .
\end{equation}
If $(d,V)$ is non-parabolic, then $(d,V)$ is of Fredholm type
if and only if, for some $\psi\in\Lip_c(\R_+)$
which is equal to $1$ near $t=0$,
\begin{equation}\label{freddirac2}
  ||(1-\psi)\sigma||_{L^{2}(\mathcal H)}
  \le C_{\psi} ||\sigma||_{W}
  \quad\text{for all $\sigma\in\mathcal{L}_{c}(e^{0})$} .
\end{equation}
In the geometric setting considered by Carron, it is
enough to work with smooth sections supported near infinity,
hence the space $\check H$ does not enter his discussion.
However, the two formulations of non-parabolicity
here and there are equivalent in the following sense.

\begin{lem}\label{nonpareq}
The inequality of Axiom \ref{ax:nonpar} holds
for all $\sigma\in\mathcal{L}_{c}(e^{0})$
if it holds for all $\sigma\in\mathcal{L}_{0,c}(e^{0})$.
\end{lem}

\begin{proof}
Choose $\psi\in \Lip_{c}(\R_{+})$ with $\psi(0) = 1$.
Let $D^0$ be the Dirac operator
and $\mathcal{E}^{0}$ be the extension operator
for $d^{0}$, see \eqref{eq:extdef}.
Let $\sigma\in\mathcal{L}_{c}(e^{0})$ and set
\[
  \sigma_{0} := \psi\mathcal{E}^{0}\sigma(0)
  \quad\text{and}\quad
  \sigma_{1} := \sigma - \sigma_{0} .
\]
Since $\sigma(0)\in H_{A}$,
we have $\sigma_0\in \mathcal{L}_{c}(e^{0})$;
hence $\sigma_{1}\in\mathcal{L}_{0,c}(e^{0})$.
Now we can estimate, using the assumption,
\lref{cc:ne}, and \tref{dcomp},
\begin{align*}
    ||\sigma||_{L^{2}([0,T],H)}
    &\le  ||\sigma_{1}||_{L^{2}([0,T],H)}
    +||\sigma_{0}||_{L^{2}([0,T],H)} \\
     \le C_{T,\psi}
     &(||D_{\max}\sigma_{1}||_{L^{2}(\mathcal{H})}
     +  ||\sigma(0)||_{-1/2})\\
     \le C_{T,\psi}
     &(||D_{\max}\sigma||_{L^{2}(\mathcal{H})}
     + ||D_{\max}\psi\mathcal{E}^{0}
         \sigma(0)||_{L^{2}(\mathcal{H})}
     + ||\sigma(0)||_{-1/2}) \\
     \le C_{T,\psi}
     &(||D_{\max}\sigma||_{L^{2}(\mathcal{H})}
     + ||D^0_{\max}\psi\mathcal{E}^{0}
         \sigma(0)||_{L^{2}(\mathcal{H})}
     + ||\sigma(0)||_{-1/2}) \\
     = C_{T,\psi}
     &( ||D_{\max}\sigma||_{L^{2}(\mathcal{H})} \\
     &\phantom{xx}
      + ||(A_{0} - |A_{0}| - Q_{0})\psi\mathcal{E}^{0}
         \sigma(0)||_{L^{2}(\mathcal{H})}
      + ||\sigma(0)||_{-1/2} ) \\
     \le C_{T,\psi}
     &(||D_{\max}\sigma||_{L^{2}(\mathcal{H})}
     + ||\sigma(0)||_{\check H}) ,
\end{align*}
where we allow the constant $C_{T,\psi}$
to change from line to line.
\end{proof}

As a first implication of non-parabolicity we note
that the seminorm $||\cdot||_{W}$,
as defined in \eqref{eq:west2gen},
is actually a norm on $\mathcal{D}_{\max}$.
Thus we can introduce again the space $W$ as
the completion of $\mathcal D_{\max}$ under this norm.
Since $\mathcal{L}_{c}(e^{0})$ is dense
in $\mathcal D_{\max}$ with respect to the graph norm
of $D_{\max}$, $\mathcal{L}_{c}(e^{0})$
is dense in $W$ with respect to the $W$-norm.

\begin{lem}\label{lem:west3gen}
If $(d,V)$ is a non-parabolic Dirac-Schr\"odinger system,
then we have:
\begin{enumerate}
\item
The restriction map $\mathcal{R}$ and $D_{\max}$
extend to continuous maps $\mathcal{R}_{\ext}$
and $D_{\ext}$ on $W$, respectively;
$\mathcal{R}_{\ext}$ induces an isometry
from $\ker D_{\ext}$ into $\check H$.
\item
If $\psi\in \Lip_{c}(\R_{+})$ and $\sigma\in W$,
then $\psi\sigma\in\mathcal{D}_{\max}\subset W$.
Moreover, there is a constant $C_{\psi}$ such that
\[
  ||\psi\sigma||_{D_{\max}}
  \le C_\psi ||\sigma||_W .
\]
In particular, $W$ can be viewed as
a space of locally integrable functions
and $W\cap L^{2}(\mathcal{H}) = \mathcal{D}_{\max}$.
\item
$W=\mathcal{D}_{\max}$ if and only if $(d,V)$ is a
Dirac-Schr\"odinger system of Fredholm type; that is,
there is a constant $C$ such that
\begin{equation*}
  ||\sigma||_{L^{2}(\mathcal{H})}
  \le C ||\sigma||_W
  \quad\text{for all $\sigma\in\mathcal{L}_{c}(e^{0})$} .
\end{equation*}
\end{enumerate}
\end{lem}

\begin{proof}
(1) and (3) are immediate from the definition of $W$.
As for (2), we note that, by non-parabolicity,
there is a constant $C_\psi$ such that
\begin{equation*}
    ||\psi\sigma||_{\mathcal{D}_{\max}}
    \le C_{\psi}||\sigma||_{W}
\end{equation*}
for all $\sigma\in\mathcal{L}_{c}(e^{0})$,
hence for all $\sigma\in W$ by the density
of $\mathcal{L}_{c}(e^{0})$.

Let $\sigma\in W\cap L^{2}(\mathcal{H})$ 
and $\tau \in \mathcal L_c(e^0)$.
Choose $\psi\in\Lip_c(\R_+)$ with $\psi\tau=\tau$.
Then, by the first part of (2) and the choice of $\psi$,
\begin{equation}\label{wl2dm}
\begin{split}
  ( D_{\ext}\sigma, \tau )_{L^2(\mathcal H)}
  &= ( D_{\ext}(\psi\sigma), \tau )_{L^2(\mathcal H)}
  = ( D_{\max}(\psi\sigma), \tau )_{L^2(\mathcal H)} \\
  &= ( \psi\sigma, D\tau )_{L^2(\mathcal H)}
  = ( \sigma, D\tau )_{L^2(\mathcal H)} ,
\end{split}
\end{equation}
and hence $\sigma\in\mathcal D_{\max}$.
The converse inclusion is clear.
\end{proof}

\begin{lem}\label{wpclem}
Let $U$ be a bounded subset of $W$.
Then $U$ is precompact if and only if
$D_{\ext}(U)\subset L^2(\mathcal H)$
and $Q_{\ge}\mathcal R(U)\subset\check H$
are both precompact.
\end{lem}

\begin{proof}
If $U$ is precompact, then also its image
under the continuous maps $D_{\ext}$
and $Q_\ge\mathcal R$.

Vice versa, assume that
$D_{\ext}(U)\subset L^2(\mathcal H)$
and $Q_{\ge}\mathcal R(U)\subset\check H$
are both precompact.
By the definition of $W$, it suffices to show
that $\mathcal R(U)$ is precompact in $\check H$.

Let $D^0$ be the Dirac operator associated to $d^0$.
Let $\varphi,\psi\in\Lip_c(\R_+)$
such that $\varphi\psi=\psi$.
The operator $S_{D^0}\varphi$ is
the norm limit of the Hilbert-Schmidt operators
$S_{D^0}\varphi Q_{[-n,n]}$ on $L^2(\R_+,H_0)$,
hence $S_{D^0}\varphi$ is a compact operator.
On the other hand,
$\psi U\subset\mathcal D^0_{\max}$
and $D^0(\psi U)$ is bounded
in $L^2(\R_+,H_0)$,
see \tref{dcomp} and \lref{lem:west3gen}.2.
It follows that
$S_{D^0}\varphi(D^0(\psi U))$
is precompact in $L^2(\R_+,H_0)$.
By \cref{cc:rf},
\[
  \psi U \subset
  \psi(0) \mathcal E^0Q_>\mathcal R(U)
  + S_{D^0}\varphi(D^0(\psi U))
  + Q_0 (\psi U) ,
\]
hence $\psi U$ is precompact in $L^2(\R_+,H_0)$.

Now choose $\varphi,\psi$ as above with $\psi$ smooth
and equal to $1$ in a neighborhood of $0$.
We have
\[
  D_{\ext}(\psi U) \subset
  \gamma\psi' U + \psi D_{\ext}(U) .
\]
Since $\psi'$ is in $\Lip_c(\R_+)$
with $\varphi\psi'=\psi'$,
$\psi' U$  is precompact in $L^2(\R_+,H_0)$,
by the first part of the proof.
By assumption, $\psi D_{\ext}(U)$ is precompact
in $L^2(\R_+,H_0)$.
Hence $\psi U$ and $D_{\ext}(\psi U)$
are precompact in $L^2(\R_+,H_0)$,
hence $\psi U$ is precompact in $\mathcal D_{\max}$.
We conclude that $\mathcal R(U)=\mathcal R(\psi U)$
is precompact in $\check H$,
and hence that $U$ is precompact in $W$.
\end{proof}

For a boundary condition $B\subset\check H$, set
\begin{equation}\label{bcwb}
  W_B := \{\sigma\in W \,:\, \sigma(0) \in B \}
  \quad\text{and}\quad
  D_{B,\ext} := D_{\ext}|W_B .
\end{equation}

\begin{nameit}{Theorem and Definition}
\label{thm:wfredgen}
Assume that $(d,V)$ is non-parabolic
and that $B$ is regular.
Then $D_{B,\ext}:W_{B}\to L^{2}(\mathcal{H})$
is a left-Fredholm operator
with $(\im D_{B,\ext})^{\perp} = \ker D_{B^a,\max}$
and index
\begin{equation*}
  \ind D_{B,\ext}
  = \dim\ker D_{B,\ext} - \dim\ker D_{B^a,\max} ,
\end{equation*}
called the {\em extended index} of $D_B$,
also denoted $\ind_{\ext}D_B$.
\end{nameit}

\begin{proof}
Let $(\sigma_n)$ be a bounded sequence in $W_B$
such that the sequence $(D_{\ext}\sigma_n)$
converges in $L^2(\mathcal H)$.
By the continuity of $\mathcal R$,
the sequence $(\mathcal R\sigma_n(0))$
is bounded in $B\subset\check H$.
By the regularity of $B$,
the sequence $(Q_\ge\mathcal R\sigma_n(0))$
has a convergent subsequence in $H^{-1/2}$
and hence in $B$.
Therefore, $(\sigma_n)$ has a convergent
subsequence in $W$, by \lref{wpclem}.
Finally,
since $\mathcal D_{B,\max}$ is dense in $W_B$
and ($D_{B,\max})^*=D_{B^a,\max}$,
we also have 
$(\im D_{B,\ext})^{\perp}=\ker D_{B^a,\max}$.
\end{proof}

We note some important consequences
of \tref{thm:wfredgen}.

\begin{nameit}{Corollary and Definition}
\label{cor:l2ind}
If $(d,V)$ is non-para\-bolic and $B$ is elliptic,
then the kernels of $D_{B}$ and $D_{B^a}$
have finite dimension,
and we can define the {\em $L^{2}$-index} of $D_{B}$
to be the number
\begin{equation*}
    L^{2}\text{-}\ind D_{B}
    := \dim\ker D_{B} - \dim\ker D_{B^a} .
    \qed
\end{equation*}
\end{nameit}

Suppose that $(d,V)$ is non-parabolic.
For $\Lambda\in\R$,
let $D_{<\Lambda,\max}:=D_{B,\max}$
and $D_{<\Lambda,\ext}:=D_{B,\ext}$,
where $B=\check H_{<\Lambda}=H_{<\Lambda}^{1/2}$,
and similarly with $\le$ substituted for $<$.
The boundary conditions $B=\check H_{<\Lambda}$
and $B=\check H_{\le\Lambda}$ are elliptic
with $B^a=H_{\le-\Lambda}^{1/2}$ 
and $B^a=H_{<-\Lambda}^{1/2}$,
respectively.
Hence $D_{<\Lambda}=D_{<\Lambda,\max}$ and,
furthermore,
$D_{<\Lambda,\ext}$ and $D_{\le\Lambda,\ext}$
are Fredholm operators with
\begin{equation}\label{eqextperp}
  (\im D_{\le\Lambda,\ext})^\perp
  = \ker D_{<-\Lambda}
  \subset \ker D_{<-\Lambda,\ext} ,
\end{equation}
see \tref{thm:wfredgen}.

\begin{prop}\label{injesu}
If $(d,V)$ is non-parabolic, 
then there is $\Lambda_{0}\ge0$ such that
$D_{<-\Lambda,\ext}$ is injective
and $D_{\le\Lambda,\ext}$ is surjective
for all $\Lambda\ge\Lambda_{0}$.
\end{prop}

\begin{proof}
For any $\Lambda\in\R$, $D_{<\Lambda,\ext}$
is a Fredholm operator.
In particular, 
\[
  E := \mathcal R_{\ext}(\ker D_{<0,\ext})
  \subset H^{1/2}
\]
has finite dimension,
and hence all $H^{s}$-norms are equivalent on $E$ 
for $|s|\le 1/2$. 
Let $\Lambda\ge0$,
$\sigma\in\ker D_{<-\Lambda,\ext}\subset\ker D_{<0,\ext}$,
and suppose that $\sigma(0)\ne0$.
Since $\sigma(0)\in E\cap\check H_{<-\Lambda}$, we can estimate
\begin{align*}
  0 &\ne ||\sigma(0)||_{1/2}^{2} 
  \le C_{E}^{2}||\sigma(0)||_{-1/2}^{2} \\
  &= C^{2}_{E}\la(I +A_{0}^{2})^{-1/2}\sigma(0),\sigma(0)\ra
  < C^{2}_{E}(1+\Lambda^2)^{-1}||\sigma(0)||_{1/2}^{2} ,
\end{align*}
a contradiction if
\[
  \Lambda \ge \Lambda_0
  := (C_E^2 - 1)^{1/2} .
\]
Therefore $\sigma(0)=0$ if $\Lambda \ge \Lambda_0$,
and then $\sigma=0$, by the non-parabolicity of $(d,V)$.
Hence $D_{<-\Lambda,\ext}$ is injective
for $\Lambda \ge \Lambda_0$.
\end{proof}

Next we would like to write the index formula
in \tref{thm:wfredgen}
in a way analogous to \cref{bcbsr}.
For this, we need the {\em Calder\'on spaces}
\begin{equation}\label{caldercheck}
  \check{\mathcal C}_{\max}
  := \mathcal{R}(\ker D_{\max})
  \quad\text{and}\quad
  \check{\mathcal C}_{\ext}
  := \mathcal{R}(\ker D_{\ext}) .
\end{equation}
Since $\mathcal{R}:\ker D_{\ext}\to\check H$
is isometric, $\check{\mathcal C}_{\ext}$
is a closed subspace of $\check H$.
For $|s|\le1/2$, we let
\begin{equation}\label{eq:sols}
  \mathcal{C}_{\max}^s
  := \check{\mathcal{C}}_{\max} \cap H^{s}
  \quad\text{and}\quad
  \mathcal{C}_{\ext}^s
  := \check{\mathcal{C}}_{\ext} \cap H^{s} .
\end{equation}
If $B$ is a regular boundary condition,
then $\mathcal R$ induces isomorphisms
\begin{equation}
\begin{split}
  \ker D_{B,\max}
  &\cong B\cap\check{\mathcal C}_{\max}
  = B\cap\mathcal C_{\max}^{1/2} , \\
  \ker D_{B,\ext}
  &\cong B\cap\check{\mathcal C}_{\ext}
  = B\cap\mathcal C_{\ext}^{1/2} .
\end{split}
\end{equation}
We will write $\mathcal{C}_{\max}$ and $\mathcal{C}_{\ext}$
instead of $\mathcal{C}^{0}_{\max}$ and
$\mathcal{C}^{0}_{\ext}$, respectively.

\begin{cor}\label{cor:wfredgen}
If $(d,V)$ is non-parabolic and $B$ is elliptic,
then $D_{B,\ext}$ is a Fredholm operator with 
$(\im D_{B,\ext})^{\perp} = \ker D_{B^a}$
and index
\begin{align*}
  \ind D_{B,\ext}
  &= \dim B\cap\mathcal C_{\ext}^{1/2}
  - \dim B^\perp\cap\gamma\mathcal C_{\max}^{1/2}. \\
  &= \dim B\cap\mathcal C_{\ext}
  - \dim B^\perp\cap\gamma\mathcal C_{\max} .
\end{align*}
\end{cor}

\begin{proof}
The assertions follow easily from \tref{thm:wfredgen}
and \lref{bcbann},
except for the last identity.
Since $B$ is elliptic, 
we have $B\subset H^{1/2}$ and
\[
  B^a = (\gamma B^\perp) \cap \check H
  \subset H^{1/2} \subset H .
\]
Therefore
\[
  B^a \cap \mathcal C_{\max}^{1/2}
  = B^a \cap \check{\mathcal C}_{\max}
  = (\gamma B^\perp) \cap \check{\mathcal C}_{\max}
  = (\gamma B^\perp) \cap \mathcal C_{\max} .
  \qedhere
\]
\end{proof}

\begin{cor}\label{cor:wfredgenp}
Assume that $(d,V)$ is non-parabolic and 
that $P$ is an orthogonal elliptic projection in $H$.
Then $D_{P,\ext}$ is a Fredholm operator
with $(\im D_{P,\ext})^{\perp} = \ker D_{P_\gamma}$
and index
\begin{align*}
  \ind D_{P,\ext}
  &= \dim\ker P\cap\mathcal C_{\ext}^{1/2}
    - \dim\im P\cap\mathcal \gamma\mathcal C_{\max}^{1/2} \\
  &= \dim\ker P\cap\mathcal C_{\ext}
    - \dim\im P\cap\mathcal \gamma\mathcal C_{\max} .
\end{align*}
\end{cor}

\begin{proof}
The boundary condition associated to $P$
is $B_P=\ker\tilde P\cap\check H$, see \eqref{rpbp}.
Since $B_P$ is regular, 
$B_P=\ker P\cap\check H=\ker\hat P$
and therefore
\[
  B_P\cap\check{\mathcal C}_{\ext}
  = \ker P\cap\check{\mathcal C}_{\ext}
  = \ker P\cap\mathcal C_{\ext}
  = \ker P\cap\mathcal C_{\ext}^{1/2} .
\]
The remaining identities follow from \cref{cor:wfredgen}
since $\im P$ is the orthogonal complement of $\ker\hat P$
in $H$.
\end{proof}

\subsection{Some examples}\label{somexa}
The first two examples are Dirac systems 
on $\R_+$ which are not non-parabolic.
In the first example,
$\ker D_{P_{APS},\max}$ is infinite-dimensional
so that $D_{P_{APS},\ext}$ cannot be a Fredholm operator.
In the second example, 
the assumption of non-paraboli\-city would lead to the
contradiction that $\ker D_{P_{APS},\ext}$
has infinite dimension.
These examples are modelled on the Gauss-Bonnet operators
of real hyperbolic spaces of even and odd dimension.

\begin{exa}\label{nopaexa1}
For $t\in\R_+$ and $k\in\Z$, let
\[
  B_{t}(k) =
  \begin{pmatrix}
  1 & ike^{-t} \\ -ike^{-t} & 1 \\
  \end{pmatrix} ,
\]
and consider the evolution equation
\[
  \sigma' + B_t(k) \sigma = 0 .
\]
Solutions $\sigma$ of this equation satisfy
$(||\sigma||^2)'\le-2(1-|k|e^{-t})||\sigma||^2$,
hence belong to $L^2(\R_+,\C^2)$.
Eigenvalues and eigenvectors of $B_0(k)$
are given by
\[
  B_0(k)\begin{pmatrix} 1 \\ i \end{pmatrix}
  = (1-k) \begin{pmatrix} 1 \\ i \end{pmatrix}
  \quad\text{and}\quad
  B_0(k)\begin{pmatrix} 1 \\ -i \end{pmatrix}
  = (1+k) \begin{pmatrix} 1 \\ -i \end{pmatrix} .
\]
On $L^2(\R_+,\C^2\oplus \C^2)$,
consider the Dirac system
\begin{equation*}
\begin{split}
  D_k \sigma
  &= \begin{pmatrix}
  -\sigma_2' + B_t(k)\sigma_2 \\
  \sigma_1'+B_t(k)\sigma_1
  \end{pmatrix} \\
  &= \begin{pmatrix} 0 & -I \\ I & 0 \end{pmatrix}
  \left( \partial_t  + \begin{pmatrix}
  B_t(k) & 0 \\ 0 & -B_t(k) \end{pmatrix} \right)
  \begin{pmatrix} \sigma_1 \\ \sigma_2 \end{pmatrix} \\
  &=: \gamma ( \partial_t  + A_t(k) ) \sigma
\end{split}
\end{equation*}
For any $k\in\Z$, let
\[
  \sigma_k
  :=  \begin{pmatrix} \tau_k \\ 0 \end{pmatrix}
  \quad\text{with}\quad
  \tau_k' + B_t(k)\tau_k = 0
   \quad\text{and}\quad
  \tau_k(0) = \begin{pmatrix} 1 \\ i \end{pmatrix} .
\]
Then $\sigma_k\in L^2(\R_+,\C^2\oplus \C^2)$, $D_k\sigma_k=0$,
and $A_0(k)\sigma_k(0)=(1-k)\sigma_k$.
Hence $\sigma_k$ belongs to the negative eigenspace of $A_0(k)$
for $k\ge 2$.

We can now sum these Dirac systems to obtain a Dirac system
\[
  d=\left(\mathcal H, \partial_t,
  A_t = \oplus A_t(k), \gamma \right)
  \quad\text{on}\quad
  \mathcal H
  = \R_+\times l^2(\Z,\C^2\oplus \C^2)
\]
with associated Dirac operator $D=\oplus D_k$.
For this Dirac system, there is a family $(\sigma_k)$
of orthogonal non-zero $L^2$-sections of $\mathcal H$
with $D\sigma_k=0$ and $A_0\sigma_k=(1-k)\sigma_k$.
Hence, with $Q_{\ge\Lambda}$
the corresponding spectral projection of $A_0$,
the $L^2$-kernel of $D_{Q_{\ge\Lambda}}$
has infinite dimension, for any $\Lambda\in\R$.
In particular, $d$ is not non-parabolic.
\end{exa}

\begin{exa}\label{nopaexa2}
For $k\in \Z$, consider the Dirac system on $\R_+\times\C^2$
with Dirac operator
\begin{equation*}
\begin{split}
  D_k\sigma
  &= \begin{pmatrix} 0 & -1 \\ 1 & 0 \end{pmatrix}
  \left(\partial_t + \begin{pmatrix}
  ke^{-t} & 0 \\ 0 & -ke^{-t} \end{pmatrix} \right)
  \begin{pmatrix} \sigma_1 \\ \sigma_2 \end{pmatrix} \\
  &=: \gamma (\partial_t + A_t(k))\sigma .
\end{split}
\end{equation*}
Solutions of the equation $D_k\sigma=0$ are obviously 
uniformly bounded and, therefore, admit an upper bound
\[
  \int_0^T||\sigma(t)||^2 dt \le C_k T ||\sigma(0)||^2 .
\]
Moreover, for $k \ge 1$,
\[
  \sigma_k(t)
  = \begin{pmatrix} 0 \\ e^{-ke^{-t}} \\ \end{pmatrix}
\]
satisfies $D_k\sigma_k=0$
and $A_0(k)\sigma_k(0)=-k\sigma_k(0)$.
Again, we sum all these Dirac systems to get a
Dirac system on $L^2(\R_+, l^2(\Z,\C^2))$ given by
$\partial_t$, $A_t=\oplus A_t(k)$ and $\oplus\gamma$.

Let $Q_{\ge 0}$ be the spectral projection of $A_0$
onto the non-negative eigenspaces of $A_0$.
We obtain that the space
of $\sigma\in L^2(\R_+, l^2(\Z,\C^2))$ with
\begin{equation*}
  D\sigma=0 , \quad  Q_{\ge 0}\sigma(0)=0
  \quad\text{and}\quad
  \int_0^T ||\sigma(t)||^2 dt =O(T)
\end{equation*}
has infinite dimension.
The following lemma implies that this Dirac system
is not non-parabolic.

\begin{lem}\label{nopaexa2lem}
Let $d$ be a non-parabolic Dirac system.
If $\sigma\in H^1_{\loc}(e)$ satisfies $D\sigma=0$ and
\[
  \lim_{T\to\infty}
  \frac{ \int_0^T ||\sigma(t)||^2dt}{T^2} = 0 ,
\]
then $\sigma\in W$.
\end{lem}

\begin{proof}
It suffices to find a sequence $(\sigma_n)$ in $H^1_{c}(e^0)$
such that
\[
  \lim_{n\to\infty}
  || D(\sigma-\sigma_n) ||_{L^2(\mathcal H)}
  + || \sigma(0)-\sigma_n(0) ||_{\check H} = 0 .
\]
Let $\psi$ be a Lipschitz function on $\R_+$ with compact
support such that $\psi=1$ in a neighborhood of $0$,
and set $\psi_n(t):=\psi(t/n)$ and $\sigma_n:=\psi_n\sigma$.
Since $\sigma(0)=\sigma_n(0)$ and
\[
  D(\sigma-\sigma_n)(t)
  = -\frac{1}{n}\gamma\psi'(t/n)\sigma,
\]
we obtain that $\|D(\sigma-\sigma_n)\|^2_{L^2(\mathcal H)}=o(1)$
as $n$ tends to infinity.
\end{proof}
\end{exa}

\begin{exa}\label{nopaexa3}
For $\mu\in\R$, let $d_\mu$ be the Dirac system
on $L^2([1,\infty),\C^2)$ with Dirac operator
\[
  D(\sigma_+,\sigma_-) 
  := (-\sigma_-'+\frac{\mu}{t}\sigma_-,
      \sigma_+'+\frac{\mu}{t}\sigma_+) .
\]
Clearly, $d_\mu$ is not of Fredholm type.
On the other hand, 
since the equation $D\sigma=\tau$ corresponds to a linear
ODE in the finite dimensional space $H=\C^2$,
$d_\mu$ is non-parabolic.
We have
\[
  |-\sigma_-' + \frac{\mu}{t}\sigma_- |^2
  = |\sigma_-'|^2 + \frac{\mu(\mu-1)}{t^2}|\sigma_-|^2
    - ( \frac{\mu}{t} |\sigma_-|^2 )' ,
\]
and similarly for $\sigma_+$, 
where all the minus signs turn into plus signs.
Now $W$ is the closure of the
space of Lipschitz sections with compact support
with respect to the $W$-norm.
Hence, if $\mu>1$ and $\sigma=(\sigma_+,\sigma_-)$ is in $W$,
then $|\sigma/t|^2$ is integrable with integral
uniformly bounded by the $W$-norm of $\sigma$.
(This also shows non-parabolicity in the case $\mu>1$.)

The space of solutions of the equation
$D\sigma=0$ is given by the space of sections 
$(at^{-\mu},bt^\mu)$ with $a,b\in\C$.
For $\mu>1$ and $b\ne0$, 
$(at^{-\mu},bt^\mu)$ does not belong to $W$
since $(at^{-\mu},bt^\mu)/t$ is not square integrable.
It follows that, for $\mu>1$, 
$W$-solutions of the equation $D\sigma=0$ 
are square integrable,
hence that $\mathcal C_{\max}=\mathcal C_{\ext}$,
although $d_\mu$ is not of Fredholm type.

The above analysis can be refined.
By (5.3) in \cite{Ca2} and by what is said in
the two lines above it,
\[
  \int_1^\infty \big(
  |\tau'|^2 - \frac{1}{4t^2}|\tau|^2 \big)
  \ge \int_1^\infty \frac1{4t^2(\ln t)^2}|\tau|^2 ,
\]
for all $\tau\in\Lip_c([1,\infty))$ with $\tau(1)=0$.
Since
\[
  |\tau' - \frac{\mu}{t}\tau |^2
  = |\tau'|^2 - \frac{1}{4t^2}|\tau|^2
    + \frac{(\mu-1/2)^2}{t^2}|\tau|^2
    - ( \frac{\mu}{t} |\tau|^2 )' ,
\]
we get the following inequality 
\[
  \int_1^\infty |\tau' - \frac{\mu}{t}\tau |^2
  \ge \int_1^\infty \frac{(\mu-1/2)^2}{t^2}|\tau|^2
   + \int_1^\infty \frac1{4t^2(\ln t)^2}|\tau|^2 ,
\]
for all $\tau\in\Lip_c([1,\infty))$ with $\tau(1)=0$.
It follows that $\mathcal C_{\max}=\mathcal C_{\ext}$ 
if $|\mu|>1/2$
(and, again, that $d_\mu$ is non-parabolic for all $\mu\in\R$).
\end{exa}

\newpage
\section{Calder\'on projections and index formulas}\label{calderon}

\subsection{The Calder\'on projections}
Recall the definition of the Calder\'on spaces
in \eqref{caldercheck} and \eqref{eq:sols}.

\begin{thm}\label{cald1}
Let $(d,V)$ be non-parabolic.
If $\Lambda_{0}\ge0$ is the constant from \pref{injesu}
and $\Lambda\ge\Lambda_{0}$,
then we have a direct sum decomposition
\begin{equation*}
  \check{\mathcal C}_{\ext} =
  K_{\Lambda} \oplus \check G_{\Lambda} ,
\end{equation*}
where
$K_{\Lambda}=\{x\in\check{\mathcal C}_{\ext} \,:\,
 Q_{>\Lambda}x=0\}\subset H^{1/2}$ is of finite dimension
and $\check G_{\Lambda}$ is the graph
of a continuous linear map
 \[
   T_{\Lambda}: H_{>\Lambda}^{-1/2}
            \to H_{\le\Lambda}^{1/2} ,
 \]
where
$T_{\Lambda}=T_{\Lambda_0}|H_{>\Lambda}^{-1/2}$.
The finite rank and remainder parts
$Q_{[-\Lambda,\Lambda]}T_{\Lambda}$
and $Q_{<-\Lambda}T_{\Lambda}$, respectively,
satisfy
\begin{equation*}
  ||Q_{[-\Lambda,\Lambda]}T_{\Lambda}||_{s}
  \le C \Lambda^{-1/2-s}
  \quad\text{and}\quad
  ||Q_{<-\Lambda}T_{\Lambda}||_{s}
  \le C \Lambda^{-1} ,
\end{equation*}
where $C$ is a constant independent
of $\Lambda\ge\Lambda_{0}$ and $s\in[-1/2,1/2]$.
In particular,
\[
  \mathcal C_{\ext}^s =
  K_{\Lambda} \oplus G_{\Lambda}^s ,
\]
where $G_{\Lambda}^s=\check G_{\Lambda}\cap H^s$ is
the graph of $T_{\Lambda}|H_{>\Lambda}^{s}$,
and hence $\mathcal C_{\ext}^s$
is a closed subspace of $H^s$,
for all $\Lambda\ge\Lambda_{0}$ and $s\in[-1/2,1/2]$.
\end{thm}

\begin{proof}
Throughout the proof, we assume $\Lambda\ge\Lambda_0$,
where $\Lambda_0$ is the constant from \pref{injesu}.

Let $x\in\check H_{>\Lambda}=H^{-1/2}_{>\Lambda}$.
Choose a function $\psi\in\Lip_{c}(\R_+)$
which is equal to $1$ in a neighborhood of $0$
and set $\sigma:=\psi\mathcal E^0x$.
Then $\sigma\in\mathcal D_{\max}\subset W$,
by \tref{dcomp}.
Since $D_{\le\Lambda,\ext}$ is surjective,
there is $\tau\in\mathcal D_{\le\Lambda,\ext}$
with $D_{\ext}\tau=D_{\ext}\sigma$.
Hence $\sigma-\tau\in\ker D_{\ext}$ and
\[
  x = Q_{>\Lambda}((\sigma-\tau)(0))
  \in Q_{>\Lambda}(\check{\mathcal C}_{\ext}) .
\]
Therefore $Q_{>\Lambda}:\check{\mathcal C}_{\ext}
\to H^{-1/2}_{>\Lambda}$ is surjective.
We have
\[
  K_\Lambda =
  \mathcal R_{\ext}(\ker D_{\ext}) \cap H^{1/2}_{\le\Lambda}
  =   \mathcal R_{\ext}(\ker D_{\le\Lambda,\ext}) ,
\]
hence $K_{\Lambda}$ is of finite dimension,
by \tref{thm:wfredgen}.
Let $\check G_{\Lambda_0}$ be a complement
of $K_{\Lambda_0}$ in $\check{\mathcal C}_{\ext}$.
Then $Q_{>\Lambda_0}:\check G_{\Lambda_0}
   \to H^{-1/2}_{>\Lambda_0}$ is an isomorphism,
hence $\check G_{\Lambda_0}$ is the graph
of a continuous linear map
\[
  T_{\Lambda_0}: H^{-1/2}_{>\Lambda_0}
  \to \check H_{\le\Lambda_0}
  = H^{1/2}_{\le\Lambda_0} .
\]
This is the place where we gain regularity:
By the very structure of $\check H$,
$T_{\Lambda_0}$ extends naturally to a smoothing operator.

Let $|s|\le1/2$.
Since the image of $T_{\Lambda_0}$
is contained in $H^{1/2}$,
$x=y+T_{\Lambda_0}y\in\check G_{\Lambda_0}$ is in $H^s$
if and only if $y=Q_{>\Lambda_0}x$ is in $H^s$, i.e.,
\[
  G_{\Lambda_0}^s =
  \check G_{\Lambda_0} \cap H^s
  = \{ y+T_{\Lambda_0}y \,:\, y \in H^{s}_{>\Lambda_0} \} .
\]
For $\Lambda\ge\Lambda_0$, we define
\begin{align*}
  \check G_{\Lambda} :&=
  \{x\in \check G_{\Lambda_0}
    \,:\, Q_{(\Lambda_0,\Lambda]}x=0\}
  \\ &= \{y + T_{\Lambda_0}y
        \,:\, y \in H^{-1/2}_{>\Lambda} \} .
\end{align*}
Let $|s|\le1/2$.
Then $G_{\Lambda}^s=\check G_{\Lambda}\cap H^s$
is the graph of
$T_{\Lambda}:=T_{\Lambda_0}|H_{>\Lambda}^s$.
Hence $G_{\Lambda}^s$ is
a closed subspace of $H^s$.

We show next that $G_{\Lambda}^{s}$ is a complement
of $K_{\Lambda}$ in 
$\mathcal C_{\ext}^s=\check{\mathcal C}_{\ext}\cap H^s$.
Since $K_{\Lambda_0}\subset K_{\Lambda}$ and,
clearly, $K_{\Lambda}\cap G_{\Lambda}^s=0$,
it is enough to show 
that $G_{\Lambda_0}^s\subset K_{\Lambda}+G_{\Lambda}^s$.
Now for $y\in G_{\Lambda_0}^s$ 
there is $z\in\check G_{\Lambda_0}$
with $Q_{>\Lambda_0}z=Q_{(\Lambda_0,\Lambda]}y$,
by the surjectivity
of $Q_{>\Lambda_0}|\check G_{\Lambda_0}$.
It follows that $z\in K_{\Lambda}$
and $y-z\in G_{\Lambda}^{s}$.

For $x\in H^{-1/2}_{>\Lambda}$,
\[
  ||Q_{[-\Lambda,\Lambda]}T_{\Lambda}x||_{1/2}
  \le ||T_{\Lambda}x||_{1/2} \le C ||x||_{-1/2} ,
\]
where $C=||T_{\Lambda_0}||_{\check H}$,
and similarly for $Q_{<-\Lambda}T_{\Lambda}x$.
For $r<t$ and $y\in H^t$ with $Q_{(-\Lambda,\Lambda)}y=0$,
we have 
\begin{equation}\label{compnorm}
  ||y||_{r} \le \Lambda^{r-t} ||y||_{t} .
\end{equation}
Hence
\begin{align}
  ||Q_{[-\Lambda,\Lambda]}T_{\Lambda}x||_{s}
  &\le ||Q_{[-\Lambda,\Lambda]}T_{\Lambda}x||_{1/2}
  \notag \\
  &\le C ||x||_{-1/2}
  \le C \Lambda^{-1/2-s} ||x||_{s} ,
  \label{compnorm2} \\
  ||Q_{<-\Lambda}T_{\Lambda}x||_{s}
  &\le \Lambda^{s-1/2} ||Q_{<-\Lambda}T_{\Lambda}x||_{1/2}
  \notag \\
  &\le C \Lambda^{s-1/2} ||x||_{-1/2}
  \le C \Lambda^{-1} ||x||_{s} ,
  \label{compnorm4}
\end{align}
for all $|s|\le1/2$ and $x\in H^{s}_{>\Lambda}$.
\end{proof}

\begin{dfn}\label{caldprodef}
The orthogonal projections in $H$ onto (the closure of)
$\mathcal{C}_{\max}=\mathcal{C}_{\max}^0$
and onto $\mathcal{C}_{\ext}=\mathcal{C}_{\ext}^0$
will be called the {\em Calder\'on projection}
and the {\em extended Calder\'on projection}
associated to the Dirac-Schr\"odinger system $(d,V)$
and will be denoted by $C_{\max}$ and $C_{\ext}$,
respectively.
\end{dfn}

\begin{thm}\label{cald2}
Let $(d,V)$ be non-parabolic.
Then there are constants $\Lambda_{0},C\ge0$ such that,
for $\Lambda\ge\Lambda_{0}$,
\begin{equation*}
  C_{\ext} =
  Q_{>} + R_{\Lambda} + S_{\Lambda},
\end{equation*}
where $R_{\Lambda}$ and $S_{\Lambda}$ are smoothing,
$R_{\Lambda}$ has finite rank, and
\begin{equation*}
  ||S_{\Lambda}||_{s}
  + ||S_{\Lambda}^{*}||_{s}
  \le C\Lambda^{-1}
\end{equation*}
for all $|s|\le1/2$.
In particular, $C_{\ext}$ is $1/2$-smooth,
$C_{\ext} - Q_{>}$ is compact in $H^{s}$
for all $|s|\le1/2$,
and $C_{\ext}$ is elliptic.
Furthermore, 
$D_{C_{\ext},\ext}:W_{C_{\ext}}\to L^2(\mathcal H)$
is an isomorphism.
\end{thm}

\begin{proof}
We use notation and results from \tref{cald1}.
Since $T_{\Lambda}$ maps $H_{>\Lambda}^{-1/2}$
to $H_{\le\Lambda}^{1/2}$,
the dual operator of $T_{\Lambda}$
maps $H_{\le\Lambda}^{-1/2}$ to $H_{>\Lambda}^{1/2}$.
We recall that the dual operator of $T_{\Lambda}$
is the extension of the adjoint $T_{\Lambda}^*$
of $T_{\Lambda}|H_{>\Lambda}$.
In particular, $T_{\Lambda}^*$ is smoothing as well
and, considered as a linear map
from $H_{\le\Lambda}^{-1/2}$ to $H_{>\Lambda}^{1/2}$,
it satisfies
\[
  ||T_{\Lambda}^*||
  =||T_{\Lambda}||
  \le||T_{\Lambda_0}|| =C .
\]
Arguing as in \eqref{compnorm4}, we obtain
that $T_{\Lambda}^*T_{\Lambda}:
 H_{>\Lambda}^{-1/2}\to H_{>\Lambda}^{1/2}$ satisfies
\begin{align*}
  ||T_{\Lambda}^*T_{\Lambda}x||_s
  &\le \Lambda^{s-1/2}
      ||T_{\Lambda}^*T_{\Lambda}x||_{1/2} \\
  &\le C \Lambda^{s-1/2} ||x||_{-1/2}
  \le C \Lambda^{-1} ||x||_{s} ,
\end{align*}
for all $|s|\le1/2$ and $x\in H^{s}_{>\Lambda}$.
Hence $||T_{\Lambda}^*T_{\Lambda}||_s\le C \Lambda^{-1}$
for all $|s|\le1/2$.
In particular,
if $I$ denotes the identity of $H_{>\Lambda}$,
then $I+T_{\Lambda}^*T_{\Lambda}$
is invertible with $1/2$-smooth inverse
as soon as $\Lambda>C$, 
and for $\Lambda\ge2C$ we find
\[
  ||(I+T_{\Lambda}^*T_{\Lambda})^{-1}||_s \le 2 .
\]
Clearly,
\[
  (I+T_{\Lambda}^*T_{\Lambda})^{-1}
  = I - T_{\Lambda}^*T_{\Lambda}
    (I+T_{\Lambda}^*T_{\Lambda})^{-1}
  =: I + T_{\Lambda}^x ,
\]
where $T_{\Lambda}^x$ is smoothing
with $||T_{\Lambda}^x||_s\le 2C \Lambda^{-1}$
and the superscript $x$ means that this object
will not survive the end of the proof.

In accordance with with our convention $H=H^0$,
we let $G_{\Lambda}=G_{\Lambda}^0$.
Then $G_{\Lambda}$ is the graph of the restriction
of $T_{\Lambda}$ to $H_{>\Lambda}$,
for short also denoted by $T_{\Lambda}$.
We recall that
\[
  G_{\Lambda}^\perp
  = \{ (-T_{\Lambda}^*y,y) \,:\, y \in H_{\le\Lambda} \} .
\]
Hence the orthogonal projection $P_{\Lambda}$
onto $G_{\Lambda}$ in $H$ is given by
\[
  P_{\Lambda} =
  \begin{pmatrix}
  (I+T_{\Lambda}^*T_{\Lambda})^{-1} &
  (I+T_{\Lambda}^*T_{\Lambda})^{-1}T_{\Lambda}^* \\
  T_{\Lambda}(I+T_{\Lambda}^*T_{\Lambda})^{-1} &
  T_{\Lambda}(I+T_{\Lambda}^*T_{\Lambda})^{-1}T_{\Lambda}^*
  \end{pmatrix} ,
\]
where the operator matrix arises from the decomposition
$H_{>\Lambda}\oplus H_{\le\Lambda}$ of $H$
and $I$ denotes the identity of $H_{>\Lambda}$
as above.
We now get a representation
\[
  P_{\Lambda}
  = Q_{>\Lambda} + R_{\Lambda}^x + S_{\Lambda}
\]
analogous to the asserted representation for $C_{\ext}$,
where
\begin{align*}
  R_{\Lambda}^x &=
  \begin{pmatrix}
  0 &
  (I+T_{\Lambda}^x)(Q_{[-\Lambda,\Lambda]}T_{\Lambda})^* \\
  Q_{[-\Lambda,\Lambda]}T_{\Lambda}
  (I+T_{\Lambda}^*T_{\Lambda})^{-1} &
  Q_{[-\Lambda,\Lambda]}T_{\Lambda}
  (I+T_{\Lambda}^*T_{\Lambda})^{-1}T_{\Lambda}^*
  \end{pmatrix} , \\
  S_{\Lambda} &=
  \begin{pmatrix}
  T_{\Lambda}^x &
  (I+T_{\Lambda}^x)(Q_{<-\Lambda}T_{\Lambda})^* \\
  Q_{<-\Lambda}T_{\Lambda}
  (I+T_{\Lambda}^*T_{\Lambda})^{-1} &
  Q_{<-\Lambda}T_{\Lambda}
  (I+T_{\Lambda}^*T_{\Lambda})^{-1}T_{\Lambda}^*
  \end{pmatrix} .
\end{align*}
Obviously, $R_{\Lambda}^x$ and $S_{\Lambda}$
are smoothing, $R_{\Lambda}^x$ has finite rank,
and the operator norms of $S_{\Lambda}$ satisfy
the desired inequalities.

The orthogonal complement of $G_{\Lambda}$
in $\mathcal C_{\ext}$ 
is $(I-P_{\Lambda})(K_{\Lambda})\subset H^{1/2}$
so that $C_{\ext}-P_{\Lambda}$ is smoothing
of finite rank.
This implies the asserted formula for $C_{\ext}$
with $R_{\Lambda}=
 R_{\Lambda}^x+C_{\ext}-P_{\Lambda}-Q_{(0,\Lambda]}$.

By \pref{injesu}, $D_{\ext}:W\to L^2(\mathcal H)$
is surjective.
By definition, the kernel of $D_{C_{\ext},\ext}$
is trivial.
The theorem follows.
\end{proof}

\begin{thm}\label{cald3}
Assume that $(d,V)$ is non-parabolic.
Then
\[
  C_{\max} =C_{\ext,\gamma}
  = \gamma^*(I-C_{\ext})\gamma .
\]
In particular, $C_{\max}$ is elliptic,
$C_{\max} - Q_{>}$ is compact in $H^{s}$
for all $|s|\le1/2$,
and $\ind D_{C_{\max},\ext}=\rk(C_{\ext}-C_{\max})$.
\end{thm}

\begin{proof}
Let $x\in\check{\mathcal C}_{\max}$
and $y\in\check{\mathcal C}_{\ext}$.
Choose $\sigma\in\ker D_{\max}$ with $\sigma(0)=x$
and $\tau\in\ker D_{\ext}$ with $\tau(0)=y$.
Let $(\tau_n)$ be a sequence in $\mathcal L_c(e^0)$
which converges to $\tau$ in $W$.
Then $D\tau_n\to 0=D_{\ext}\tau$ in $L^2(\mathcal H)$
and $\tau_n(0)\to\tau(0)$ in $\check H$.
By \tref{thm:diprops}.5,
\[
  \omega(x,y)
  \leftarrow \omega(\sigma(0),\tau_n(0))
  = ( D_{\max}\sigma,\tau_n )_{L^2(\mathcal H)}
  - ( \sigma,D\tau_n )_{L^2(\mathcal H)}
  \to 0 .
\]
We conclude that
$\check{\mathcal C}_{\ext}
 \subset(\check{\mathcal C}_{\max})^a$
and hence that $\mathcal C_{\ext}
    \subset\gamma(\mathcal C_{\max})^\perp$.

Suppose now that $\mathcal C_{\ext}$ is not equal
to $\gamma(\mathcal C_{\max})^\perp$.
Then there is a vector $z$ of norm $1$
in $\gamma(\mathcal C_{\max})^\perp$
which is perpendicular to $\mathcal C_{\ext}$.
Choose $y\in H^{1/2}$
with $||y-z||_{H}\le 1/2$
and set $x := (I-C_{\ext})y$.
Then $x$ is non-zero,
$x\notin\gamma\mathcal C_{\max}^{1/2}$,
and is perpendicular to $\mathcal C_{\ext}$.
Furthermore, $x\in H^{1/2}$
since $C_{\ext}$ is $1/2$-smooth.
Let $P := C_{\ext} + R$, where $R$ is the orthogonal
projection onto $\C x$ in $H$.
Then $P$ is an elliptic orthogonal projection,
by \lref{cc:stabsym},
since $C_{\ext}$ is elliptic
and $x$ is in $H^{1/2}$.
By \cref{cor:wfredgen},
\begin{align*}
  \ind D_{P,\ext}
  &= \dim (\ker P\cap\mathcal C_{\ext}^{1/2})
  - \dim (\im P\cap\gamma\mathcal C_{\max}^{1/2}) \\
  &= - \dim (\im P\cap\gamma\mathcal C_{\max}^{1/2}) .
\end{align*}
Let $y\in\mathcal C_{\ext}$ and $\alpha\in\C$,
and suppose that
$y+\alpha x\in\gamma\mathcal{C}_{\max}^{1/2}$.
Since $\mathcal C_{\ext}$ is perpendicular
to $\gamma\mathcal{C}_{\max}^{1/2}$ and $x$,
we get $y=0$.
This implies that
$\alpha x\in\gamma\mathcal{C}_{\max}^{1/2}$
and hence that $\alpha=0$, by the choice of $x$.
Hence $\im P\cap\gamma\mathcal C_{\max}^{1/2}=0$
and, therefore, $\ind D_{P,\ext} = 0$.

On the other hand, the inclusion
$i_{P}:\mathcal D_{P,\ext}
  \to\mathcal{D}_{C_{\ext},\ext}$
is a Fredholm operator of index -1.
Since $D_{P,\ext}=D_{C_{\ext},\ext}\circ i_P$
and $D_{C_{\ext},\ext}$ is an isomorphism,
we get $\ind D_{P,\ext} = -1$, a contradiction.
We conclude that
$\mathcal C_{\ext}=\gamma(\mathcal C_{\max})^\perp$
and hence that $C_{\max}=C_{\ext,\gamma}$.

Since $C_{\ext}$ is elliptic
$C_{\max}=C_{\ext,\gamma}$ is elliptic as well.
Moreover,
\[
  C_{\max} - Q_>
  = \gamma^*(I-C_{\ext} - Q_<)\gamma
  = \gamma^*(Q_\ge - C_{\ext})\gamma ,
\]
hence $C_{\max} - Q_>$ is compact, 
by \tref{cald2}.

Finally, since $C_{\max}$ is elliptic,
$D_{C_{\max},\ext}$ is a Fredholm operator.
Now $\im C_{\max}\subset\mathcal C_{\ext}$, hence
\begin{align*}
  \ind D_{C_{\max},\ext}
  &= \dim (\ker C_{\max}\cap\mathcal C_{\ext})
  - \dim (\im C_{\max}\cap
    \gamma\mathcal C_{\max}^{1/2}) \\
  &= \dim (\ker C_{\max}\cap\mathcal C_{\ext})
  = \rk(C_{\ext}-C_{\max}) .
  \qedhere
\end{align*}
\end{proof}

\begin{cor}\label{fredtosa}
If $(d,V)$ is of Fredholm type,
then $C_{\ext}=C_{\ext,\gamma}$, that is,
$\ker\hat C_{\ext}$ is an elliptic
self-adjoint boundary condition.
\end{cor}

\begin{proof}
Since $(d,V)$ is of Fredholm type, 
we have $W=\mathcal D_{\max}$
and hence $C_{\ext}=C_{\max}$.
\end{proof}

\begin{thm}\label{cmaxthm}
Assume that $(d,V)$ is non-parabolic.
Then
\[
  \mathcal C_{\max}^{1/2}
  = \im\hat C_{\max}
  = \im C_{\max}\cap H^{1/2} . \tag{1}
\]
If $B$ is an elliptic boundary condition
and $B^{(s)}$ denotes the closure of $B$ in $H^s$,
where $|s|\le1/2$,
then $(B^{(s)},\mathcal C_{\ext}^s)$ 
is a Fredholm pair in $H^s$ 
with nullity and deficiency independent of $s$.
More precisely, we have
\begin{align*}
  \nuli(B^{(s)},\mathcal C_{\ext}^s)
  &= \dim (B \cap \mathcal C_{\ext}^{1/2}) , \tag{2} \\
  \defi (B^{(s)},\mathcal C_{\ext}^s)
  &= \dim (B^a \cap \mathcal C_{\max}^{1/2}) . \tag{3}
\end{align*}
\end{thm}

\begin{proof}
Clearly,
$\mathcal{C}_{\max}^{1/2}\subset
 \im\hat C_{\max}\cap H^{1/2}$.
If they are not equal, there is a vector
$y\in\im\hat C_{\max}\setminus\mathcal{C}_{\max}^{1/2}$,
and then $x=\gamma y\in H^{1/2}$ is non-zero,
$x\notin\gamma\mathcal C_{\max}^{1/2}$,
and is perpendicular to $\mathcal C_{\ext}$.
Arguing as in the proof of \tref{cald3},  
we arrive at a contradiction.

Let $B$ be an elliptic boundary condition
and $|s|\le1/2$.
Choose $\Lambda_0$ according to \tref{cald1}
and let $\Lambda\ge\Lambda_0$.
Write
\[
  B = \{ x+y+by \,:\, x\in F, y\in U\cap H^{1/2} \}
\]
as in \pref{bcbsrdata},
where $F\subset H_{>\Lambda}^{1/2}$
is of finite dimension,
$U\subset H_{\le\Lambda}$ is the orthogonal
complement of a subspace
$E\subset H_{\le\Lambda}^{1/2}$ of finite dimension,
and $b:U\to V=F^\perp\cap H_{>\Lambda}$ is $1/2$-smooth.
In particular,
\[
  B^{(s)} = \{ x+y+by \,:\, x\in F, y\in U^{(s)} \} ,
\]
where $U^{(s)}$ is the closure of $U\cap H^{1/2}$ in $H^s$
and, simultaneously, the annihilator of $E$ 
in $H_{\le\Lambda}^s$.
By \tref{cald1} we have, on the other hand,
\[
  \check{\mathcal C}_{\ext} =
  \{ u+v+Tv \,:\, u\in K_\Lambda, v\in H_{>\Lambda}^{-1/2} \} ,
\]
where $K_\Lambda\subset H_{\le\Lambda}^{1/2}$
and $T:H_{>\Lambda}^{-1/2}\to H_{\le\Lambda}^{1/2}$.
In particular,
$Q_{\le\Lambda}z\in H^{1/2}$
for any $z\in\check{\mathcal C}_{\ext}$.
We conclude that $B^{(s)}\cap\mathcal C_{\ext}^s$
is contained in $H^{1/2}$
and hence that
\begin{equation}\label{nulibc}
  B^{(s)} \cap \mathcal C_{\ext}^s
  = B \cap \mathcal C_{\ext}^{1/2} .
\end{equation}
By the above characterization of $B^{(s)}$,
$(B^{(s)},H_{>\Lambda}^s)$ is a Fredholm pair.
Now $I-C_{\ext}-Q_{\le\Lambda}$ is compact,
by \tref{cald2},
hence $(B^{(s)},\mathcal C_{\ext}^s)$
is a left-Fredholm pair, by \pref{bp4}.
By \tref{cald3}, 
we have 
\[
  (\mathcal C_{\ext}^s)^{\rm pol}
  = \gamma\im C_{\max}^{-s} ,
\]
where the superscript \lq pol' indicates the annihilator
of a subset of $H^s$ in $H^{-s}$.
Using \eqref{fgd}, we obtain
\begin{align*}
  (B^{(s)} + \mathcal C_{\ext}^s)^{\rm pol}
  &= (B^{(s)})^{\rm pol} 
   \cap (\mathcal C_{\ext}^s)^{\rm pol} \\
  &= (B^{(s)})^{\rm pol} \cap \gamma\im C_{\max}^{-s} \\
  &= \gamma\big( \gamma(B^{(s)})^{\rm pol} 
                \cap \im C_{\max}^{-s} \big) ,
\end{align*}
We also have
\[
  \im C_{\max}^{-s} \subset \im C_{\ext}^{-s}
  = \mathcal C_{\ext}^{-s}
  \subset \check{\mathcal C}_{\ext}
  \subset \check H .
\]
By the ellipticity of $B$,
\[
  \gamma (B^{(s)})^{\rm pol} \cap \check H
  \subset \gamma B^0 \cap \check H
  = B^a \subset H^{1/2} ,
\]  
where $B^0$ denotes the annihilator of $B$ in $H^{-1/2}$.
In conclusion,
\begin{equation}\label{defibc}
  (B^{(s)} + \mathcal C_{\ext}^s)^{\rm pol}
  = \gamma (B^a \cap \im C_{\max}^{-s})
  = \gamma (B^a \cap \mathcal C_{\max}^{1/2}) .
  \qedhere
\end{equation}
\end{proof}

\subsection{Some index formulas}\label{coth}
\tref{cmaxthm} and \cref{cor:wfredgen} 
have the following consequence:

\begin{thm}\label{windgen}
If $(d,V)$ is non-parabolic and $B$ is elliptic, then 
\begin{equation*}
  \ind_{\ext} D_B
  = \ind D_{B,\ext} 
  = \ind (\bar B,\mathcal C_{\ext}) ,
\end{equation*}
where $\bar B$ denotes the closure of $B$ in $H$.
\qed
\end{thm}

\begin{thm}\label{relindbsa}
If $(d,V)$ is non-parabolic and $B$ is elliptic, then
\[
  \ind D_{B,\ext} + \ind D_{B^a,\ext}
  = \dim (\mathcal C_{\ext}/\im C_{\max}) ,
\]
where $\im C_{\max}$ is the closure
of $\mathcal C_{\max}$ in $H$.
\end{thm}

\begin{proof}
Since $C_{\ext}-C_{\max}$ is compact in $H$
and $\mathcal C_{\ext}=\im C_{\ext}$,
we have
\begin{align*}
  \ind(\bar B,\mathcal C_{\ext})
  = \ind(\bar B,\im C_{\max}) 
     + \ind(\ker C_{\max},\mathcal C_{\ext}) ,
\end{align*}
by \pref{bp4}.
Since $\mathcal C_{\ext}=\gamma(\im C_{\max})^\perp$,
we have, by Theorem IV.4.8 in \cite{Ka},
\begin{align*}
  \ind(\bar B,\im C_{\max})
  &= - \ind(B^\perp,(\im C_{\max})^\perp) \\
  &= - \ind (\gamma B^\perp,\im C_{\ext}) 
  = - \ind D_{B^a,\ext} .
\end{align*}
Furthermore, 
$\im C_{\max} = \bar{\mathcal C}_{\max}
  \subset \mathcal C_{\ext}$,
hence
\[
  \ind(\ker C_{\max},\mathcal C_{\ext})
  = \dim (\mathcal C_{\ext}/\im C_{\max}) .
  \qedhere
\]
\end{proof}

\tref{windgen} implies the following index formula
of Agranovi\v{c}-Dynin type,
which corresponds to Theorem 23.1 in \cite{BW}.

\begin{thm}\label{relindadt}
If $(d,V)$ is non-parabolic, $B$ is elliptic,
and $\Lambda\in\R$, then
\begin{equation*}
   \ind D_{B,\ext}
   = \ind D_{\le\Lambda,\ext} 
     + \ind(\bar B,H_{>\Lambda}) .
\end{equation*}
\end{thm}

\begin{proof}
Since $C_{\ext}-Q_{>\Lambda}$ is compact, 
we can apply \tref{windgen} and \pref{bp4}
and get
\begin{equation*}
  \ind D_{B,\ext} = \ind(\bar B,\im C_{\ext})
  = \ind(\bar B,H_{>\Lambda}) 
    + \ind(H_{\le\Lambda},\im C_{\ext})
  \qedhere
\end{equation*}
\end{proof}

Note that in the notation of \pref{bcbsrdata},
\begin{equation}\label{relindadt2}
  \ind(\bar B,H_{>\Lambda})=\dim F-\dim E .
\end{equation}
In the corresponding form, 
the index formula in \tref{relindadt}
was also observed in \cite{BB}
(in the case of Dirac operators on smooth manifolds).

\begin{cor}[Discontinuity formula]\label{relindcor}
If $(d,V)$ is non-parabolic and $\Lambda\in\R$, then
\begin{equation*}
   \ind D_{\le\Lambda,\ext}
   = \ind D_{<\Lambda,\ext} 
     + \dim H_{\Lambda} .
   \qed
\end{equation*}
\end{cor}

In one of its versions,
the Cobordism Theorem for Dirac operators 
states that the index of the Dirac operator $D^+$ of a closed 
spin manifold $M$ of even dimension vanishes if $M$ bounds 
a compact spin manifold.
As an application of our results, 
we derive a general form of this.
Let $(d,V)$ be a Dirac-Schr\"odinger system.
Set
\begin{equation}\label{cothbc0}
  H^\pm :=\{x\in H \,:\, i\gamma x = \pm x\} .
\end{equation}
Since $\gamma$ and $A$ anti-commute,
\begin{equation}\label{cothbc}
  B^\pm := H^\pm \cap \check H = H^\pm \cap H^{1/2} .
\end{equation}
Since $H^+$ is the orthogonal complement of $H^-$
in $H$, we conclude that $B^+$ and $B^-$ are mutually
adjoint elliptic boundary conditions.

\begin{nameit}{Cobordism Theorem}\label{cothm}
If the system $(d,V)$ is of Fredholm type, 
then the restriction $A^+:H_A^+\to H^-$ of $A=A_0$
satisfies $\ind A^+ = 0$.
\end{nameit}

\begin{proof}
Since $(d,V)$ is of Fredholm type, $\ker C_{\ext}$
is an elliptic self-adjoint boundary condition,
by \cref{fredtosa}.
Now \tref{bcsrsa} implies $\ind A^+=0$  
(compare also \cref{bcsrsacor}).
\end{proof}

We now consider Dirac-Schr\"odinger systems
together with a boundary value problem
which models the decomposition of a manifold $M$
into two pieces $M_1$ and $M_2$ along a closed 
hypersurface $N=M_1\cap M_2$.
This requires the transmission boundary condition
for sections of bundles over $M$
and Dirac-Schr\"odinger operators acting on them;
compare Example \ref{bcsrsaexa}.

Let $(d_1,V_1)$ and $(d_2,V_2)$ be Dirac-Schr\"odinger
systems with the same initial Hilbert space $H$ at $t=0$
(after some appropriate identification).
Suppose that, at $t=0$,
\begin{equation}\label{decoinit}
  A_{1,0} = - A_{2,0} =: A
  \quad\text{and}\quad
  \gamma_{1,0} = - \gamma_{2,0} =: \gamma .
\end{equation}
We consider the Dirac-Schr\"odinger system 
$(d,V)=(d_1,V_1)\oplus(d_2,V_2)$
with the boundary condition 
\begin{equation}\label{decobc}
  B = \{ (x,x) \,:\, x\in H^{1/2} \} ,
\end{equation}
where we use $A$ to define $H^{1/2}$.
We already observed in Example \ref{bcsrsaexa}
that $B$ is elliptic and self-adjoint.
The Calder\'on space of $d$ is the direct sum 
of the Calder\'on spaces of $d_1$ and $d_2$,
\begin{equation}\label{decocald}
  \check{\mathcal C}_{\ext}
  = \check{\mathcal C}_{1,\ext} 
    \oplus \check{\mathcal C}_{2,\ext}
  \quad\text{and}\quad
  \mathcal C_{\ext}
  = \mathcal C_{1,\ext} \oplus \mathcal C_{2,\ext} .
\end{equation}
We then arrive at the following index formula of Bojarski type.

\begin{thm}\label{botythm}
If $(d_1,V_1)$ and $(d_2,V_2)$ are non-parabolic,
then $(d,V)$ is non-parabolic,
$(\mathcal C_{1,\ext},\mathcal C_{2,\ext})$
is a Fredholm pair in $H$, and
\[
  \ind D_{B,\ext} = 
  \ind(\mathcal C_{1,\ext},\mathcal C_{2,\ext}) .
\]
\end{thm}

\begin{proof}
The first assertion is clear.
By \tref{cmaxthm},
$(C_{1,\ext},H_\le)$ is a Fredholm pair,
where we use spectral projections and spaces
associated to $A$.
By \tref{cald2}, 
$C_{2,\ext}-Q_\le$ is a compact operator.
Hence $(\mathcal C_{1,\ext},\mathcal C_{2,\ext})$
is a Fredholm pair, by \pref{bp4}.
As for the index formula, 
we note that
\begin{align*}
  \bar B \cap \mathcal C_{\ext}
  &= \{(x,x) \in H\oplus H\,:\, \text{
     $x\in\mathcal C_{1,\ext}$ 
     and $x\in\mathcal C_{2,\ext}$} \} \\
  &\cong \mathcal C_{1,\ext} \cap \mathcal C_{2,\ext} \\
  B^\perp \cap \mathcal C_{\ext}^\perp
  &= \{(x,-x) \in H\oplus H\,:\, \text{
     $x\perp\mathcal C_{1,\ext}$ 
     and $x\perp\mathcal C_{2,\ext}$} \} \\
  &\cong (\mathcal C_{1,\ext} + \mathcal C_{2,\ext})^\perp .
\end{align*}
Therefore
\[
  \ind D_{B,\ext}
  = \ind(\bar B,\mathcal C_{\ext})
  =  \ind(\mathcal C_{1,\ext},\mathcal C_{2,\ext}) .
  \qedhere
\]
\end{proof}

Using \tref{relindadt} and \cref{relindcor},
we get a splitting formula for the index,
which generalizes Theorems 23.3 of \cite{BW}
and 4.3 of \cite{BL1}.

\begin{thm}[Splitting formula]\label{decothm}
If $(d_1,V_1)$ and $(d_2,V_2)$ are non-parabolic,
$B_1$ is an elliptic boundary condition
with respect to $A$, 
and $B_2$ is an elliptic boundary condition
with respect to $-A$, then
\begin{align*}
  \ind D_{B,\ext} 
  = \ind D_{1,B_1,\ext} &+ \ind D_{2,B_2,\ext}  \\
  &- \ind(H_>,\bar B_1) - \ind(H_\le,\bar B_2) .
\end{align*}
In particular, if $B_1$ is any elliptic boundary 
condition with respect to $A$
and $B_2=B_1^\perp\cap H^{1/2}$, then
\[
  \ind D_{B,\ext} 
  = \ind D_{1,B_1,\ext} + \ind D_{2,B_2,\ext} .
\]
\end{thm}

\begin{proof}
By \tref{relindadt} and \cref{relindcor},
\begin{align*}
  \ind D_{B,\ext}
  &= \ind D_{1,\le,\ext} 
     + \ind D_{2,\ge,\ext}
     + \ind (\bar B,H_>\oplus H_<) \\
  &= \ind D_{1,\le,\ext} 
     + \ind D_{2,>,\ext} \\
  &= \ind D_{1,B_1,\ext} - \ind(H_>,\bar B_1)
     + \ind D_{2,B_2,\ext} - \ind(H_\le,\bar B_2).
\end{align*}
If $B_2=B_1^\perp\cap H^{1/2}$,
then the second and last term on the right hand side cancel
each other.
\end{proof}

Besides modeling the case mentioned in the beginning
of this section,
the above results also apply to a Dirac-Schr\"odinger 
system $d$ defined over the whole real line, 
decomposed into pieces $d_1:=d|\R_+$ and $d_2:=d|\R_-$,
where we need to turn the latter into a Dirac-Schr\"odinger 
system over $\R_+$ in the appropriate and obvious way.

\newpage
\section{Supersymmetric systems}\label{supsym}

Our treatment so far does not allow to treat
the usual index theorems since $D_{0,c}$ is symmetric.
To adjust this we formulate a further axiom,
introducing a supersymmetry, i.e. an involution
which anticommutes with $D_{\max}$.

\begin{axia}\label{ax:alpha}
There is a section
\[
  \alpha \in
  \Lip_{\loc}(\R_+,\mathcal{L}(H)) \cap
  L_{\loc}^{\infty}(\R_+,\mathcal{L}(H_{A})),
\]
such that the following relations hold:
\begin{alignat}{2}
    \alpha_t = \alpha_t^* &= \alpha_{t}^{-1} \quad
    &&\mbox{on $H_t$} ,
    \tag{1} \\
    \alpha_t\gamma_t + \gamma_t \alpha_t &= 0
    &&\mbox{on $H_{A}$} ,
    \tag{2} \\
    \lbrack\partial,\alpha\rbrack &= 0
    &&\mbox{on $\Lip_{\loc}(\mathcal H)$} ,
    \tag{3} \\
    \lbrack A_{t},\alpha_{t}\rbrack &= 0
    &&\mbox{on $H_{A}$} , 
    \tag{4} \\
    \alpha_t V_t + V_t \alpha_t &= 0
    &&\mbox{on $H_{t}$} .
    \tag{5} 
\end{alignat}
\end{axia}

A {\em supersymmetric Dirac-Schr\"odinger system} 
is a Dirac-Schr\"odinger system $(d,V)$
together with a supersymmetry $\alpha$
as in Axiom \ref{ax:alpha}.

Let $(d,V,\alpha)$ be a supersymmetric Dirac-Schr\"odinger
system.
Then we have, for each $t\ge0$, an orthogonal decomposition
\begin{equation}\label{supsympm}
  H_t = H_{t}^+ \oplus H_{t}^- ,
  \quad
  H_{t}^\pm := \{ x\in H \,:\, \alpha_tx=\pm x \} .
\end{equation}
Since $A_t$ commutes with $\alpha_t$,
we get an associated decomposition
\begin{equation}\label{supsympma}
  H_A = H_{A,t}^+ \oplus H_{A,t}^- ,
  \quad
  H_{A,t}^\pm := H_A \cap H_{t}^\pm ,
\end{equation}
which is orthogonal
with respect to the graph norm of $A_t$
and such that $A_t$ maps $H_{A,t}^\pm$ to $H_{t}^\pm$.
There are analogous decompositions of the associated
Sobolev and function spaces.
We also have
\begin{equation}\label{supsymdodd}
     \alpha D + D\alpha = 0
\end{equation}
on $\mathcal{L}_{\loc}(e^0)$.
It follows that $D$ is an {\em odd operator}, that is,
maps locally Lipschitz sections of $\mathcal H^\pm$
to locally essentially bounded measurable sections
of $\mathcal H^\mp$.
We let $D^\pm$ be the corresponding parts of $D$
so that $D$ is represented by the matrix
\begin{equation}\label{supsymdodd4}
  \begin{pmatrix} 0 & D^- \\ D^+ & 0 \end{pmatrix}
\end{equation}
with respect to the above decomposition
of $\mathcal{L}_{\loc}(e^0)$.
We obtain orthogonal decompositions
\begin{equation}\label{supsymmw}
  \mathcal D_{\max}
  = \mathcal D_{\max}^+ \oplus \mathcal D_{\max}^-
  \quad\text{and}\quad
  W = W^+ \oplus W^- ,
\end{equation}
and $D_{\max}$ and $D_{\ext}$ are odd operators with
respect to these with corresponding parts
$D_{\max}^\pm$ and $D_{\ext}^\pm$, respectively.
Since $\ker D_{\max}$ and $\ker D_{\ext}$
are $\alpha$-invariant, we have
\begin{equation}\label{supsymker}
\begin{split}
  \ker D_{\max}
  &= \ker D_{\max}^+ \oplus \ker D_{\max}^- , \\
  \ker D_{\ext}
  &= \ker D_{\ext}^+ \oplus \ker D_{\ext}^- ,
\end{split}
\end{equation}
respectively.
Since $\mathcal R$ commutes with $\alpha_0$,
$\check C_{\max}$ and $\check C_{\ext}$
are $\alpha_0$-invariant
and hence decompose accordingly,
\begin{equation}\label{supsymcal}
  \check{\mathcal C}_{\max}
  = \check{\mathcal C}_{\max}^+ 
    \oplus \check{\mathcal C}_{\max}^-  
  \quad\text{and}\quad
  \check{\mathcal C}_{\ext}
  = \check{\mathcal C}_{\ext}^+ 
    \oplus \check{\mathcal C}_{\ext}^- .
\end{equation}
We are interested in boundary value problems
that are compatible with the supersymmetry.
That is, we require that boundary conditions $B$
are $\alpha_0$-invariant,
and then we have a decomposition $B=B^+\oplus B^-$
as above.
In other words,
we pose the boundary conditions separately
for the $+$ and $-$ parts of the elements
in the corresponding domains and get
corresponding domains and operators
\begin{equation}\label{supsymbpm}
  D_{B^\pm,\ext}^\pm : W_{B^\pm}^\pm
  \to  L^2(\mathcal H^\mp) ,
\end{equation}
and similarly for $D$ and $D_{\max}$.

\begin{prop}\label{susyind}
Let $(d,V)$ be a non-parabolic supersymmetric 
Dirac-Schr\"odinger system
with supersymmetry $\alpha$
and $B$ be an $\alpha_0$-invariant elliptic
boundary condition.
Then
\begin{equation*}
  \ind D_{B,\ext}
  = \ind D_{B^+,\ext}^+ + \ind D_{B^-,\ext}^- .
  \qed
\end{equation*}
\end{prop}

If $C$ is an $\alpha_0$-invariant subspace of $H$,
then $C^\perp$ and $\gamma C$ are invariant under
$\alpha_0$ as well and we have 
\begin{equation}\label{susyinv}
  (\gamma C^\perp)^\pm
  = \gamma (C^{\perp,\mp})
  = \gamma ((C^\mp)^\perp\cap H^\mp) .
\end{equation}
In particular, from \tref{cald3},
\begin{equation}\label{susyinv2}
   (\im C_{\max})^\mp
   = (\gamma\mathcal C_{\ext}^\perp)^\mp
   = \gamma(\mathcal C_{\ext}^{\perp,\pm}) . 
\end{equation}
If $P$ is a projection in $H$,
then $\ker P$ and $\im P$ are invariant under $\alpha_0$
if and only if $[P,\alpha_0]=0$,
and then $P$ decomposes as
\begin{equation}\label{supsymp}
  P = \frac{1}{2}(\alpha + I) P
      + \frac{1}{2}(\alpha - I) P
    =: P^{+} + P^{-} .
\end{equation}
Clearly $[P_\gamma,\alpha_0]=0$
if $[P,\alpha_0]=0$, and then
\begin{equation}\label{supsymppm}
  P_{\gamma}^\pm
  = \gamma^*(I^\mp-P^{*,\mp})\gamma .
\end{equation}
The following index formulas are immediate from
Theorems \ref{windgen}, \ref{relindbsa},
and \ref{relindadt}. 

\begin{thm}\label{susyind2}
Let $(d,V)$ be a non-parabolic supersymmetric 
Dirac-Schr\"odinger system
with supersymmetry $\alpha$
and $B$ be an $\alpha_0$-invariant elliptic
boundary condition.
Then
\begin{equation*}
\begin{split}
  \ind D_{B^+,\ext}^+
  &= \ind (\bar B^+,\mathcal C_{\ext}^+) \\
  &= \ind D_{H_\le^+,\ext}^+ 
      + \ind (\bar B^+,H_>^+) , \\
  \ind D_{B^+,\ext}^+ + \ind D_{B^{a,-},\ext}^-
  &= \dim(\mathcal C_{\ext}^+/\im C_{\max}^+) \\
  &= \dim (\mathcal C_{\ext}^-/\im C_{\max}^-) .
  \qed
\end{split}
\end{equation*}
\end{thm}

Recall the setup in Theorems \ref{botythm}
and \ref{decothm}.
Let $\alpha_1$ and $\alpha_2$ be supersymmetries
of the Dirac-Schr\"odinger systems $(d_1,V_1)$ 
and $(d_2,V_2)$, respectively,
that agree at $t=0$.
Consider the Dirac-Schr\"odinger system 
$(d,V)=(d_1,V_1)\oplus(d_2,V_2)$
with the induced supersymmetry $(\alpha_1,\alpha_2)$.
The boundary condition $B$ from \eqref{decobc}
is $(\alpha_1,\alpha_2)$-invariant with
\begin{equation}\label{susydecobc}
  B^\pm = \{ (x,x) \,:\, x \in H^{\pm} \} \cap H^{1/2} .
\end{equation}
We also have
\begin{equation}\label{susydecocald}
  \check{\mathcal C}_{\ext}^\pm = 
    \check{\mathcal C}_{1,\ext}^\pm 
    \oplus \check{\mathcal C}_{2,\ext}^\pm .
\end{equation}
Arguing as in the proofs of Theorems 
\ref{botythm} and \ref{decothm} we get the 
following index formulas.

\begin{thm}\label{susyind4}
Assume that $(d_1,V_1)$ and $(d_2,V_2)$ are non-parabolic.
Then
\[
  \ind D_{B^+,\ext}^+
  = \ind (\mathcal C_{1,\ext}^+,\mathcal C_{2,\ext}^+) .
\]
If $B_1$ is any $\alpha_1$-invariant elliptic
boundary condition for $d_1$ and $B_2$ any
$\alpha_2$-invariant elliptic boundary condition
for $d_2$, then
\begin{align*}
  \ind D_{B^+,\ext}^+
  = \ind D_{1,B_1^+,\ext}^+ &+ \ind D_{2,B_2^+,\ext}^+ \\
  &- \ind (H_>^+, \bar B_1^+) - \ind (H_\le^+,\bar B_2^+) . 
\end{align*}
In particular,
if $B_1$ is any $\alpha_1$-invariant elliptic
boundary condition for $d_1$ 
and $B_2=B_1^\perp\cap H^{1/2}$, then
\[
  \ind D_{B^+,\ext}^+
  = \ind D_{1,B_1^+,\ext}^+ + \ind D_{2,B_2^+,\ext}^+  .
  \qed
\]
\end{thm}

\newpage
\section{Manifolds with boundary}
\label{nsb}

In this last chapter, 
we explain how our results can be applied 
to obtain formulas for the index of Dirac type operators
on manifolds with boundary.
Such formulas are well known in the case of compact 
manifolds with smooth boundary 
and Dirac operators with smooth coefficients,
see for instance \cite{BW}.
However, in applications one often faces the problem
that the boundary of the manifold is not smooth or
that the coefficients of the operator are not smooth.
We will encounter such a situation in a forthcoming article
on $L^2$-index formulas on manifolds with finite volume 
and pinched negative curvature 
in which we extend the results of \cite{BB2}.
Here we concentrate on a rather general case
which sets the stage for the applications we have in mind,
but should also be useful in other situations.

\subsection{The geometric setup}\label{nsbgsu}
Let $M$ be a $C^{1,1}$ manifold with compact boundary $N=\partial M$
and with a Lipschitz continuous Riemannian metric. 
Let $E\to M$ be a $C^{0,1}$ Hermitian vector bundle
and $D$ be a differential operator on $E$ of order one
with $L^\infty_{\loc}$ coefficients.
Then we obtain a linear operator
\begin{equation}\label{nsbd}
  D: \Lip_{\loc}(M,E)\rightarrow L^{\infty}_{\loc}(M,E) .
\end{equation}
Let $\Lip_{0,c}(M,E)$ be the space of Lipschitz sections of $E$ 
with compact support in $M$,
which vanish along the boundary $N$,
and set $D_{0,c}:=D|\Lip_{0,c}(M,E)$,
considered as an unbounded operator on $L^2(M,E)$.
We assume that $D_{0,c}$ is symmetric, that is,
\begin{equation}\label{nsbfsa}
  (D\sigma_1,\sigma_2)_{L^2(M,E)} 
  = (\sigma_1,D\sigma_2)_{L^2(M,E)} 
\end{equation}
for all $\sigma_1,\sigma_2\in\Lip_{0,c}(M,E)$.
We let $D_{\min}$ be the closure of $D_{0,c}$
and $D_{\max}$ be the adjoint of $D_{0,c}$ 
in $L^2(M,E)$.
We denote by $\mathcal D_{\min}$ and $\mathcal D_{\max}$
the domains of $D_{\min}$ and $D_{\max}$, respectively.

\begin{axia}\label{ax:nsb}
There is a Lipschitz function $\rho:M\to\R_+$ 
and a constant $r>0$ such that $N=\rho^{-1}(0)$
and $O:=\rho^{-1}([0,r))$ is relatively compact in $M$.
Moreover, there is a Dirac-Schr\"odinger system 
$(d,V)=(\mathcal H,\partial,A,\gamma,V)$
with Lipschitz coefficients, and a unitary isomorphism 
$U: L^2(O,E) \to L^2(\mathcal H|[0,r))$
such that
\begin{enumerate}
\item\label{nsbnatur}
$U((\varphi\circ\rho)\sigma)=\varphi U\sigma$
for all $\sigma\in L^2(O,E)$ 
and $\varphi\in L^\infty_{\loc}(\R_+)$.
\item\label{nsbmaxmin}
$(1-\varphi\circ\rho)\sigma\in\mathcal D_{\min}$
for all $\sigma\in\mathcal D_{\max}$
and $\varphi\in\Lip(\R_+)$ with compact support in $[0,r)$
and equal to one close to zero.
\item\label{nsbdense}
$U(\Lip_c(O,E))$ is contained and dense 
in $\mathcal L_c(\mathcal H|[0,r))$
with respect to the graph norm of $D^d$.
\item\label{nsbdense0}
$U(\Lip_{0,c}(O,E))$ is contained and dense 
in $\mathcal L_{0,c}(\mathcal H|[0,r))$
with respect to the graph norm of $D^d$.
\item\label{nsbtwine}
$D^d(U\sigma)=U(D\sigma)$ for all $\sigma\in\Lip_{c}(O,E)$.
\end{enumerate}
\end{axia}

Here $\mathcal L_c(\mathcal H|[0,r))$ denotes the space
of sections in $\mathcal L_c(\mathcal H)$ 
with compact support in $[0,r)$.
We also use a superscript $d$ to distinguish quantities 
belonging to $(d,V)$ if necessary.

\begin{rem}\label{nsbrem}
Axiom \ref{ax:nsb} is tailored to fit the geometric 
examples which we will discuss in our next article,
notably the case of Dirac-Schr\"odinger operators
over the ends of complete Riemannian manifolds
with finite volume 
and pinched negative sectional curvature, see \cite{BB2}.
\end{rem}

For $\sigma\in\Lip_{\loc}(M,E)$, let 
$\mathcal R\sigma:=\mathcal R^d(U((\varphi\circ\rho)\sigma))$,
where $\mathcal R^d$ denotes the restriction map of $d$
and $\varphi\in\Lip(\R_+)$ has compact support in $[0,r)$
and is equal to one close to zero.
By Axiom \ref{ax:nsb}.\ref{nsbnatur} above, 
$\mathcal R\sigma$  does not depend 
on the choice of $\varphi$.
As before, we also write $\sigma(0)=\mathcal R\sigma$.
Using \eqref{eq:intpart1gen}, \eqref{nsbfsa}, 
and Axiom \ref{ax:nsb}.\ref{nsbtwine} we get
\begin{equation}\label{nsbpi}
  (D\sigma_1,\sigma_2) - (\sigma_1,D\sigma_2)
  = \omega^d (\sigma_1(0),\sigma_2(0)) 
  =: \omega (\sigma_1(0),\sigma_2(0)) ,
\end{equation}
for all $\sigma_1,\sigma_2\in\Lip(M,E)$
with compact support.

\begin{lem}\label{nsbdommax}
Suppose $\sigma\in L^2(M,E)$ has compact support in $O$.
Then $\sigma\in\mathcal D_{\max}$ if and only if
$U\sigma\in\mathcal D^d_{\max}$, and then
$D^d_{\max}(U\sigma) = UD_{\max}\sigma$.
\end{lem}

\begin{proof}
We need only to test against Lipschitz sections of $E$
with compact support in $O$ and vanishing along $N$
respectively Lipschitz sections of $\mathcal H$ 
with compact support in $[0,r)$ and vanishing at $0$.
To these, \eqref{nsbdense0} and \eqref{nsbtwine} 
of Axiom \ref{ax:nsb} apply,
and the lemma follows.
\end{proof}

Using Axiom \ref{ax:nsb} and \lref{nsbdommax},
we get the following characterization of 
the maximal domain $\mathcal D_{\max}$.

\begin{cor}\label{nsbdommax2}
For any $\varphi\in\Lip(\R_+)$ with compact support in $[0,r)$
and equal to one close to zero and any $\sigma\in L^2(M,E)$,
\[
  \sigma\in\mathcal D_{\max} \Longleftrightarrow 
  \text{$\varphi U\sigma \in \mathcal D^d_{\max}$ and 
        $(1-(\varphi\circ\rho))\sigma \in \mathcal D_{\min}$} .
  \qed
\]
\end{cor}

\begin{prop}[Regularity]\label{nsbreg}
The maximal domain $\mathcal D_{\max}$ satisfies:
\begin{enumerate}
\item
$\Lip_c(M,E)$ is dense in $\mathcal D_{\max}$.
\item
The restriction map on $\Lip_c(M,E)$ extends to \\
a continuous surjective 
map $\mathcal R:\mathcal D_{\max}\to\check H$.
\item
For $\sigma_1,\sigma_2\in\mathcal D_{\max}$ we have
\begin{equation*}
  (D_{\max}\sigma_1,\sigma_2) - (\sigma_1,D_{\max}\sigma_2)
  = \omega (\sigma_1(0),\sigma_2(0)) .
\end{equation*}
\end{enumerate}
\end{prop}

\begin{proof}
Apply \pref{thm:diprops}, \cref{nsbdommax2}, 
and Axiom \ref{ax:nsb}.
\end{proof}

For a boundary condition $B\subset\check H$, we set
\begin{align}
  \mathcal D_{B,\max}
  :&= \{ \sigma\in\mathcal D_{\max} 
         \,:\, \mathcal R\sigma\in B \} , \\
  D_{B,\max} :&= D_{\max}|\mathcal D_{B,\max} \notag .
\end{align}
Then $D_{B,\max}$ is closed 
with adjoint $D_{B^a,\max}$,
see Section \ref{boco}.

 \subsection{Fredholm properties}\label{nsbfp}
We now discuss Fredholm properties of and index formulas 
for the operators $D_B$. 
As in the case of Dirac-Schr\"odinger systems,
we need the non-parabolicity condition of the third named author:

\begin{axia}\label{geononpara} 
We say that $D$ is {\em non-parabolic} if for any compact subset 
$K\subset M$ there is a positive constant $C_K$ such that
any $\sigma\in\mathcal D_{\max}$ satisfies
\begin{equation*}
  \|\sigma\|_{L^2(K)} \le C_K \left(
  \|\mathcal R\sigma\|_{\check H} +
  \|D_{\max}\sigma\|_{L^2(M,E)} \right) .
\end{equation*}
\end{axia}

Assume from now on that $D$ is non-parabolic.
Let $W$ be the completion of $\mathcal D_{\max}$ 
with respect to the norm appearing on the right hand side
of the equation in Axiom \ref{geononpara}. 
There is the following analogue of \lref{lem:west3gen}.

\begin{lem}\label{lem:west3geo} 
If $D$ is non-parabolic, then we have:
\begin{enumerate}
\item
The restriction map $\mathcal{R}$ and $D_{\max}$ extend to 
continuous maps 
\[
  \mathcal R_{\ext}: W \to \check H
  \quad\text{and}\quad
  D_{\ext}: W \to L^2(M,E),
\]
respectively;
$\mathcal R_{\ext}$ induces an isometry from $\ker D_{\ext}$
into $\check H$.
\item
If $\psi\in \Lip_{c}(M)$ and $\sigma\in W$,
then $\psi\sigma\in\mathcal{D}_{\max}\subset W$.
Moreover, there is a constant $C_{\psi}$ such that
\[
  ||\psi\sigma||_{D_{\max}}
  \le C_\psi ||\sigma||_W .
\]
In particular, $W$ can be viewed as
a space of locally integrable functions
and $W\cap L^{2}(M,E) = \mathcal{D}_{\max}$.
\item
$W=\mathcal{D}_{\max}$ if and only if
there is a constant $C$ such that
\begin{equation*}
  ||\sigma||_{L^{2}(\mathcal{H})}
  \le C ||\sigma||_W
  \quad\text{for all $\sigma\in\Lip_c(M,E)$} .
\end{equation*}
\end{enumerate}
\end{lem}

\begin{proof}
(1) and (3) are clear.
As for (2), use \lref{nsbdommax} and argue 
as in the proof of (2) of \lref{lem:west3gen}.
\end{proof}

Similarly, there is an analogue of \lref{wpclem}:

\begin{lem}\label{wpclemgeo}
Let $V$ be a bounded subset of $W$.
Then $V$ is precompact if and only if
$D_{\ext}(V)\subset L^2(M,E)$
and $Q_{\ge}\mathcal R_{\ext}(V)\subset\check H$
are both precompact.
\end{lem}

\begin{proof}
It is easy to adapt the arguments in the proof
of \lref{wpclem} to the present situation.
\end{proof}

For a boundary condition $B\subset\check H$, set
\begin{equation}\label{bcwbgeo}
  W_B := \{\sigma\in W \,:\, \sigma(0) \in B \}
  \quad\text{and}\quad
  D_{B,\ext} := D_{\ext}|W_B .
\end{equation}
We arrive at the following generalization
of \tref{thm:wfredgen}, \cref{cor:l2ind},
and \pref{injesu}

\begin{thm}\label{thm:wfreddirac} 
Assume that $D$ is non-parabolic and that $B$ is regular.
Then $D_{B,\ext}:W_{B}\to L^{2}(\mathcal{H})$
is a left-Fredholm operator
with $(\im D_{B,\ext})^{\perp} = \ker D_{B^a,\max}$
and extended index
\begin{equation*}
  \ind_{\ext}D_B:= \ind D_{B,\ext}
  = \dim\ker D_{B,\ext} - \dim\ker D_{B^a,\max} .
\end{equation*}
If $B$ is elliptic, then the kernels of $D_B$ 
and $D_{B^a}$ have finite dimension
and the {\em $L^{2}$-index} of $D_{B}$ is well defined,
\begin{equation*}
    L^{2}\text{-}\ind D_{B}
    := \dim\ker D_{B} - \dim\ker D_{B^a} .
\end{equation*}
Moreover, there is $\Lambda_{0}\ge0$ such that
$D_{<-\Lambda,\ext}$ is injective
and $D_{\le\Lambda,\ext}$ is surjective
for all $\Lambda\ge\Lambda_{0}$.
\qed
\end{thm}

We define Calder\'on spaces and projections
as in the case of Dirac-Schr\"odinger systems,
see \eqref{caldercheck}, \eqref{eq:sols},
and Definition \ref{caldprodef}.
If $B$ is a regular boundary condition,
then $\mathcal R$ induces isomorphisms
\begin{equation}\label{kerdbgeo}
\begin{split}
  \ker D_{B,\max}
  &\cong B\cap\check{\mathcal C}_{\max}
  = B\cap\mathcal C_{\max}^{1/2} , \\
  \ker D_{B,\ext}
  &\cong B\cap\check{\mathcal C}_{\ext}
  = B\cap\mathcal C_{\ext}^{1/2} .
\end{split}
\end{equation}
As before,
we write $\mathcal{C}_{\max}$ and $\mathcal{C}_{\ext}$
instead of $\mathcal{C}^{0}_{\max}$ and
$\mathcal{C}^{0}_{\ext}$, respectively.
We have the following analogue of \cref{cor:wfredgen}:

\begin{cor}\label{cor:wfredgeo}
Assume that $D$ is non-parabolic and 
that $B$ is elliptic.
Then $D_{B,\ext}$
is a Fredholm operator with 
$(\im D_{B,\ext})^{\perp} = \ker D_{B^a,\max}$
and index
\begin{align*}
  \ind D_{B,\ext}
  &= \dim B\cap\mathcal C_{\ext}^{1/2}
  - \dim B^\perp\cap\gamma\mathcal C_{\max} ^{1/2}. \\
  &= \dim B\cap\mathcal C_{\ext}
  - \dim B^\perp\cap\gamma\mathcal C_{\max} .
  \qed
\end{align*}
\end{cor}

It is a routine matter to check that the arguments
developed in Chapter \ref{calderon} also work
under Axioms \ref{ax:nsb} and \ref{geononpara}
imposed here;
hence all the results obtained there have their
analogues here.
We arrive at the following version of Theorems 
\ref{cald2}, \ref{cald3}, and \ref{cmaxthm}.

\begin{thm}\label{nsbcald}
Assume that $D$ is non-parabolic. Then:
\begin{enumerate}
\item
The Calder\'on projections
$C_{\ext}$ and $C_{\max}$ are elliptic 
with $C_{\max}=C_{\ext,\gamma}$.
\item
$C_{\max}-Q_>$ and $C_{\ext}-Q_>$ are compact
in $H^s$ for all $|s|\le1/2$.
\item
If $B$ is an elliptic boundary condition,
then $(\bar B,\mathcal C_{\ext})$ 
is a Fredholm pair in $H$ and
\[
  \bar B\cap \check{\mathcal C}_{\ext}
  = B \cap \mathcal C_{\ext}^{1/2}
  \quad\text{and}\quad
  (\bar B+\mathcal C_{\ext})^\perp
  = B^\perp \cap \mathcal\gamma C_{\max}^{1/2} .
  \qed
\] 
\end{enumerate}
\end{thm}

With the same arguments as in Chapter \ref{calderon},
we get the analogues of the index formulas in
Theorems \ref{windgen}, \ref{relindbsa}, 
and \ref{relindadt}:

\begin{thm}\label{nsbind2}
Assume that $D$ is non-parabolic and that $B$
is an elliptic boundary condition.
Then
\begin{align*}
  \ind D_{B,\ext}
  &= \ind (\bar B,\mathcal C_{\ext}) \\
  &= \ind D_{H_\le,\ext}
      + \ind (H_>,\bar B) , \\
  \ind D_{B,\ext} + \ind D_{B^{a},\ext}
  &= \dim(\mathcal C_{\ext}/\im C_{\max}).
  \qed
\end{align*}
\end{thm}

\begin{rem}\label{nsbanalog}
The further results from Chapter \ref{calderon}
and Chapter \ref{supsym} are consequences of 
the results on the Calder\'on projections
and the index formulas from Theorems 
\ref{windgen}, \ref{relindbsa},
and \ref{relindadt}.
Therefore they have their exact analogs here,
and we refrain from repeating the 
corresponding statements.
\end{rem}

\newpage
\begin{appendix}\section{Fredholm pairs}\label{fredpair}

T. Kato has developed the notion of
{\em Fredholm pairs of closed subspaces},
cf. \cite[Ch.IV, Section 4]{Ka}.
Consider a Banach space $E$
and a pair of closed subspaces $F$ and $G$.
Introduce {\em nullity} and {\em deficiency}
of the pair $(F,G)$,
\begin{align}
     \nuli(F,G) :&= \dim(F\cap G) , \label{fgnull} \\
     \defi(F,G) :&= \codim(F + G) , \label{fgdef}
\end{align}
and recall that $\text{def\,}(F,G)<\infty$
implies that $F+G$ is closed.
We say that the pair $(F,G)$ is a {\em left-}
or {\em right-Fredholm pair}, respectively, if
\begin{equation}\label{fgclosed}
    \text{$F + G$ is closed}
\end{equation}
and
\begin{equation}\label{fgfred}
   \text{null\,}(F,G) < \infty
   \quad\text{or}\quad
   \text{def\,}(F,G) <\infty ,
\end{equation}
respectively.
We say that $(F,G)$ is a {\em semi-Fredholm pair}
if it is a left- or right-Fredholm pair,
and that it is a {\em Fredholm pair}
if it is a left- and right-Fredholm pair.
For any semi-Fredholm pair $(F,G)$, its {\em index},
\begin{equation}\label{fgind}
     \ind (F,G) := {\rm null\,}(F,G) - {\rm def\,}(F,G) ,
\end{equation}
is well defined as an extended real number.
The index of $(F,G)$ is a rough measure
of the non-complementarity of $F$ and $G$ in $E$.

Let $E'$ be the dual space of $E$
and $F^0,G^0\subset E'$ be the annihilators
(or polar sets) of $F$ and $G$, respectively.
By \cite[Theorem IV.4.8]{Ka},
$F^0+G^0$ is closed if and only if $F+G$ is closed,
\begin{align}
   (F\cap G)^0 &= F^0 + G^0 ,
   \quad
   (F+G)^0 = F^0\cap G^0 , \label{fgd} \\
   {\rm null\,}(F^0,G^0) &= {\rm def\,}(F,G) ,
   \quad
   {\rm def\,}(F^0,G^0) = {\rm null\,}(F,G) .
\end{align}

For Banach spaces $E_{1}, E_{2}$
and an operator $T\in\mathcal L(E_1,E_2)$,
we recover the Fredholm properties
of $T$ by considering
\begin{equation}\label{tefg}
  E = E_1\times E_2 , \quad
  F = E_1\times\{0\} , \quad
  G = \graph T .
\end{equation}
To that end we note that $F+G$ is closed in $E$
if and only if $\im T$ is closed in $E_2$
and that the canonical inclusions $E_1\to E$
and $E_2\to E$ induce isomorphisms
\begin{equation}\label{tfred}
  \ker T \cong F\cap G
  \quad \text{and}\quad
  \coker T \cong E/(F+G) .
\end{equation}
In particular, if $T$ is semi-Fredholm,
then the index of $T$ is
\begin{equation}\label{tind}
  \ind T = \dim\ker T - \dim\coker T  = \ind (F,G) ,
\end{equation}
where $F$ and $G$ are as above.
Next we quote a criterion for left-Fredholmness of $T$
which is used several times in this work;
for a proof,
see for example \cite[Proposition 19.1.3]{Ho}.

\begin{lem}\label{trfred}
The following conditions are equivalent:
\begin{enumerate}
\item
$T\in\mathcal L(E_1,E_2)$ is a left-Fredholm operator.
\item
If $(x_{n})$ is a bounded sequence in $E_{1}$
with $(Tx_{n})$ convergent in $E_2$,
then $(x_{n})$ possesses a convergent subsequence.
\qed
\end{enumerate}
\end{lem}

Traditionally, the results on Fredholm pairs we have mentioned
are applied to subspaces with topological complements, i,e. 
to pairs of spaces of the form $F = \im P$, $G = \im Q$,
where $P,Q$ are projections (continuous idempotents) in $E$.
We need the more general case of a pair formed by a closed
subspace and the image of a projection.

\begin{prop}\label{bp2}
Let $B$ be a closed subspace
and $P$ be a projection in $E$.
Then
\begin{equation*}
  (I-P)(B) = \ker P \cap (B+\im P) ,
\end{equation*}
and $(I-P)(B)$ is closed in $E$ 
if and only if $B+\im P$ is closed in $E$.
Furthermore, the codimension of $(I-P)(B)$
in $\ker P$ is equal to the codimension
of $B+\im P$ in $E$.
In particular, $(I-P):B\to\ker P$ is a
left-Fredholm operator if and only if
$(B,\im P)$ is a left-Fredholm pair,
and then 
\begin{equation*}
  \ind((I-P):B\to\ker P)
  = \ind (B,\im P) .
\end{equation*}
\end{prop}

\begin{proof}
Let $x\in\ker P$ and suppose that $x=y+Pz$
for some $y\in B$. 
Then $x=(I-P)x=(I-P)y\in(I-P)(B)$.
Conversely, if $x=(I-P)y$ for some $y\in B$,
then $x=y-Py\in B+\im P$.
This shows the first assertion.

If $B+\im P$ is closed in $E$, 
then also $(I-P)(B)=\ker P\cap(B+\im P)$.
Vice versa, suppose that $(I-P)(B)$ is closed
and let $(x_n=y_n+z_n)$ be a sequence in $B+\im P$
converging to $x\in E$.
Then 
\[
  (I-P)y_n = (I-P)x_n \to (I-P)x ,
\]
hence there is a $y\in B$ with $(I-P)y=(I-P)x$,
by assumption.
Hence 
\[
  x = (I-P)y+Px = y+P(x-y)\in B+\im P .
\]
It follows that $(I-P)(B)$ is closed
if and only if $B+\im P$ is closed.

The natural linear map $\ker P\to E/(B+\im P)$
is surjective with kernel 
$\ker P\cap(B+\im P)=(I-P)(B)$.
Hence the codimension of $(I-P)(B)$ in $\ker P$
is equal to the codimension of $B+\im P$ in $E$.
The remaining assertions follow.
\end{proof}

\begin{prop}[Stability]\label{bp4}
Let $P,Q$ be projections in $E$
such that $P-Q$ is compact.
Then $(\im P,\ker Q)$ is a Fredholm pair.

If $B$ is a closed subspace of $E$,
then $(B,\im P)$ is a left-Fredholm pair
if and only if $(B,\im Q)$ is a left-Fredholm pair,
and then 
\[
  \ind(B,\im P)
  = \ind(B,\im Q) + \ind(\ker Q,\im P) .
\]
\end{prop}

\begin{proof}
It is immediate from \lref{trfred}
that $(I-P):\ker Q\to E$ is a left-Fredholm operator.
Applying \pref{bp2} to $B=\ker Q$
we get that $\im P+\ker Q$ is closed in $E$.
The annihilator of $\im P+\ker Q$ in the dual space $E'$
is $\ker P'\cap\im Q'$.
Now $P'-Q'$ is compact, hence $\ker P'\cap\im Q'$
is of finite dimension.
It follows that $(\ker Q,\im P)$ is a Fredholm pair.
We also have
\[
  (I-P)(I-Q) = (I-P)-(I-P)Q = (I-P)+C ,
\]
where $C=(P-Q)Q$ is compact.
Now \pref{bp2} applies.
\end{proof}

\section{An inequality}\label{ineq}

In the proof of \lref{lem:h1restr}, 
we need a special case of the Sobolev inequality, 
cf. Theorem 3.9 in \cite{Ag}.
For the sake of completeness, we give a very simple proof here.

Let $\sigma$ be a complex valued Lipschitz function
on some interval $I\subset\R$.
Then
\begin{equation}\label{ineq2}
  (|\sigma|^2)' = 2\Re(\sigma'\bar\sigma) .
\end{equation}
Hence, if $I=[s,\infty)$ and $\sigma$ has compact support, 
then
\begin{equation}\label{ineq4}
  a |\sigma(s)|^2 
  \le ||\sigma'||_{L^2([s,\infty))}^2 
      + a^2 ||\sigma||_{L^2([s,\infty))}^2 ,
\end{equation}
for any constant $a>0$.
A corresponding estimate holds for bounded intervals:
Let $s<t$ and $\sigma$ be a complex valued Lipschitz function
on $[s,t]$.
Then, for any constant $a>0$,
\begin{equation}\label{ineq6}
  a |\sigma(s) - \sigma(t)|^2 \le
  2 ||\sigma'||_{L^2([s,t])}^2 + 2a^2 ||\sigma||_{L^2([s,t])}^2 .
\end{equation}

\begin{proof}[Proof of \eqref{ineq6}]
By shifting $[s,t]$ if necessary we can assume $s=-t$.
Since even functions are perpendicular to odd functions
in $H^1([-t,t])$,
we can assume that $\sigma$ is odd.
Then the left hand side of the inequality is equal 
to $4a|\sigma(t)|^2$.
Using \eqref{ineq2} and $\sigma(0)=0$,
we derive the asserted estimate.
\end{proof}

\end{appendix}


\end{document}